\pdfoutput=1
\RequirePackage{ifpdf}
\ifpdf 
\documentclass[pdftex]{sigma}
\else
\documentclass{sigma}
\fi

\usepackage{mathrsfs}

\numberwithin{equation}{section}

\newtheorem{Theorem}{Theorem}[section]
\newtheorem{Corollary}[Theorem]{Corollary}
\newtheorem{Lemma}[Theorem]{Lemma}
\newtheorem{Proposition}[Theorem]{Proposition}
 { \theoremstyle{definition}
\newtheorem{Remark}[Theorem]{Remark} }

\newtheorem{Assumption}[Theorem]{Assumption}
\newtheorem*{mainthm}{Main Theorem}
\newtheorem{Claim}[Theorem]{Claim}

\newcommand\ie{\operatorname{ie}}

\newcommand\Ie{{}^{\ie}}
\newcommand\CI{{\mathcal{C}}^{\infty}}
\newcommand\bh{\mathbf{h}}
\newcommand\bn{\mathbf{n}}
\newcommand\bv{\mathbf{v}}
\newcommand{\lrpar}[1]{\left( #1 \right)}
\newcommand{\lrspar}[1]{\left[ #1 \right]}
\newcommand\ang[1]{\langle #1 \rangle}

\newcommand\ev{\operatorname{even}}
\newcommand\Hom{\operatorname{Hom}}
\newcommand\Id{\operatorname{Id}}
\newcommand\pt{\operatorname{pt}}
\newcommand\bbC{\mathbb{C}}
\newcommand\bbR{\mathbb{R}}
\newcommand\bbS{\mathbb{S}}
\newcommand\cO{\mathcal{O}}
\newcommand\cR{\mathcal{R}}
\newcommand\cS{\mathcal{S}}
\newcommand\cV{\mathcal{V}}
\newcommand\sC{\mathscr{C}}

\newcommand \p {\partial}
\newcommand \pa {\partial}
\newcommand \SPB {\mathcal{S}}
\newcommand \Gr {\operatorname{Gr}}
\newcommand \lra {\longrightarrow}
\newcommand\xlra[1]{\xrightarrow{\phantom{x} #1 \phantom{x}}}
\newcommand \wt [1]{\widetilde{#1}}
\newcommand \pdo {\Psi}
\newcommand \Nop {B_{\epsilon}}
\newcommand \YY {\mathcal{Y}}
\newcommand \cl {c}

\DeclareMathOperator \sign {sgn}
\renewcommand{\hat}[1]{\widehat{#1}}
\renewcommand{\bar}[1]{\overline{#1}}

\makeatletter

\newcommand{\Rmnum}[1]{\expandafter\@slowromancap\romannumeral #1@}
\makeatother

\newcommand{\la}{\langle}
\newcommand{\ra}{\rangle}
\newcommand{\lp}{\left(}
\newcommand{\rp}{\right)}
\newcommand{\set}[1]{\left\{ #1 \right\} }
\newcommand \absv [1]{\left \lvert #1 \right \rvert }
\newcommand{\norm}[2][]{\left \| #2 \right \|_{#1} }
\newcommand{\rest}[1]{\big\rvert_{#1}} 

\DeclareMathOperator \ff {f\/f}
\DeclareMathOperator \Tr {Tr}
\DeclareMathOperator \End {End}
\DeclareMathOperator \spec {spec}

\DeclareMathOperator \ext {Ext}
\DeclareMathOperator \diag {diag}

\DeclareMathOperator \spn {span}
\DeclareMathOperator \Diff {Dif\/f}
\DeclareMathOperator \Ran {Ran}
\DeclareMathOperator \supp {supp}
\DeclareMathOperator \ind {Ind}
\DeclareMathOperator \lf {lf}
\DeclareMathOperator \rf {rf}
\newcommand\eps{\varepsilon}
\renewcommand\epsilon{\varepsilon}

\begin{document}


\newcommand{\arXivNumber}{1312.4241}

\renewcommand{\PaperNumber}{089}

\FirstPageHeading

\ShortArticleName{The Index of Dirac Operators on Incomplete Edge Spaces}

\ArticleName{The Index of Dirac Operators\\ on Incomplete Edge Spaces}

\Author{Pierre ALBIN~$^\dag$ and Jesse GELL-REDMAN~$^\ddag$}

\AuthorNameForHeading{P.~Albin and J.~Gell-Redman}

\Address{$^\dag$~University of Illinois, Urbana-Champaign, USA}
\EmailD{\href{mailto:palbin@illinois.edu}{palbin@illinois.edu}}
\URLaddressD{\url{http://www.math.uiuc.edu/~palbin/}}

\Address{$^\ddag$~Department of Mathematics, University of Melbourne, Melbourne, Australia}
\EmailD{\href{mailto:j.gell@unimelb.edu.au}{j.gell@unimelb.edu.au}}
\URLaddressD{\url{http://www.math.jhu.edu/~jgell/}}

\ArticleDates{Received November 02, 2015, in f\/inal form August 30, 2016; Published online September 08, 2016}

\Abstract{We derive a formula for the index of a Dirac operator on a compact, even-dimensional incomplete edge space satisfying a ``geometric Witt condition''. We accomplish this by cutting of\/f to a smooth manifold with boundary, applying the Atiyah--Patodi--Singer index theorem, and taking a limit. We deduce corollaries related to the existence of positive scalar curvature metrics on incomplete edge spaces.}

\Keywords{Atiyah--Singer index theorem; Dirac operators; singular spaces; positive scalar curvature}

\Classification{58G10; 58A35; 58G05}

\vspace{-3mm}
{\small \tableofcontents}

\section{Introduction}

Ever since Cheeger's celebrated study of the spectral invariants of singular spaces \cite{C1979, C1980, C1983} there has been a great deal of research to extend our understanding of geometric analysis from smooth spaces. Index theory in particular has been extended to spaces with isolated conic singularities quite successfully (beyond the papers of Cheeger see, e.g., \cite{Bruning-Seeley:Index, FST99, GLM01, L1997}) and was used by Bismut and Cheeger to establish their families index theorem on manifolds with boundary~\cite{BC1990, BC1990II, BC1991}.

The fact that Bismut and Cheeger used, \cite[Theorem 1.5]{BC1990}, is that for a Dirac operator on a spin space with a conic singularity, the null space of $L^2$ sections naturally corresponds to the null space of the Dirac operator on the manifold with boundary obtained by excising the singularity and imposing the `Atiyah--Patodi--Singer boundary condition'~\cite{APSI}, provided an induced Dirac operator on the link has no kernel. (Here `Dirac operator' refers to the `classical' Dirac operator $\Tr \cl \circ \nabla$ where $\nabla$ is the Levi-Cevita connection on the spin
bundle and $\cl$ is Clif\/ford multiplication.) This is not true for more general f\/irst order operators; indeed, even for twisted Dirac operators there is no such equivalence. On the other hand, Cheeger points out in~\cite{C1979} that one can leverage the related APS boundary value problems to prove the Gauss--Bonnet and Hirzebruch signature theorems for spaces with conic singularities.

In this paper we consider the Dirac operator on a spin space with non-isolated conic singularities, also known as an `incomplete edge space', and the Dirac operator on the manifold with boundary obtained by excising a tubular neighborhood of the singularity and imposing the Atiyah--Patodi--Singer boundary condition. Although the relation between the domains of these two Dirac operators is much more complicated than in the case of isolated conic singularities, we show that under a~``geometric Witt assumption'' analogous to that used by Bismut--Cheeger, the indices of these operators
coincide. Thus we obtain a~formula for the index of the Dirac operator on the singular space as the `adiabatic limit' of the index of the Dirac operator with Atiyah--Patodi--Singer boundary conditions.

 \begin{figure}[t] \centering
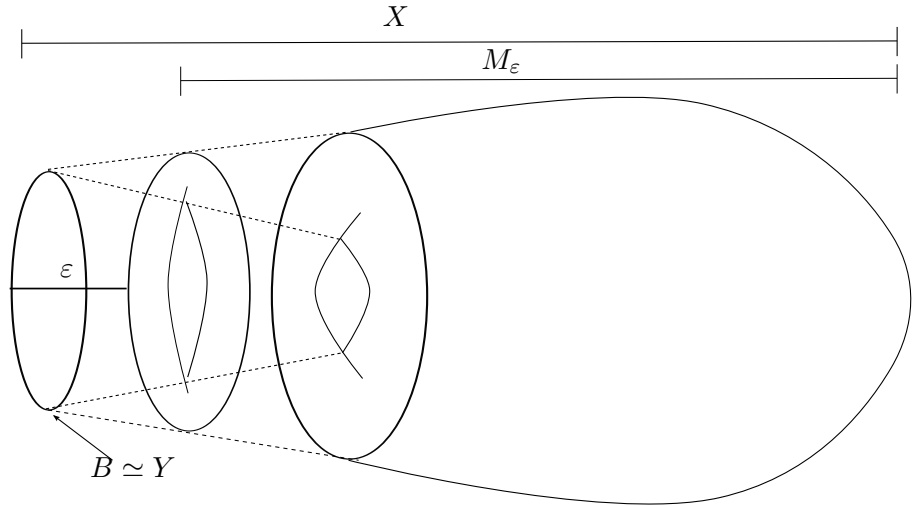\vspace{-0.5mm}
 \caption{The singular space $X$ obtained by collapsing the f\/ibers of the boundary f\/ibration of~$M$. The spaces $M_{\epsilon}$ play a central role in our proofs.}\label{fig:geometric}\vspace{-1.5mm}
 \end{figure}

An incomplete edge (ie) space is a stratif\/ied space $X$ with a single singular stratum~$Y$. In keeping with Melrose's paradigm for analysis on singular spaces (see, e.g., \cite{Melrose1992, Melrose1996}) we resolve~$X$ by `blowing-up' $Y$ and obtain a smooth manifold with boundary~$M$, whose boundary is the total space of a f\/ibration of smooth manifolds with typical f\/iber~$Z$,
\begin{gather*}
	Z \hookrightarrow \p M \xlra{\phi} Y.
\end{gather*}
A (product-type) incomplete edge metric is a metric that, in a collar neighborhood of the boundary, takes the form
\begin{gather}\label{eq:exactmetric}
	g = dx^2 + x^2 g_Z + \phi^*g_Y
\end{gather}
with $x$ a def\/ining function for $\p M$, $g_Y$ a metric on $Y$ and $g_Z$ a family of two-tensors that restrict to a metric on each f\/iber of~$\phi$. Thus we see that metrically the f\/ibers of the boundary f\/ibration are collapsed, as they are in~$X$.

We also replace the cotangent bundle of $M$ by a bundle adapted to the geometry, the `incomplete edge cotangent bundle' $T_{\ie}^*M$, see~\eqref{eqn:incomplete-edge-tangent-def} below. This bundle is locally spanned by forms like~$dx$, $x\,dz$, and $dy$, and the main dif\/ference with the usual cotangent bundle is that the form~$x\, dz$ is a non-vanishing section of $T_{\ie}^*M$ all the way to~$\pa M$.

We assume that $M^\circ$ is spin and denote a spin bundle on $M^\circ$ by $\SPB \lra M$ and the associated Dirac operator by $\eth$. The operator $\eth$ does not induce an operator on the boundary in the usual sense, due to the degeneracy of the metric there, but we do have
\begin{gather*}
	x\eth\rest{\p M} = c(dx)\big(\tfrac12\dim Z + \eth_Z\big),
\end{gather*}
where $\eth_Z$ is a family of operators on the f\/ibers of $\p M$. It turns out, as in the conic case mentioned above, and the analogous study of the signature operator in~\cite{ALMPI, ALMPII}, that much of the functional analytic behavior of $\eth$ is tied to that of $\eth_Z$. Indeed, in Section~\ref{sec:proofofessential} below we prove the following theorem.

\begin{Theorem}\label{thm:essential}
Assume that $\eth$ is a Dirac operator on a compact, spin incomplete edge space $(M,g)$, satisfying the ``geometric Witt-assumption''
\begin{gather}\label{eq:GeoWittAss}
	\operatorname{Spec} ( \eth_Z) \cap (-1/2,1/2) = \varnothing.
\end{gather}
Then the unbounded operator $\eth$ on $L^{2}(M; \SPB)$ with core domain $C^{\infty}_{c}(M; \SPB)$ $($sections supported in $M \setminus
\p M)$ is essentially self-adjoint. Moreover, letting $\mathcal{D}$ denote the domain of this self-adjoint extension, the map
\begin{gather*}
 \eth \colon \ \mathcal{D} \lra L^{2}(M; \SPB)
\end{gather*}
is Fredholm.
\end{Theorem}

When $M$ is even-dimensional, the spin bundle admits the standard $\mathbb{Z}/2\mathbb{Z}$ grading into even and odd spinors
\begin{gather*}
	\SPB = \SPB^{+} \oplus \SPB^{-},
\end{gather*}
and thus we have the chirality spaces $\mathcal{D}^{\pm} = \mathcal{D} \cap L^{2}(M; \SPB^{\pm})$ and the restriction of the Dirac operator satisf\/ies
\begin{gather*}
\eth \colon \ \mathcal{D}^{+} \lra L^{2}(M; \SPB^{-}).
\end{gather*}
This map is Fredholm and our main result is an explicit formula for its index.

The metric \eqref{eq:exactmetric} naturally def\/ines a bundle metric on $T_{\ie}^*M$, non-degenerate at $\pa M$, and the Levi-Civita
connection of $g$ naturally def\/ines a connection $\nabla$ on~$T_{\ie}^*M$. Our index formula involves the transgression of a~characteristic class between two related connections. The restriction of~$T_{\ie}M$ to~$\pa M$ can be identif\/ied with $N_M\pa M \oplus T\pa M/Y \oplus \phi^*TY$. Let
\begin{gather*}
	\bn\colon \ T_{\ie}M\rest{\pa M} \lra N_M\pa M, \qquad \bv\colon \ T_{\ie}M\rest{\pa M} \lra T\pa M/Y,	
\end{gather*}
be the orthogonal projections onto the normal bundle of $\pa M$ in $M$, $N_M\pa M = \ang{\pa_x}$, and the vertical bundle of $\phi$, respectively, and let $\bv_+ = \bn\oplus \bv$. Both
\begin{gather}\label{eq:TwoConn}
	\nabla^{v_+} = \bv_+ \circ \nabla\rest{\pa M} \circ \bv_+, \qquad \text{and} \qquad
	\nabla^{\pt} = \bn \circ \nabla\rest{\pa M} \circ \bn \oplus \bv \circ j_0^*\nabla \circ \bv,
\end{gather}
where $\pt$ stands for `product', are connections on $N_M\pa M \oplus T\pa M/Y \lra \pa M$.

\begin{mainthm}
Let $X$ be a compact, even-dimensional stratified space with a single singular stratum endowed with an incomplete edge metric $g$ and let $M$ be its resolution. If $\eth$ is a Dirac operator associated to a spin bundle $\SPB\lra M$ and $\eth$ satisfies the geometric Witt condition~\eqref{eq:GeoWittAss}, then
\begin{gather} \label{eq:finalindexformula}
	\ind\big(\eth \colon \mathcal{D}^{+} \lra L^{2}(M; \SPB^{-})\big)
	= \int_{M}\! \hat{A}(M) + \int_Y\! \hat A(Y)\lrpar{{-}\frac 12\hat\eta(\eth_Z) + \int_{Z}\! T\hat A(\nabla^{v_+},\nabla^{\pt}) }\!,\!\!\!
\end{gather}
where $\hat A$ denotes the $\hat A$-genus, $T\hat A(\nabla^{v_+}, \nabla^{\pt})$ denotes the transgression form of the~$\hat A$ genus associated to the connections~\eqref{eq:TwoConn}, and $\hat\eta$ the $\eta$-form of Bismut--Cheeger~{\rm \cite{BC1989}}.
\end{mainthm}

The simplest setting of incomplete edge spaces occurs when $Z$ is a~sphere, as then $X$ can be given a smooth structure and the singularity at $Y$ is entirely in the metric. Atiyah and LeBrun have recently studied the case where $Z=\bbS^1$ and $X$ is four-dimensional, so that $Y$ is an embedded surface, and the metric $g$ asymptotically has the form
\begin{gather*}
	dx^2 + x^2 \beta^{2}d\theta^{2} + \phi^*g_Y.
\end{gather*}
The cone angle $2\pi\beta$ is assumed to be constant along $Y$.
In~\cite{AL2013} they f\/ind formulas for the signature and the Euler characteristic of~$X$ in terms of the curvature of this incomplete edge metric. Kronheimer and Mrowka also study invariants of such spaces, in particular working out the dimension of the moduli space of naturally associated f\/lat connections \cite{MR1241873,MR1308489}, and Lock and Viaclovsky~\cite{LV2013} compute the index of the `anti-self-dual deformation complex'. Using work of Dai~\cite{Dai:adiabatic} and Dai--Zhang~\cite{DZ1995}, in Theorem~\ref{thm:signature} below, we recover the formula for the signature in~\cite{AL2013}, and moreover, we show that our formula for the index of the Dirac operator~\eqref{eq:finalindexformula} simplif\/ies substantially in the case $Z = {\mathbb S}^1$ and $\dim X = 4$. In the context of general Witt spaces, on the other hand, work on the signature operator and $L$-class in the incomplete edge setting includes~\cite{Bruning-signature} and~\cite{C1997}.

In studying the Dirac operator on incomplete edges with sphere f\/ibers, it is natural to assume, and we do so below, that the manifold~$X$ itself has a spin structure. This induces a spin structure on the interior of~$M$ (which is dif\/feomorphic to $X$ minus the singular locus), which
extends to~$M$. Then the induced spin structure on the f\/iber $Z \simeq \mathbb{S}^f$ is the spin structure induced from thinking of $\mathbb{S}^f$ as the boundary of the ball $\mathbb{B}^{f + 1}$ in Euclidean space. In particular when $f = 1$, so the f\/iber~$Z$ is a circle, the induced spin structure on $\mathbb{S}^1 = \mathbb{R}/2\pi\mathbb{Z}$ with the round metric $d\theta^2$, induce the Dirac operator $\eth_\theta = -i \p_\theta + 1/2$ whose spectrum satisf\/ies $\spec(\eth_\theta) = \{ 1/2 + \mathbb{Z} \}$. (Note that the spectrum of the Dirac operator depends on the choice of spin structure, so it is important that we have the bounding spin structure on the circle.) As we show in the proof of the following corollary, this implies that for cone angles $2 \pi \beta \le 2\pi$, the geometric Witt assumption is satisf\/ied and we arrive at a~simpler index formula in this case.

\begin{Corollary}\label{thm:cor4dim}
If $\eth$ is a Dirac operator on a smooth, compact four-dimensional spin mani\-fold~$X$, associated to an incomplete edge metric with constant cone angle $2\pi\beta\leq 2\pi$ along an embedded surface $Y$, then $\eth$ is essentially self-adjoint and its index is given by
\begin{gather}\label{eq:finalindexformuladim4}
	\ind\big(\eth \colon \mathcal{D}^{+} \lra L^{2}(M; \SPB^{-})\big)
	= - \frac{1}{24} \int_{M} p_{1}(M) + \frac{1}{24} \big(\beta^{2} - 1\big) [Y]^{2},
\end{gather}
where $[Y]^{2}$ is the self-intersection number of $Y$ in $X$.
\end{Corollary}

The formulas in \eqref{eq:finalindexformula} and \eqref{eq:finalindexformuladim4}, and indeed our proof, are obtained by taking the limit of the index formula for the Dirac operators on the manifolds with boundary
\begin{gather*}
	M_{\epsilon} = \{ x \geq \epsilon \},
\end{gather*}
so in particular the contribution from the singular stratum $Y$ is the adiabatic limit \cite{ BC1989, Dai:adiabatic, W1985} of the $\eta$-invariant
from the celebrated classical theorem of Atiyah, Patodi, and Singer \cite{APSI}, which we review in Section~\ref{sec:maintheorem}. It is important to note that the analogous statement for a general twisted Dirac operator is {\em false} and the general index formula requires an extra contribution from the singularity. We will return to this in a subsequent publication. We also point out to the reader that there exist other derivations of index formulas on manifolds with structured ends in which the computation is reduced to taking a limit in the Atiyah--Patodi--Singer index formula; see for example~\cite{Bruning-signature} where the author discovers a self-adjointness criterion and proves an $L^2$-signature theorem on incomplete edge
spaces, or~\cite{Leichtnam-Mazzeo-Piazza} where the authors prove an index formula for twisted Dirac operators on spin manifolds with f\/ibered boundary metrics, a complete Riemannian manifold with a~structured end that is a f\/iber bundle over an asymptotically conical (big end of a cone) manifold.

One very interesting aspect of the spin Dirac operator is its close relation to the existence of positive scalar curvature metrics. Most directly, the Lichnerowicz formula shows that the index of the Dirac operator is an obstruction to the existence of such a metric. This is still true among metrics with incomplete edge singularities. Analogously to the results of Chou for conic singularities~\cite{Ch1985} we prove the following theorem in Section~\ref{sec:posscal}.

\begin{Theorem} \label{thm:PosScal}
Let $(M,g)$ be a spin incomplete edge space. The geometric Witt assump\-tion'~\eqref{eq:GeoWittAss} holds if either:
\begin{enumerate}\itemsep=0pt
\item[$1.$] $\dim Z \ge 2$ and the scalar curvature of $g$ is non-negative in a neighborhood of~$\p M$.
\item[$2.$] $\dim Z = 1$, the spin structure on $M$ is the lift of a spin structure on $X$, and the cone angle satisfies $2 \pi \beta \le 2 \pi$.
\end{enumerate}

If the geometric Witt assumption holds and in addition the scalar curvature of $g$ is non-negative on all of $M$, and positive somewhere, then $\ind(\eth)=0$.
\end{Theorem}

Now let us indicate in more detail how these theorems are proved. For convenience we work throughout with a product-type incomplete edge metric as described above, but removing this assumption would only result in slightly more intricate computations below. The proof of Theorem~\ref{thm:essential} follows the arguments employed in \cite{ALMPI, ALMPII} to prove the analogous result for the signature operator. Thus we start with the two canonical closed extensions of~$\eth$ from $C^{\infty}_{c}(M)$, the sections with support a compact set in the interior $M \setminus \p M$, namely
\begin{gather}
 \mathcal{D}_{\max} := \set{ u \in L^{2}(M ; \SPB) \colon \eth u \in L^{2}(M; \SPB)}, \nonumber\\
 \mathcal{D}_{\min} := \set{ u \in \mathcal{D}_{\max} \colon \exists\, u_{k} \in C^{\infty}_{c}(M) \mbox{ with } u_{k} \to u,\, \eth u_{k} \to \eth
 u \mbox{ as } k \to \infty}, \label{eq:maxandmin}
\end{gather}
where the convergence in the second
def\/inition is in $L^{2}(M; \SPB)$, and we show that under assumption~\eqref{eq:GeoWittAss}, these domains coincide
\begin{gather}\label{eq:minismax}
	\mathcal{D}_{\min} = \mathcal{D}_{\max} = \mathcal{D}.
\end{gather}
Since $\eth$ is a symmetric operator, this shows that it is essentially self-adjoint.

One dif\/ference between the case of isolated conic singularities ($\dim Y = 0$) and the general incomplete edge case is that in the former, even if Assumption \eqref{eq:GeoWittAss} does not hold, $\mathcal D_{\max}/\mathcal D_{\min}$ is a f\/inite-dimensional space. In contrast, when $\dim Y>0$, this space is generally inf\/inite-dimensional.

We prove \eqref{eq:minismax} by constructing a parametrix $\overline{Q}$ for $\eth$ in Section \ref{sec:parametrix}. From the mapping properties of $\overline{Q}$, we deduce both that $\eth$ is essentially self-adjoint, and that it is a Fredholm operator from the domain of its unique self-adjoint extension to $L^{2}$. The relationship between the mapping properties of $\overline{Q}$ and the stated conclusions can be seen largely through~\eqref{eq:extravanishing} below, which states that the maximal domain has `extra' vanishing, i.e., sections in~$\mathcal{D}_{\max}$ lie in weighted spaces $x^{\delta}L^{2}$ with weight $\delta$ higher than generically expected. This shows that the inclusion of the domain into $L^2$ is a~compact operator, which in particular gives that the kernel of $\eth$ on the maximal domain is f\/inite-dimensional.

\begin{figure}\centering
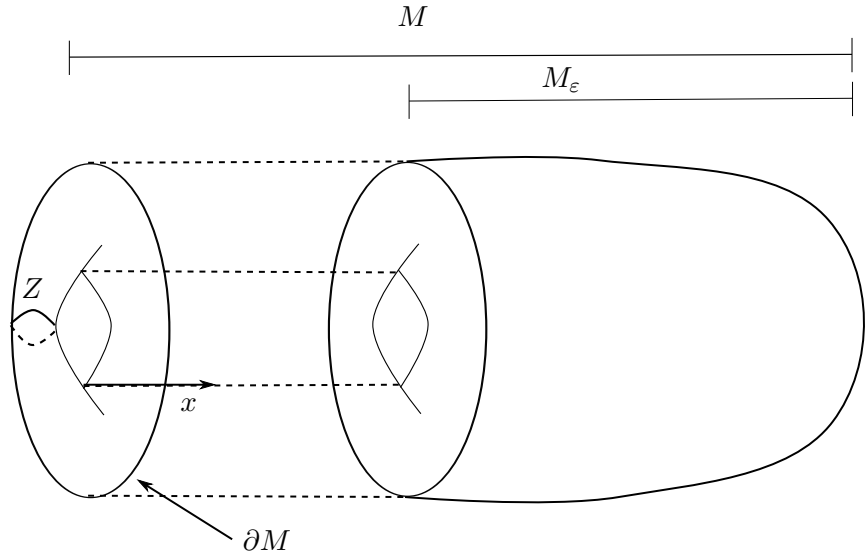
\caption{$M$ as a smooth manifold with boundary whose boundary $\p M$ is a f\/iber bundle. Here $x$ is a~boundary def\/ining function and the space $M_{\epsilon}$ are given by $M_{\epsilon} = \set{ x \ge \epsilon}$.}\label{fig:smooth}
\end{figure}

Once this is established we give a precise description of the Schwartz kernel of the generalized inverse $Q$ of $\eth$ using the technology of \cite{Ma1991, MV2012}. In Section \ref{sec:boundaryvalues}, we use $Q$ and standard methods from layer potentials to construct a family of pseudodif\/ferential projectors
\begin{gather*}
	\mathcal E_{\epsilon} \in \Psi^{0}(\p M_{\epsilon}; \SPB)
\end{gather*}
such that
\begin{gather*}
\eth\rest{M_{\epsilon}} \text{ with domain } \mathcal D_{\epsilon} = \set{ u \in H^{1}(M_{\epsilon}; \SPB)\colon ({\operatorname{Id}} - \mathcal{E}_{\epsilon})(u \rvert_{\p M_{\epsilon}}) = 0 }.
\end{gather*}
is Fredholm and has the same index as $(\eth, \mathcal D)$. This domain is constructed so that the boundary values coincide with boundary values of
`$\eth$-harmonic' $L^{2}$-sections over the excised neighborhood of the singularity, $M - M_{\epsilon}$.

To compute this index, we consider the operators
\begin{gather*}
	\eth\rest{M_{\epsilon}} \text{ with domain }
	\mathcal D_{{\rm APS}, \epsilon} = \set{ u \in H^{1}(M_{\epsilon}; \SPB)\colon ({\operatorname{Id}} -
 \pi_{{\rm APS}, \epsilon})(u \rvert_{\p M_{\epsilon}}) = 0 }.
\end{gather*}
where $\pi_{{\rm APS},\epsilon}$, is the projection onto the positive spectrum of $\eth\rest{\p M_{\epsilon}}$. From \cite{APSI} we know that these are Fredholm operators and
\begin{gather}\label{eq:APSIndex}
	\ind\big(\eth, \mathcal D^{+}_{{\rm APS},\epsilon} \lra L^{2}(M; \SPB^{-})\big) = \int_M \hat A(M) - \tfrac12\eta(\eth_{\p M_{\epsilon}}) + \mathcal L_{\epsilon},
\end{gather}
where $\mathcal L_{\epsilon}$ is a local integral over $\p M_{\epsilon}$ compensating for the fact that the metric is not of product-type at $\p M_{\epsilon}$. These domains depend fundamentally on $\epsilon$. Not only does $\mathcal D_{{\rm APS},\epsilon}$ vary as $\epsilon \to 0$, it does not limit to a f\/ixed subspace of $L^{2}(\p M)$ with any natural metric. (More precisely, the boundary value projectors $\pi_{{\rm APS}, \epsilon}$ which def\/ine the boundary condition do not converge in norm.)

Through a semiclassical analysis, which we carry out using the adiabatic calculus of Mazzeo--Melrose \cite{MaMe1990}, we show that the projections $\mathcal{E}_{\epsilon}$ and $\pi_{{\rm APS}, \epsilon}$ are homotopic for small enough~$\epsilon$, with a homotopy through operators with the same principal symbol. The adiabatic calculus technology boils this down to an explicit analysis of modif\/ied Bessel functions, which we carry out in the appendix. Then we can appeal to arguments from Booss--Bavnbek--Wojciechowski~\cite{BB1993} to see that the two boundary value problems have the same index. Having shown that the index of $(\eth, \mathcal D)$ is equal to the adiabatic limit of the index formula~\eqref{eq:APSIndex}, the Main Theorem follows as shown in Section~\ref{sec:maintheorem}.

\looseness=-1
Incomplete edge dif\/ferential operators, such as the Dirac operator for on an incomplete edge spin manifold, are closely related with (complete) edge dif\/ferential operators, which include the natural geometric PDEs associated with edge metrics. Indeed, our analysis of~$\eth$ here relies partly on pre-existing analysis of edge dif\/ferential operators~\cite{MV2013}, and in the case of the Dirac operator of a given spin structure (and \textit{not} in the general case of twisted Dirac operators) the complete and incomplete edge cases are related by conformal invariance properties. As we mention in Remark~\ref{thm:edge-theorem}, if the geometric Witt assumption is strengthened to exclude f\/iber spectrum from the \textit{closed} interval $[-1/2, 1/2]$, then our results imply an index formula for Dirac operators on edge manifolds (which under this strengthened assumption are Fredholm on their natural domain.)

Though we do not attempt to give a complete overview here, there is a~large body of work on index theory of non-compact Riemannian manifolds which is closely related to our work here. In particular, the geometric microlocal perspective, in which one compactif\/ies complete manifolds and uses radial blowups to resolve a geometrically natural Lie algebra a vector f\/ields and def\/ine a calculus of pseudodif\/ferential operators, has led to index theorems on complete Riemannian manifolds with structured ends \cite{tapsit}. In particular there is a wealth of work (though not many index theorems) on (complete) edge dif\/ferential operators \eqref{eq:ieoperator'}, see~\cite{Ma1991}. There is, moreover, work on dif\/ferential operators on singular spaces and their Fredholm and index theory from the perspective of groupoids in particular for edge dif\/ferential operators \cite{ALN2007, DLN2009,N2015}. See also~\cite{MN1996}. The f\/irst author proved a renormalized index formula for `Dirac-type' operators, specif\/ically operators arising as $\Tr \cl \circ \nabla$ for~$\nabla$ a Clif\/ford-connection on a Hermitian Clif\/ford bundle, on \textit{complete} edge manifolds. The complete setting involves the subtlety that such operators, though essentially self-adjoint, are not guaranteed to be Fredholm on their natural domains. As far as the authors are aware, ours is the f\/irst index formula for the classical Dirac operator in either the incomplete edge or (complete) edge context.

\section{Connection and Dirac operator}\label{sec:Connectionandoperator}

Let $(M, g)$ be an incomplete edge space which is spin, $\SPB \lra M$ the spinor bundle for a f\/ixed spin structure with connection $\nabla$, and
let $\eth$ be the corresponding Dirac operator. Given an orthonormal frame $e_{i}$ of the tangent bundle of~$M$, the Dirac operator satisf\/ies
\begin{gather*}
 \eth = \sum_{i} c(e_{i}) \nabla_{e_{i}},
\end{gather*}
where $c(v)$ denotes Clif\/ford multiplication by the vector $v$. See \cite{LM1989, Roe} for background on spinor bundles and Dirac operators. The main goal of this section is to prove Lemma~\ref{thm:diracwithgoodconnection} below, where we produce a tractable form of the Dirac operator on a~collar neighborhood of the bounda\-ry~$\p M$, or equivalently of the singular stratum~$Y \subset X$.

\subsection{Incomplete edge metrics and their connections}\label{sec:Conn}

Let $M$ be the interior of a compact manifold with boundary. Assume that $\pa M=N$ participates in a f\/iber bundle
\begin{gather*}
	Z \hookrightarrow N \xlra\phi Y.
\end{gather*}
Let $X$ be the singular space obtained from $M$ by collapsing the f\/ibers of the f\/ibration $\phi$. If we want to understand the dif\/ferential forms on $X$ while working on $M$, it is natural to restrict our attention to
\begin{gather}\label{eq:IeOneForms}
	\{ \omega \in \CI(M;T^*M)\colon i_N^*\omega \in \phi^*\CI(Y;T^*Y) \}.
\end{gather}
Following Melrose's approach to analysis on singular spaces \cite{Melrose:Survey} let
\begin{gather}\label{eqn:incomplete-edge-tangent-def}
	T_{\ie}^*M \lra M
\end{gather}
be the vector bundle whose space of sections is \eqref{eq:IeOneForms}. We call $T_{\ie}^*M$ the `incomplete edge cotangent bundle', and its dual bundle $T_{\ie}M$, the `incomplete edge tangent bundle'. (Note that $T_{\ie}M$ is simply a rescaled bundle of the (complete) `edge tangent bundle' of Mazzeo \cite{Ma1991}.) The incomplete edge tangent bundle $T_{\ie}M$ is canonically isomorphic to $TM$ over the interior $M^\circ$, but its extension to $M$ is not canonically isomorphic to $T M$ (though they are of course isomorphic bundles) as we discuss below.
\begin{Remark}
The bundle $T^*_{\ie}M$ is the natural space for def\/ining uniform ellipticity of the natural operators, such as the Hodge--Laplacian or the
Dirac operator, on incomplete edge spaces. Indeed, for example, the Laplacian on functions $\Delta_g$ for an incomplete edge metric~$g$, has principal symbol $\sigma(\Delta_g)$ which acts on covectors $\xi \in T^*M$ by $\sigma(\Delta_g)(\xi) = \|\xi\|^2_g$. This extends up to $\p M$ as a function $\sigma(\Delta_g) \colon T^*_{\ie}M \lra \mathbb{R}$ and is non-vanishing away from the zero section (hence uniformly elliptic). It does not extend smoothly as a function on $T^*M$.
\end{Remark}

Let $x$ be a boundary def\/ining function (bdf) on $M$, meaning a smooth non-negative function $x\in \CI( M;[0,\infty))$ such that $\{x=0\} = N$ and $|dx|$ has no zeroes on $N$. Near a point $p \in \p M$, we will typically work in local coordinates on $M$ written
\begin{gather}\label{eq:coordinates}
	x, \quad y, \quad z,
\end{gather}
where $x$ is the bdf above, $y$ are coordinates along $Y$ on a neighborhood of $\phi(p)$ and $z$ are coordinates along $Z$. In particular, we choose a local trivialization of the f\/ibration~$\p M$ on an open set $V \ni \phi(p)$, giving a local dif\/feomorphism $\phi^{-1}(V) \simeq V \times Z$. The~$y$ and~$z$ are local coordinates on the left and right factors, respectively. In local coordinates the sections of $T_{\ie}^*M$ are spanned by
\begin{gather*}
	dx, \quad x\,dz, \quad dy.
\end{gather*}
The crucial fact is that $x \,dz$ vanishes at $N$ as a section of $T^* M$, but it does not vanish at $N$ as a section of $T_{\ie}^*M$ because the `$x$' is here part of the basis element and not a coef\/f\/icient. There is an obvious map from $T_{\ie}^*M$ to $T^*M$ which takes a dif\/ferential form expressed as a linear combination of the above basis forms to the exact same form in $T^*M$; this is an isomorphism on the interior $M^{\circ}$ but takes forms $x\,dz$ to zero over the boundary. Similarly, in local coordinates the sections of $T_{\ie}M$ are spanned by
\begin{gather*}
\pa_x, \quad \tfrac1x\pa_z, \quad \pa_y,
\end{gather*}
and, in contrast to $T M$, the vector f\/ield $\tfrac1x\pa_z$ is def\/ined at $N$ as a section of $T_{\ie}M$.

Next consider a metric on $M$ that ref\/lects the collapse of the f\/ibers of $\phi$. Let $\sC$ be a collar neighborhood of $N$ in $M$ compatible with $x$, $\sC \cong [0,1]_x \times N$.

A product-type incomplete edge metricis a Riemannian metric on $M$ that on $\sC$ (i.e., for some boundary def\/ining function $x$) has the form
\begin{gather}\label{eq:incompleteedgemetric}
	g_{\ie} = dx^2 + x^2 g_Z + \phi^*g_Y,
\end{gather}
where $g_Z + \phi^*g_Y$ is a submersion metric for $\phi$ independent of $x$. Note that this metric naturally induces a bundle metric on $T_{\ie}M$
with the advantage that it extends non-degenerately up to~$\p M$. We will consider this as a metric on~$T_{\ie}M$ from now on. (A~general incomplete edge metric is simply a bundle metric on $T_{\ie}M \lra M$.)

An exact $\ie$-metric induces an orthogonal splitting
\begin{gather*}
	T\sC = \ang{\pa_x} \oplus TN/Y \oplus \phi^*TY
\end{gather*}
def\/ined by noting that $\ang{\pa_x} \oplus TN/Y$ is the kernel of the $\phi_* \colon T\sC \lra TY$ and thus the perpendicular space is isomorphic to $\phi^*TY$

To describe the asymptotics of the Levi-Civita connection of $g_{\ie}$, let us start by recalling the behavior of the Levi-Civita connection of a~submersion metric. Endow $N = \pa M$ with a~sub\-mersion metric of the form $g_N = \phi^*g_Y + g_Z$. Given a~vector f\/ield $U$ on $Y$, let us denote its horizontal lift to~$N$ by~$\wt U$. Also let us denote the projections onto each summand by
\begin{gather*}
	\bh\colon \ TN \lra \phi^*TY, \qquad \bv\colon \ TN \lra TN/Y.
\end{gather*}
The connection $\nabla^N$ dif\/fers from the connections $\nabla^Y$ on the base and the connections $\nabla^{N/Y}$ on the f\/ibers through two tensors. The second fundamental form of the f\/ibers is def\/ined by
\begin{gather*}
	\cS^{\phi}\colon \ TN/Y \times TN/Y \lra \phi^*TY, \qquad \cS^{\phi}(V_1, V_2) = \bh\big(\nabla^{N/Y}_{V_1}V_2\big)
\end{gather*}
and the curvature of the f\/ibration is def\/ined by
\begin{gather*}
	\cR^{\phi}\colon \ \phi^*TY \times \phi^*TY \lra TN/Y, \qquad \cR^{\phi}\big(\wt U_1, \wt U_2\big) = \bv\big(\big[\wt U_1, \wt U_2\big]\big).
\end{gather*}
The behavior of the Levi-Civita connection (cf.~\cite[Proposition 13]{HHM2004}) is then summed up in the table:
\begin{gather*}
\begin{tabular}{|c||c|c|} \hline
$g_{N}\big(\nabla^{N}_{W_1} W_2, W_3\big)$ & $V_0$ & $\wt U_0$\tsep{3pt}\bsep{3pt} \\ \hline\hline
$\nabla^{N}_{V_1}V_2$ &
	$g_{N/Y}\big(\nabla^{N/Y}_{ V_1} V_2, V_0\big)$ &
	$\phi^*g_Y\big(\cS^{\phi}( V_1, V_2), \wt U_0\big)$ \tsep{3pt}\bsep{3pt}\\ \hline
$\nabla^{N}_{\wt U} V$ &
	$g_{N/Y}\big([\wt U, V], V_0\big) - \phi^*g_Y\big(\cS^{\phi}( V, V_0), \wt U\big)$ &
	$-\frac12g_{N/Y}\big(\cR^{\phi}\big(\wt U, \wt U_0\big), V\big)$ \tsep{3pt}\bsep{3pt}\\ \hline
$\nabla^{N}_{ V} \wt U$ &
	$-\phi^*g_Y\big(\cS^{\phi}(V, V_0),\wt U\big)$ &
	$\frac12g_{N/Y}\big(\cR^{\phi}\big(\wt U, \wt U_0\big), V\big)$ \tsep{3pt}\bsep{3pt}\\ \hline
$\nabla^{N}_{\wt U_1}\wt U_2$ &
	$\frac12g_{N/Y}\big(\cR^{\phi}\big(\wt U_1, \wt U_2\big), V_0\big)$ &
	$g_Y\big( \nabla^Y_{U_1} U_2, U_0\big)$ \tsep{3pt}\bsep{3pt}\\ \hline
\end{tabular}
\end{gather*}

We want a similar description of the Levi-Civita connection of an incomplete edge metric. The splitting of the tangent bundle of $\sC$ induces a splitting
\begin{gather}\label{eq:SplittingIE}
	T_{\ie}\sC = \ang{\pa_x} \oplus \tfrac1xTN/Y \oplus \phi^*TY,
\end{gather}
in terms of which a convenient choice of vector f\/ields is
\begin{gather*}
	\pa_x, \quad \tfrac1xV, \quad \wt U,
\end{gather*}
where $V$ denotes a vertical vector f\/ield at $\{x=0\}$ extended trivially to $\sC$ and $\wt U$ denotes a~vector f\/ield on $Y$, lifted to $\pa M$ and then extended trivially to $\sC$. Note that, with respect to $g_{\ie}$, these three types of vector f\/ields are orthogonal, and that their commutators satisfy
\begin{gather*}
	\big[ \pa_x, \tfrac1xV \big] = -\tfrac1{x^2}V \in x^{-1}\CI\big(\sC, \tfrac1xT\pa M/Y\big), \qquad
	\big[ \pa_x, \wt U \big] = 0, \\
	\big[ \tfrac1xV_1, \tfrac1xV_2\big] = \tfrac1{x^2} [V_1, V_2] \in x^{-1}\CI\big(\sC, \tfrac1xT\pa M/Y\big), \\
	\big[\tfrac1xV, \wt U\big] = \tfrac1x\big[ V, \wt U\big] \in \CI\big(\sC, \tfrac1xT\pa M/Y\big), \\
	\big[ \wt U_1, \wt U_2\big] \in x\CI\big(\sC, \tfrac1xT\pa M/Y\big) + \CI(\sC, \phi^*TY).
\end{gather*}
The Levi-Civita connection $\nabla$ for $g_{\ie}$ satisf\/ies the Koszul formula
\begin{gather*}
	2g_{\ie}(\nabla_{W_0}W_1, W_2) = W_0 g_{\ie}(W_1, W_2) + W_1 g_{\ie}(W_0, W_2) - W_2 g_{\ie}(W_0, W_1) \\
\hphantom{2g_{\ie}(\nabla_{W_0}W_1, W_2) =}{}	+ g_{\ie}(\lrspar{W_0, W_1}, W_2) - g_{\ie}(\lrspar{W_0, W_2}, W_1) - g_{\ie}(\lrspar{ W_1, W_2}, W_0),
\end{gather*}
and it is easy to see that the expression on the right hand side is smooth on all of $M$ for any any smooth vector f\/ield $W_0 \in C^\infty(M; TM)$ and $W_1, W_2 \in C^\infty(M; T_{\ie}M)$. In fact we will now describe the action of $\nabla$ in relation to the splitting of $T\sC$ above.

If $W_0 \in \{ \pa_x, V, \wt U\}$ and $W_1, W_2 \in \{\pa_x, \tfrac1x V, \wt U \}$ then we f\/ind
\begin{gather*}
	g_{\ie}(\nabla_{W_0}W_1, W_2) =0 \qquad \text{if}\quad \pa_x \in \{ W_0, W_1, W_2 \} \\
	\text{except for} \quad 	g_{\ie}\big(\nabla_{V_1}\pa_x, \tfrac1x V_2\big) = - g_{\ie}\big(\nabla_{V_1}\tfrac1xV_2, \pa_x\big) = g_{Z}(V_1, V_2),
\end{gather*}
and otherwise
\begin{gather*}
\begin{tabular}{|c||c|c|} \hline
$g_{\ie}(\nabla_{W_1} W_2, W_3)$ & $\tfrac1xV_0$ & $\wt U_0$ \tsep{3pt}\bsep{3pt}\\ \hline\hline
$\nabla_{V_1}\tfrac1x V_2$ &
	$g_{N/Y}\big(\nabla^{N/Y}_{ V_1} V_2, V_0\big)$ &
	$x\phi^*g_Y\big(\cS^{\phi}( V_1, V_2), \wt U_0\big)$ \tsep{3pt}\bsep{3pt}\\ \hline
$\nabla_{\wt U} \tfrac1xV$ &
	$g_{N/Y}\big(\big[\wt U, V\big], V_0\big) - \phi^*g_Y\big(\cS^{\phi}( V, V_0), \wt U\big)$ &
	$-\frac x2g_{N/Y}\big(\cR^{\phi}\big(\wt U, \wt U_0\big), V\big)$ \tsep{3pt}\bsep{3pt}\\ \hline
$\nabla_{ V} \wt U$ &
	$-x\phi^*g_Y\big(\cS^{\phi}(V, V_0),\wt U\big)$ &
	$\frac{x^2}2g_{N/Y}\big(\cR^{\phi}\big(\wt U, \wt U_0\big), V\big)$ \tsep{3pt}\bsep{3pt}\\ \hline
$\nabla_{\wt U_1}\wt U_2$ &
	$\frac x2g_{N/Y}\big(\cR^{\phi}\big(\wt U_1, \wt U_2\big), V_0\big)$ &
	$g_Y\big( \nabla^Y_{U_1} U_2, U_0\big)$ \tsep{3pt}\bsep{3pt}\\ \hline
\end{tabular}
\end{gather*}

We point out a few consequences of these computations. First note that
\begin{gather*}
	\nabla\colon \ \CI(M;T_{\ie}M) \lra \CI(M; T^*M \otimes T_{\ie}M)
\end{gather*}
def\/ines a connection on the incomplete edge tangent bundle. Also note that this connection asymptotically preserves the splitting of $T_{\ie}\sC$ into two bundles
\begin{gather}\label{eq:SecondSplittingIE}
	T_{\ie}\sC = \big[ \ang{\pa_x} \oplus \tfrac1x TN/Y\big] \oplus \phi^*TY
\end{gather}
in that if $W_1, W_2 \in \cV_{\ie}$ are sections of the two dif\/ferent summands then
\begin{gather*}
	g_{\ie}(\nabla_{W_0} W_1, W_2) = \cO(x) \qquad \text{for all} \quad W_0 \in \CI(M;TM).
\end{gather*}
In fact, let us denote the projections onto each summand of \eqref{eq:SecondSplittingIE} by
\begin{gather*}
	 \bv_+\colon \ T_{\ie}\sC \lra \ang{\pa_x} \oplus \tfrac1x TN/Y, \qquad \bh\colon \ T_{\ie}\sC \lra \phi^*TY,
\end{gather*}
and def\/ine connections
\begin{gather*}
	\nabla^{v_+} = \bv_+ \circ \nabla \circ \bv_+\colon \ \CI\big(\sC; \ang{\pa_x} \oplus \tfrac1x TN/Y\big)
		\lra \CI\big(\sC; T^*\sC \otimes \big(\ang{\pa_x} \oplus \tfrac1x TN/Y \big)\big ), \\
	\nabla^h = \phi^*\nabla^Y \colon \ \CI(\sC; \phi^*TY) \lra \CI(\sC; T^*\sC \otimes \phi^*TY ) .
\end{gather*}
Denote by
\begin{gather*}
j_{\eps}\colon \ \{x = \eps\} \hookrightarrow \sC
\end{gather*}
the inclusion, and identify $\{ x=\eps\}$ with $N=\{x=0\}$, note that the pull-back connec\-tions~$j_{\eps}^*\nabla^{v_+}$ and $j_{\eps}^*\nabla^h$ are independent of $\eps$ and
\begin{gather}\label{eq:AsympSplittingConn}
	j_0^*\nabla = j_0^*\nabla^{v_+} \oplus j_0^*\nabla^h.
\end{gather}
In terms of the local connection one-form $\omega$ and the splitting \eqref{eq:SecondSplittingIE}, we have
\begin{gather}
	P^{V^+}\omega = \begin{pmatrix} \omega_{N/Y} & \cO(x) \\ \cO(x) & \cO(x^2) \end{pmatrix},\qquad
	P^H\omega = \begin{pmatrix} \omega_{\cS} & \cO(x) \\ \cO(x) & \phi^*\omega_Y \end{pmatrix},\nonumber\\
	\omega = \begin{pmatrix} \omega_{v_+} & \cO(x) \\ \cO(x) & \phi^*\omega_Y + \cO(x^2) \end{pmatrix},\label{eq:Schemeomega}
\end{gather}
where $P^{V^+}\omega$ is the projection onto the dual bundle of $\lrspar{ \ang{\pa_x} \oplus \tfrac1x TN/Y }$, $P^H\omega$ is the projection onto the dual bundle of $\phi^*TY$, and the forms $\omega_{N/Y}$, $\omega_{\cS}$, $\omega_Y$, $\omega_{v+}$ are def\/ined by these equations. Finally, consider the curvature $R_{\ie}$ of $\nabla$. If $W_1, W_2 \in \cV_{\ie}$ are sections of two dif\/ferent summands of~\eqref{eq:SecondSplittingIE} and $W_3, W_4 \in \CI(\sC,TN/Y \oplus \phi^*TY)$ then
\begin{gather*}
	g_{\ie}(R_{\ie}(W_3,W_4)W_1, W_2) = \cO(x),
\end{gather*}
but
\begin{gather*}
	g_{\ie}(R_{\ie}(\pa_x, W_4)W_1, W_2) = \frac1x g_{\ie}(\nabla_{W_4}W_1, W_2) = \cO(1).
\end{gather*}
We will be interested in the curvature along the level sets of $x$. Schematically, if $\Omega$ denotes the $\End(T_{\ie}M)$-valued two-form corresponding to the curvature of $\nabla$, then with respect to the splitting \eqref{eq:SecondSplittingIE} we have
\begin{gather*}
	\Omega\rest{x =\eps} = \begin{pmatrix} \Omega_{v_+} & \cO(\eps) \\ \cO(\eps) & \phi^*\Omega_Y \end{pmatrix},
\end{gather*}
where $\Omega_{v_+}$ is the tangential curvature associated to $\omega_{N/Y} + \omega_{\cS}$ and $\Omega_Y$ is the curvature associated to $\omega_Y$, and analogously to \eqref{eq:AsympSplittingConn},
\begin{gather}\label{eq:AsympSplittingCurv}
	j_0^*\Omega = j_0^*\Omega_{v_+} + \phi^*\Omega_Y.
\end{gather}

Following \cite{BC1990} and \cite{HHM2004}, it will be convenient to use the block-diagonal connection $\wt\nabla$ on $T_{\ie}M$ from the splitting~\eqref{eq:SplittingIE}. Thus
\begin{gather}\label{eq:euclideanconnection}
	\wt\nabla\colon \ \CI(\sC,T_{\ie}\sC) \lra \CI(\sC; T^*\sC \otimes T_{\ie}\sC)
\end{gather}
satisf\/ies
\begin{gather*}
	\wt\nabla \pa_x =0, \qquad \wt\nabla_{\pa_x} = 0, \qquad \text{and} \\
\begin{tabular}{|c||c|c|} \hline
$g_{\ie}\big(\wt \nabla_{W_1} W_2, W_3\big)$ & $\tfrac1xV_0$ & $\wt U_0$ \tsep{3pt}\bsep{3pt}\\ \hline\hline
$\wt \nabla_{V_1}\tfrac1x V_2$ &
	$g_{N/Y}\big(\nabla^{N/Y}_{ V_1} V_2, V_0\big)$ &
	$0$ \tsep{3pt}\bsep{3pt}\\ \hline
$\wt\nabla_{\wt U} \tfrac1xV$ &
	$g_{N/Y}\big(\big[\wt U, V\big], V_0\big) - \phi^*g_Y\big(\cS^{\phi}( V, V_0), \wt U\big)$ &
	$0$ \tsep{3pt}\bsep{3pt}\\ \hline
$\wt\nabla_{ V} \wt U$ &
	$0$ &
	$0$ \tsep{3pt}\bsep{3pt}\\ \hline
$\wt \nabla_{\wt U_1}\wt U_2$ &
	$0$ &
	$g_Y\big( \nabla^Y_{U_1} U_2, U_0\big)$\tsep{3pt}\bsep{3pt} \\ \hline
\end{tabular}
\end{gather*}
The connection $\wt{\nabla}$ is a metric connection and preserves the splitting \eqref{eq:SecondSplittingIE}.

\subsection{Clif\/ford bundles and Clif\/ford actions}
The incomplete edge Clif\/ford bundle, denoted $\operatorname{Cl}_{\ie}(M, g)$, is the bundle obtained by taking the Clif\/ford algebra of each f\/iber of $T_{\ie}M$. Concretely,{\samepage
\begin{gather*} 
 \operatorname{Cl}_{\ie}(M, g) = \sum_{k = 0}^{\infty} T_{\ie}M^{\otimes k} / (x
 \otimes y + y \otimes x = - 2 \la x, y \ra_{g}).
\end{gather*}
This is a smooth vector bundle on all of $M$.}

\looseness=-1
We assume that $M$ is spin and f\/ix a spin bundle $\SPB \lra M$. Note that $\SPB$ is indeed a smooth vector bundle on all of $M$ (including
the boundary) as the orthonormal frame bundle $\cO \lra M^\circ$ extends smoothly up to the boundary (indeed, consider the local orthonormal frames from the previous section). In fact, the orhonormal frame bundle of $T_{\ie} M$ gives the extension of~$\cO$ to~$M$. Denote Clif\/ford multiplication, which also extends smoothly to all incomplete edge vector f\/ields, by
\begin{gather*}
	c\colon \ \CI(M, T_{\ie}M) \lra \CI(M;\End(\SPB)).
\end{gather*}
We denote the connection induced on $\SPB$ by the Levi-Civita connection $\nabla$ by the same symbol. Let $\eth$ denote the corresponding Dirac operator.

\begin{Lemma}\label{thm:diracwithgoodconnection}
Let $\sC \cong [0,1)_x \times N$ be a collar neighborhood of the boundary. Let
\begin{gather}\label{eq:frame}
	\pa_x, \quad
	\tfrac1x V_\alpha, \quad
	\wt U_i
\end{gather}
denote a local orthonormal frame consistent with the splitting \eqref{eq:SplittingIE}. In terms of this frame and the connection $\wt\nabla$ from~\eqref{eq:euclideanconnection}, the Dirac operator $\eth$ decomposes as
\begin{gather} \label{eq:standardoperator}
	\eth = c(\p_{x}) \p_{x} + \frac{f}{2x} c(\p_{x}) + \frac{1}{x} \sum_{\alpha = 1}^{f} c\big(\tfrac 1x V_{\alpha}\big) \wt{\nabla}_{V_{\alpha}}
	+ \sum_{i = 1}^{b} c(\wt U_{i})\wt{\nabla}_{\wt U_{i}} + B,
\end{gather}
where $f = \dim Z$, $b= \dim Y$, and $B \in \CI(M, \End(\SPB))$.
\end{Lemma}

\begin{proof} Consider the dif\/ference of connections (on the tangent bundle)
\begin{gather*}
	A = \nabla - \wt\nabla \in \CI\big(\sC;T^*M\otimes \End\big(\Ie T\sC\big)\big).
\end{gather*}
From \cite{BFII} we have
\begin{gather*}
 \eth = \sum_{i} c(e_{i})\big(\wt{\nabla}_{e_{i}} + \wt{A}(e_{i})\big),
\end{gather*}
where $\wt{A}(W) := \frac 14 \sum_{jk} g_{\ie}( A(W) e_{j}, e_{k}) c(e_{j})c(e_{k})$. From Section~\ref{sec:Conn} we have
\begin{gather*}
	g_{\ie}(A(W_0)W_1, W_2) =0 \qquad \text{if}\quad \pa_x \in \{ W_0, W_1, W_2 \} \\
	\text{except for} \quad g_{\ie}\big(A\big(\tfrac1xV_\alpha\big)\pa_x, \tfrac1x V_\beta\big) = - g_{\ie}\big(\tfrac1x A(V_\alpha)\tfrac1xV_\beta, \pa_x\big) = \frac1x g_{Z}(V_\alpha, V_\beta),
\end{gather*}
and otherwise
\begin{gather*}
\begin{tabular}{|c||c|c|} \hline
$g_{\ie}(A(W_1) W_2, W_3)$ & $\tfrac1xV_0$ & $\wt U_0$ \tsep{3pt}\bsep{3pt}\\ \hline\hline
$A\big(\tfrac1xV_1\big)\tfrac1x V_2$ &
	$0$ &
	$\phi^*g_Y\big(\cS^{\phi}( V_1, V_2), \wt U_0\big)$ \tsep{3pt}\bsep{3pt}\\ \hline
$A(\wt U) \tfrac1xV$ &
	$0$ &
	$-\frac x2g_{N/Y}\big(\cR^{\phi}\big(\wt U, \wt U_0\big), V\big)$ \tsep{3pt}\bsep{3pt}\\ \hline
$A(\tfrac1xV) \wt U$ &
	$-\phi^*g_Y\big(\cS^{\phi}(V, V_0),\wt U\big)$ &
	$\frac{x}2g_{N/Y}\big(\cR^{\phi}\big(\wt U, \wt U_0\big), V\big)$ \tsep{3pt}\bsep{3pt}\\ \hline
$A(\wt U_1)\wt U_2$ &
	$\frac x2g_{N/Y}\big(\cR^{\phi}\big(\wt U_1, \wt U_2\big), V_0\big)$ &
	$0$ \tsep{3pt}\bsep{3pt}\\ \hline
\end{tabular}
\end{gather*}
Hence $\tfrac14\sum_{s,t,u}g_{\ie}(A(e_s)e_t, e_u)c(e_s)c(e_t)c(e_u)$ has terms of order $\cO(\tfrac1x)$, $\cO(1)$, and $\cO(x)$. The terms of order $\tfrac1x$ are
\begin{gather*}
	\frac1{4} \sum_{\alpha} \big( g_{\ie}\big(A\big(\tfrac1xV_{\alpha}\big)\pa_x, \tfrac1xV_{\alpha}\big) c(\pa_x)
	+ g_{\ie}\big(A\big(\tfrac1xV_{\alpha}\big)\tfrac1xV_{\alpha},\pa_x) (-c(\pa_x))\big) \\
	\qquad{} = \frac1{4}\sum_{\alpha}2\frac{c(\pa_x)}x = \frac{f}{2x}c(\pa_x),
\end{gather*}
which establishes \eqref{eq:standardoperator}.
\end{proof}

\subsection{The APS boundary projection} \label{sec:etaAPS}
We now def\/ine the APS boundary condition discussed in the introduction. We will make use of a simplif\/ied coordinate system near the boundary of $M$, namely, let $(x, x')$ be coordinates near a point on $\p M$ for which $x' \in \mathbb{R}^{n-1}$ are coordinates on $\p M$ and $x$ is the same f\/ixed boundary def\/ining function used in~\eqref{eq:coordinates}. For the cutof\/f manifold $M_{\epsilon} =\{x \geq\eps \}$, consider the dif\/ferential operator on sections of $\SPB$ over $\p M_{\epsilon}$ def\/ined by choosing any orthonormal frame $e_{p}$, $p =1, \dots, n - 1$ of the distribution of the tangent bundle orthogonal to $\p_{x}$ and setting
\begin{gather*}
 \frac{1}{\epsilon} \wt{\eth}_{\epsilon} := \left. - c(\p_{x}) \lp \sum_{p = 1}^{n -
 1}c(e_{p}) \wt{\nabla}_{e_{p}} \rp \right|_{x = \epsilon},
 \end{gather*}
where $\wt{\nabla}$ is the connection from \eqref{eq:euclideanconnection}. The operator $\wt{\eth}_{\epsilon}$ is def\/ined independently of the choice of frame, so we may take frames as in \eqref{eq:frame} to obtain
 \begin{gather} \label{eq:boundaryoperator}
 \wt{\eth}_{\epsilon}
= \left. - x c(\p_{x}) \lp \sum_{\alpha = 1}^{f} c\big(\tfrac 1x V_{\alpha}\big) \wt{\nabla}_{\frac 1x V_{\alpha}} + \sum_{j =
 1}^{b} c\big(\wt U_{j}\big) \wt{\nabla}_{\wt U_{j}} \rp \right|_{x = \epsilon}.
 \end{gather}
We refer to $\wt{\eth}_{\epsilon}$ below as the \textit{tangential operator}, since for every $\epsilon$ it acts tangentially along the boundary $\p M_{\epsilon}$. The operator $\wt{\eth}_{\epsilon}$ is self-adjoint on $L^{2}(\p M_{\epsilon}, \SPB)$.

We denote the dual coordinates on $T^{*}_{x,x'}M$ by $(\xi, \xi')$. Using the identif\/ication of $T^{*}M$ with $TM$ induced by the metric $g$, the principal symbol of $\eth$ is given by
\begin{gather} \label{eq:symbolofDirac}
 \sigma(\eth)(x, x') = i\xi c(\p_{x}) + ic(\xi' \cdot \p_{x'}),
\end{gather}
where $\xi' \cdot \p_{x'} = \sum\limits_{i = 1}^{n - 1} \xi'_{j} \p_{x'_{j}}$. Note that using coordinates as in~\eqref{eq:symbolofDirac}, for $x = \epsilon$, $\wt{\eth}_{\epsilon}$ has principal symbol
\begin{gather*}
	\sigma(\wt{\eth}_{\epsilon})(x', \xi') = - i \epsilon c(\p_{x}) c(\xi' \cdot \p_{x'}).
\end{gather*}
Since $\sigma(\wt{\eth}_{\epsilon})(x', \xi')^{2} = \epsilon^{2} \absv{(0, \xi')}_{g}^{2}$, if we def\/ine $\hat{\xi}' = \xi'/\absv{(0, \xi')}_{g}$, then
\begin{gather} \label{eq:boundarysymbolasprojections}
 \sigma(\wt{\eth}_{\epsilon})(x', \xi') = \absv{(0, \xi')}_{g} \sigma(\wt{\eth}_{\epsilon})(x', \hat{\xi}')
 = \epsilon \absv{(0, \xi')}_{g}(\pi_{\epsilon, +, \hat{\xi}'}(x') - \pi_{\epsilon, -, \hat{\xi}'}(x')),
\end{gather}
where $\pi_{\epsilon, \pm, \hat{\xi}'}(x')$ are orthogonal projections onto $\pm$ eigenspaces of $\sigma(\wt{\eth}_{\epsilon})(x', \hat{\xi'})$. We will def\/ine a~boundary condition for $\eth$ on the cutof\/f manifolds $M_{\epsilon}$,
\begin{gather}
 \label{eq:APSprojector}
 \pi_{{\rm APS}, \epsilon} := L^{2} \quad \mbox{orthogonal projection onto} \ \ V_{-, \epsilon},
\end{gather}
where $V_{-, \epsilon}$ is the direct sum of eigenspaces of $\wt{\eth}_{\epsilon}$ with negative eigenvalues. We recall basic facts about $\pi_{{\rm APS},\epsilon}$.
\begin{Theorem}[\cite{APSI, S1966}]\label{thm:APSbasics}
For fixed $\epsilon$, the operator $\pi_{{\rm APS}, \epsilon}$ is a pseudodifferential operator of order~$0$, i.e., $\pi_{{\rm APS}, \epsilon} \in \pdo^{0}(\p M_{\epsilon}; \SPB)$. Its principal symbol satisf\/ies
\begin{gather*}
	\sigma(\pi_{{\rm APS}, \epsilon})(x', \xi') = \pi_{\epsilon, -, \hat{\xi}'}(x'),
\end{gather*}
where $\pi_{\epsilon, -, \hat{\xi'}}(x')$ is projection onto the negative eigenspace of $- i c(\p_{x}) c(\hat{\xi}' \cdot \p_{x'})$ from~\eqref{eq:boundarysymbolasprojections}.
\end{Theorem}

Consider the domains for $\eth$ on $L^{2}(M_{\epsilon}; \SPB)$ def\/ined as follows
\begin{gather}
	\mathcal{D}_{{\rm APS}, \epsilon}:= \set{u \in H^{1}(M_{\epsilon}, \SPB) \colon ({\operatorname{Id}} - \pi_{{\rm APS}, \epsilon}) u = 0 },\nonumber \\
	\mathcal{D}^{+}_{{\rm APS}, \epsilon}:= \set{u \in \mathcal{D}_{{\rm APS}, \epsilon}\colon \operatorname{image}(u) \subset \SPB^{+}}.\label{eq:APSproblem}
\end{gather}
In one of the main results of this paper, we will show that $\eth$ has a unique self-adjoint exten\-sion~$\mathcal D$, such that for $\eps>0$ suf\/f\/iciently small,
\begin{gather*}
	\ind\big(\eth \colon \mathcal{D}^{+}_{{\rm APS}, \epsilon} \lra L^{2}(M_{\epsilon}; \SPB^{-})\big)
	= \ind\big(\eth \colon \mathcal{D}^{+} \lra L^{2}(M; \SPB^{-})\big),
\end{gather*}
where $\mathcal{D}^+$ the elements of $\mathcal{D}$ valued in $\SPB^{+}$. Indeed, this will follow from Theorems~\ref{thm:naturalbvp} and~\ref{thm:APSmatchesother} below.

\section[Mapping properties of $\eth$]{Mapping properties of $\boldsymbol{\eth}$}\label{sec:parametrix}

In this section we will use the results and techniques in \cite{ALMPI, ALMPII} to prove Theorem~\ref{thm:essential}. We proceed by constructing a parametrix
for $\eth$ and analyzing the mapping properties of this parametrix.

Let $\mathcal{D}$ denote the domain of the unique self-adjoint extention of $\eth$. At the end of this section, we analyze the structure of the generalized inverse~$Q$ for~$\eth$, that is, the map
\begin{gather*}
	Q \colon \ L^{2}(M; \SPB) \lra \mathcal{D} \qquad \text{satisfying} \quad \eth Q = {\operatorname{Id}} - \pi_{\operatorname{ker}}
\quad \text{and} \quad Q = Q^{*},
\end{gather*}
where $\pi_{\operatorname{ker}}$ is $L^{2}$-orthogonal projection onto the kernel of $\eth$. Here the adjoint~$Q^{*}$ is taken with respect to the pairing def\/ined for sections $\phi$, $\psi$ by
\begin{gather*}
	\la \phi, \psi \ra_{L^{2}} = \int \la \phi, \psi \ra_{G} \, d{\rm Vol}_{g},
\end{gather*}
where $G$ is the Hermitian inner product on $\SPB$.

\subsection{The ``geometric Witt condition''}

The proof of Theorem \ref{thm:essential} relies on an assumption on an induced family of Dirac operators on the f\/iber $Z$ which we describe now. By Lemma~\ref{thm:diracwithgoodconnection}, on a~collar neighborhood of the boundary,
\begin{gather}\label{eq:collarneighborhood}
	\mathcal{U} \subset M \qquad \text{with} \quad \mathcal{U} \simeq [0,\epsilon_{0})_{x} \times \p M
\end{gather}
we can write
\begin{gather}\label{eq:operatorwithparts}
 \eth = c(\p_{x}) \lp \p_{x} + \frac{f}{2x} + \frac{1}{x} \eth^{Z}_{y} - \sum_{i = 1}^{b} c(\p_{x})c(\wt U_{i})\wt{\nabla}_{\wt U_{i}}\rp + B
\end{gather}
with $\norm{B} = O(1)$ and where, for $y$ in the base~$Y$
\begin{gather}\label{eq:fiberop}
\eth^{Z}_{y}= - c(\p_{x}) \cdot \sum_{\alpha = 1}^{f} c\big(\tfrac 1x V_{\alpha}\big) \cdot \wt{\nabla}_{V_{\alpha}}.
\end{gather}
The operator $\eth^{Z}_{y}$ def\/ines a self-adjoint operator on the f\/iber over $y \in Y$ in the boundary f\/ibration $N\xlra{\phi}Y$ acting on sections of the restriction of the spin bundle $\SPB_{y}$.

We will assume the following ``geometric Witt condition" discussed in the introduction.
\begin{Assumption}\label{thm:assumption}
The fiber operator $\eth^{Z}_{y}$ in \eqref{eq:fiberop} satisfies
\begin{gather*} 
 (-1/2, 1/2) \cap \spec\big(\eth^{Z}_{y}\big) = \varnothing \qquad \mbox{for all} \quad y.
\end{gather*}
\end{Assumption}

\subsection{Review of edge and incomplete edge operators}\label{sec:reviewedge}

A vector f\/ield on $M$ is an `edge vector f\/ield' if its restriction to $N = \pa M$ is tangent to the f\/ibers of~$\phi$~\cite{Ma1991}. A~dif\/ferential operator is an edge dif\/ferential operator if in every coordinate chart it can be written as a polynomial in edge vector f\/ields. Thus if $E$ and $F$ are vector bundles over~$M$, we say that $P'$ is an $m^{\text{th}}$ order edge dif\/ferential operator between sections of~$E$ and~$F$, denoted $P' \in \Diff_{\rm e}^m(M;E,F)$, if in local coordinates we have
\begin{gather*}
P' = \sum_{j+|\alpha|+|\gamma|\leq m} a_{j,\alpha,\gamma}(x,y,z) (x\pa_x)^j(x\pa_y)^{\alpha}(\pa_z)^{\gamma},
\end{gather*}
where $\alpha$ denotes a multi-index $(\alpha_1, \ldots, \alpha_b)$ with $|\alpha|=\alpha_1+\cdots+\alpha_b$ and similarly for $\gamma=(\gamma_1, \ldots, \gamma_f)$, and each $a_{j,\alpha,\gamma}(x,y,z)$ is a local section of $\hom(E,F)$.

A dif\/ferential operator $P$ is an `incomplete edge dif\/ferential operator' of order $m$ if $P'=x^mP$ is an edge dif\/ferential operator of order~$m$. Thus, symbolically,
\begin{gather}\label{eq:ieoperator'}
	\Diff_{\ie}^m(M;E,F) = x^{-m}\Diff_{\rm e}^m(M;E,F),
\end{gather}
and in local coordinates
\begin{gather*}
	P = x^{-m} \sum_{j+|\alpha|+|\gamma|\leq m} a_{j,\alpha,\gamma}(x,y,z) (x\pa_x)^j(x\pa_y)^{\alpha}(\pa_z)^{\gamma}.
\end{gather*}
The (incomplete edge) principal symbol of $P$ is def\/ined on the incomplete edge cotangent bundle,
\begin{gather*}
	\sigma(P) \in \CI(T_{\ie}^*M;\pi^*\hom(E,F)),
\end{gather*}
where $\pi\colon T_{\ie}^*M\lra M$ denotes the bundle projection. In local coordinates it is given by
\begin{gather*}
	\sigma(P)(x,y,z,\xi, \eta, \zeta):=	\sum_{j+|\alpha|+|\gamma| = m} a_{j,\alpha,\gamma}(x,y,z) (\xi)^j(\eta)^{\alpha}(\zeta)^{\gamma}.
\end{gather*}
We say that $P$ is elliptic if this symbol is invertible whenever $(\xi,\eta,\zeta)\neq0$.

\begin{Lemma}\label{thm:diracisincompleteedge}
The Dirac operator $\eth$ on an incomplete edge space is an elliptic incomplete edge differential operator of order $1$, i.e., is an elliptic element of $\Diff^{1}_{\ie}(M; \SPB)$. In particular, $x\eth$ is an elliptic element of $\Diff^{1}_{\rm e}(M; \SPB)$.
\end{Lemma}
\begin{proof}
This follows from equation \eqref{eq:standardoperator} in Lemma~\ref{thm:diracwithgoodconnection}.
\end{proof}

\subsection[Parametrix of $x\eth$ on weighted edge spaces]{Parametrix of $\boldsymbol{x\eth}$ on weighted edge spaces}

Lemma~\ref{thm:diracisincompleteedge} shows that $x\eth$ is an elliptic edge operator. By the theory of edge operators \cite{Ma1991}, this implies that $x\eth$ is a bounded operator between appropriate weighted Sobolev spaces, whose def\/inition we now recall.

Let $\mathcal{D}'(M; \SPB)$ denote distributional sections. Given $k \in \mathbb{N}$, let
\begin{gather*}
 H^{k}_{\rm e}(M; \SPB) := \set{u \in \mathcal{D}'(M; \SPB)\colon A^{1} \cdots A^{j} u \in L^{2}(M; \SPB) \mbox{ for } j \le k \mbox{ and } A^{i} \in \Diff^{1}_{\rm e}(M; \SPB)}.
\end{gather*}
In particular, $u \in H^{1}_{\rm e}(M; \SPB)$ if and only if $u \in L^2(M; \SPB)$ and, for any edge vector f\/ield $V \in C^{\infty}(M; T_{\rm e}M)$, $\nabla_{V} u \in L^{2}(M; \SPB)$. The weighted edge Sobolev spaces are def\/ined by
\begin{gather*}
	x^{\delta}H^{k}_{\rm e}(M; \SPB) := \big\{u \colon x^{-\delta}u \in H^{k}_{\rm e}(M; \SPB)\big\}.
\end{gather*}
Thus, the map
\begin{gather}\label{eq:epsilonforedge}
	x \eth \colon \ x^{\delta}H^{k}_{\rm e}(M; \SPB) \lra x^{\delta}H^{k-1}(M; \SPB)
\end{gather}
is bounded for all $\delta \in \mathbb{R}$, $k \in \mathbb{N}$, in fact for $k \in \mathbb{R}$ by duality and interpolation. We will prove the following
\begin{Proposition}\label{thm:edgeparam}
 Under the Witt assumption $($Assumption {\rm \ref{thm:assumption})}, the map \eqref{eq:epsilonforedge} is Fredholm for $0 < \delta < 1$.
\end{Proposition}

\begin{Remark} Recall that the space $L^{2}(M; \SPB)$ used to def\/ine the $x^{\delta}H^{k}(M ; \SPB)$ is equipped with the inner product from the Hermitian metric on $\SPB$ and the volume form of the \textit{incomplete} edge metric~$g$.
\end{Remark}
\begin{Remark}[complete edge manifolds]\label{thm:edge-theorem}
We now give a rough sketch of how one can use the main theorem of this paper to prove an index formula for Dirac operators in the complete edge case, leaving the details to the reader. Note that for each incomplete edge metric $g$ on $M$ as above, there is an edge metric $\wt{g}$ def\/ined by
\begin{gather*}
g = x^{-2} \wt{g},
\end{gather*}
where $x$ remains a boundary def\/ining function as above, and in particular is non-vanishing in the interior of~$M$. Choosing a spin structure, the Dirac operator of $g$ is related to that of~$\wt{g}$ (under a natural identif\/ication of the spin bundles) by
\begin{gather*}
\eth_{\wt{g}} = x^{(n + 1)/2} \, \eth_{\wt{g}} \, x^{- (n - 1)/2}.
\end{gather*}
One can show that, under the strengthened assumption, $\eth_{\wt{g}} \colon H^{1}_{\rm e}(M; \SPB) \lra L^{2}(M; \SPB)$ is Fredholm (i.e., one can take the weight $\delta = 0$) and that the index of this map is equal to that of $\eth_{g}$. The main theorem then gives an index formula for $\eth_{\wt{g}}$ in this case.
\end{Remark}

In fact we will need more than Proposition \ref{thm:edgeparam}; the proof of the Main Theorem requires a~detailed understanding of the structure of parametrices for $x \eth$. To understand these, we must recall of edge double space $M^{2}_{\rm e}$, depicted heuristically in Fig.~\ref{fig:edgedouble} below, and edge pseudodif\/ferential operators, def\/ined in~\cite{Ma1991} with background material in~\cite{damwc}. The edge double space $M^{2}_{\rm e}$ is a~mani\-fold with corners, obtained by radial blowup of $M \times M$, namely $M^{2}_{\rm e} := [M \times M ; \diag_{\rm f\/ib}(\p M \times \p M)]$, where the notation is that in~\cite{damwc}. Here $\diag_{\rm f\/ib}(\p M \times \p M)$ denotes the f\/iber diagonal
\begin{gather}\label{eq:fiberdiagonal}
 \diag_{\rm f\/ib}(N \times N) = \set{(p, q) \in N \times N \colon \phi(p) = \phi(q) },
\end{gather}
where $\phi\colon N \lra Y$ is f\/iber bundle projection onto $Y$. Whereas $M \times M$ is a manifold with corners with two boundary hypersurfaces, $M^{2}_{\rm e}$ has a third boundary hypersurface introduced by the blowup.
\begin{gather*}
	\mbox{Let $\ff$ be the boundary hypersurface of $M^{2}_{\rm e}$ introducedby the blowup.}
\end{gather*}
Furthermore we have a blowdown map.
\begin{gather*}
 \beta \colon \ M^{2}_{\rm e} \lra M \times M,
\end{gather*}
which is a $b$-map in the sense of \cite{damwc} and a dif\/feomorphism from the interior of $M^{2}_{\rm e}$ to that of $M \times M$. The two boundary hypersurfaces of $M \times M$, $\set{x = 0}$ and $\set{x' = 0}$, lift to boundary hypersurfaces of $M^{2}_{\rm e}$ which we denote by
\begin{gather*}
\lf := \overline{\beta^{-1}\big(\set{x = 0}^{\rm int}\big)} \qquad \text{and} \qquad \rf := \overline{\beta^{-1}\big(\set{x' = 0}^{\rm int}\big)}.
\end{gather*}

The edge front face, $\ff$, is the radial compactif\/ication of the total space of a f\/iber bundle
\begin{gather*}
	Z^{2} \times \mathbb{R}^{b} \times \mathbb{R}_{+} \hookrightarrow \ff \lra Y,
\end{gather*}
where $b = \dim Y$. This bundle is obtained by pulling-back the tangent bundle of $Y$ and trivial $\bbR^+$ bundle, to the f\/iber product of two copies of $\p M$, along the natural projection to $Y$. Choosing local coordinates $(x, y, z)$ as in \eqref{eq:coordinates}, and our f\/ixed bdf $x$, and letting $(x, y, z, x', y', z')$ denote coordinates on $M \times M$ with $(x', y', z')$ the same functions as $(x, y, z)$ but on the right factor, the functions
\begin{gather}\label{eq:projectivecoords}
	x',\quad \sigma = \frac{x}{x'},\quad \YY = \frac{y - y'}{x'},\quad y',\quad z,\quad z'
\end{gather}
def\/ine coordinates near $\ff$ in the set $0 \le \sigma < \infty$, and in these coordinates $x'$ is a boundary def\/ining function for $\ff$, meaning that $\set{x' = 0}$ coincides with $\ff$ on $\{ 0 \le \sigma < \infty \} = M^2_{\rm e} \setminus \rf$, and~$x'$ has non-vanishing dif\/ferential on~$\ff$. When $x' = 0$, $\sigma$ gives coordinates on the $\mathbb{R}_{+}$ f\/iber, $\YY$ on the $\mathbb{R}^{b}$ f\/iber, and $z$, $z'$ on the $Z^{2}$ f\/iber. Below we will also use cylindrical coordinates near~$\ff$. These have the advantage that they are def\/ined on open neighborhoods of sets in~$\ff$ which lie over open sets $V$ in the base $Y$. With $(x,y,z)$ as above, let
\begin{gather*}
	\rho_{\ff} = \big(x^{2} + (x')^{2} + \absv{y - y'}^{2}\big)^{1/2}, \qquad
	\phi = \lp \frac{x}{\rho_{\ff}}, \frac{x'}{\rho_{\ff}}, \frac{y - y'}{\rho_{\ff}}\rp,
\end{gather*}
so $(\rho_{\ff}, \phi, y', z, z')$ form cylindrical coordinates (in the sense that $\absv{\phi}^{2} = 1$) near the lift of $V$ to $\ff$ and in the domain
of validity of $y$, $z$.

\begin{figure} \centering
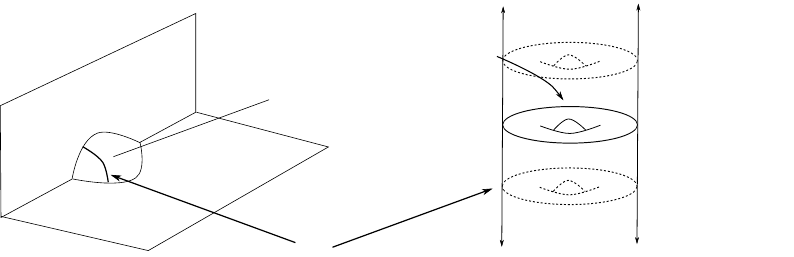
 \caption{}\label{fig:edgedouble}
\end{figure}

We now def\/ine the calculus of edge pseudodif\/ferential operators with bounds, which is similar to the large calculus of pseudodif\/ferential edge operators def\/ined in~\cite{Ma1991}. Thus, $\pdo^{m}_{{\rm e}, {\rm bnd}}(M; \SPB)$ will denote the set of operators $A$ mapping $C^{\infty}_{c}(M; \SPB)$ to distributional sections $\mathcal{D}'(M; \SPB)$, whose Schwartz kernels have the following structure. Let $\End(\SPB)$ denote the bundle over $M \times M$ whose f\/iber at $(p, q)$ is $\Hom(\SPB_{q}; \SPB_{p})$. The Schwartz kernel of~$A$, $K_{A}$ is a distributional section of the bundle $\End(\SPB)$ over $M \times M$ satisfying that for a section $\phi \in C^{\infty}_c(M ; \SPB)$,
\begin{gather}\label{eq:schwartzkernelacts}
	A \phi(w) = \int_{M} K_{A}(w, w') \phi(w')\, d{\rm Vol}_{g}(w'),
\end{gather}
where $d{\rm Vol}_{g}$ is the volume form of an \textit{incomplete} edge metric $g$ asymptotically of the form~\eqref{eq:incompleteedgemetric}. Moreover,
\begin{gather*}
	K_{A} \in \mathcal{A}_{a,b, f}I^{m}\big(M^{2}_{\rm e}, \Delta_{\rm e} ; \beta^{*}\End(\SPB)\big),
\end{gather*}
meaning that $K_{A} = K_{1} + K_{2}$ where $\rho_{\ff}^{f}K_{1}$ is in the H\"ormander conormal space \cite[Chapter~18]{Hvol3}
\begin{gather*}
	\rho_{\ff}^{f} K_{1} \in I^{m}\big(M^{2}_{\rm e}, \Delta_{\rm e} ; \beta^{*}\End(\SPB)\big),
\end{gather*}
$K_{1}$ is supported near $\Delta_{\rm e}$, and $K_{2} \in \mathcal{A}(M^2_{\rm e}; \End(\SPB))$, i.e., is smooth in the interior and conormal to the boundary, and satisf\/ies the bounds
\begin{gather}\label{eq:kernelbounds}
 K_{2}(p) = O\big(\rho_{\lf}^{a}\big) \qquad \text{as} \quad p \to \lf, \qquad K_{2}(p) = O\big(\rho_{\rf}^{b}\big) \qquad \text{as} \quad p \to \rf,
\end{gather}
where $\rho_{\lf}$, $\rho_{\rf}$, and $\rho_{\ff}$ are boundary def\/ining functions for $\lf$, $\rf$, and $\ff$ respectively. The bound is in the norm on $\End(\SPB)$ over $M \times M$, see~\cite{Ma1991} for details. Since the bounds $a$ and $b$ in~\eqref{eq:kernelbounds} will be of some importance, we let
\begin{gather*}
	\pdo^{m}_{\rm e}(M ; \SPB; a, b)
\end{gather*}
denote the subspace of $\Psi^{m}_{{\rm e}, {\rm bnd}}(M; \SPB)$ of pseudodif\/ferential edge operators whose Schwartz kernels satisfy~\eqref{eq:kernelbounds} with bounds~$a$ and~$b$.

The bounds in \eqref{eq:kernelbounds} determine the mapping properties of~$A$ on weighted Sobolev spaces. From \cite[Theorem~3.25]{Ma1991}, we have
\begin{Theorem}
An element $A \in \Psi^{m}_{\rm e}(M; \SPB; a, b)$ is bounded as a map
\begin{gather*}
	A \colon \ x^{\delta} H^{k}_{\rm e}(M; \SPB) \lra x^{\delta'}H^{k - m}(M; \SPB),
\end{gather*}
if and only if
\begin{gather*}
	a > \delta' - f/2 - 1/2 \qquad \text{and} \qquad b > - \delta - f/2 - 1/2.
\end{gather*}
\end{Theorem}

\begin{Remark} In Mazzeo's paper \cite{Ma1991} the convention used to describe the weights (orders of vanishing) of the Schwartz kernels of elements in $\Psi^{m}_{\rm e}$ is slightly dif\/ferent from ours. There one chooses a half-density $\mu$ on $M$ which looks like $\sqrt{dx dy dz}$ near $\p M$. The choice of $\mu$ gives an isomorphism between the sections of $\SPB$ and the sections of $\SPB \otimes \Omega^{1/2}(M)$ where $\Omega^{1/2}(M)$ is the half-density bundle of $M$ (simply by multiplying by $\mu$), and the Schwartz kernel of an edge pseudodif\/ferential operator, $A$, in this context is the section $\kappa_{A}$ of $\End(\SPB) \otimes \Omega^{1/2}(M^{2}_{\rm e})$ with the property that
\begin{gather}\label{eq:edgeconvention}
	A (\psi \mu) = \int_{M} \kappa_{A} \psi \mu.
\end{gather}
One nice feature of \eqref{eq:edgeconvention} is that $\kappa_{A}$ is smooth (away from the diagonal) down to $\ff$. With our convention in~\eqref{eq:schwartzkernelacts}, it is singular of order $-f$ due to the factors of~$x$ in the volume form of~$g$.
\end{Remark}

Given an elliptic edge operator $\wt{P} \in \Diff_{\rm e}^{m}(M; \SPB)$, to construct a parametrix for $\wt{P}$ one must study two models for~$\wt{P}$, the \textit{indicial family} $I_{y}(\wt{P}, \zeta)$ and the \textit{normal operator} $N(\wt{P})_{y}$, in addition to inverting the principal symbol.

First we discuss the indicial operator. For each $y$ in the base~$Y$, the indicial family $\zeta \mapsto I_{y}(\wt{P}, \zeta)$ is an elliptic operator-valued function on $\mathbb{C}$ obtained by taking the Mellin transform (see \cite[Section~2]{Ma1991}) of the leading order part of $\wt{P}$ in $x$. By \eqref{eq:operatorwithparts}, the leading order part of $\wt{P} = x \eth$ is $ c(\p_{x}) \lp x \p_{x} + f/2 + \eth^{Z}_{y} \rp$, so taking the Mellin transform and ignoring the $c(\p_{x})$ gives
\begin{gather}\label{eq:mellin}
	i\zeta + f/2 + \eth^{Z}_{y}.
\end{gather}
The meaningful values of $\zeta$ are the \textit{indicial roots}, which we def\/ine to be
\begin{gather}\label{eq:23}
	\Lambda_{y} = \set{i\zeta + f/2 + 1/2 \colon \eqref{eq:mellin} \ \text{is not invertible}}.
\end{gather}
By def\/inition, \eqref{eq:mellin} is invertible as long as $-(i\zeta + f/2) \not \in \sigma(\eth^{Z}_{y})$, so Assumption~\ref{thm:assumption} implies that
\begin{gather*}
 \Lambda_{y} \cap [0, 1] \subset \set{0, 1} \qquad \text{for all} \quad y \in Y.
\end{gather*}

\begin{Remark}
The shift by $f/2 + 1/2$ in \eqref{eq:23} comes from the following considerations. We want to understand the mapping properties of $x\eth$ on $L^{2}(M; \SPB)$ with the natural measure $d{\rm Vol}_{g}$ given by the incomplete edge metric~$g$. On the other hand, the values of $i\zeta$ for which~\eqref{eq:mellin} fails to be invertible give information about the mapping properties of~$x \eth$ on the Sobolev spaces def\/ined with respect to the $b$-measure
\begin{gather*}
	\mu_{b} := x^{-f-1} \,d{\rm Vol}_{g}.
\end{gather*}
In particular, the Fredholm property in Proposition~\ref{thm:edgeparam} is equivalent to $x \eth$ being a Fredholm map from the space $x^{\delta - f/2 - 1/2}H^{1}_{\rm e}(M ; \SPB ; \mu_{b})$ to the space $x^{\delta - f/2 - 1/2}L^{2}_{\rm e}(M; \SPB ; \mu_{b})$, where the Sobolev spaces are now def\/ined with respect
to the $b$-measure. Alternatively, as in \cite{ALMPII} we could def\/ine $\wt{P}' = x^{-f/2 - 1/2}(x \eth)x^{f/2 + 1/2}$ take the Mellin transform and use the values of $i\zeta$ as the indicial roots, but we would get the same answer as in~\eqref{eq:23}.
\end{Remark}

Now we discuss the normal operator $N(\wt{P})$. Elements of $\cV_{\rm e}$ acting on either factor of $M \times M$ are tangent to $\ff$ when lifted to $M^2_{\rm e}$. This implies that, letting $\wt{P}$ act on the left on $M \times M$ and lifting to $M^2_{\rm e}$, $\wt{P}$ def\/ines an operator on sections over $\ff$. In fact, we can see explicitly that $\wt{P}$ acts on the f\/ibers of $\ff$, so the base $Y$ enters its action only parametrically; that is, for every $y \in Y$, $\wt{P}$ def\/ines an operator
\begin{gather*}
	N\big(\wt{P}\big)_{y} \ \ \text{acting on the f\/iber $\ff_{y}$ over $y$}.
\end{gather*}
To obtain an expression for $N(\wt{P})_{y}$ in coordinates, write
\begin{gather*}
	\wt{P} = \sum_{i+ \absv{\alpha} + \absv{\beta} \le m} a_{i,\alpha,\beta}(x,y,z) (x\p_{x})^{i}
	(x\p_{y})^{\alpha}\p_{z}^{\beta}, \qquad \text{where} \quad a_{i,\alpha,\beta} \in C^{\infty}(M; \End \SPB),
\end{gather*}
and use the projective coordinates in \eqref{eq:projectivecoords} to write
\begin{gather*}
	N\big(\wt{P}\big)_{y} = \sum_{i+ \absv{\alpha} + \absv{\beta} \le m} a_{i,\alpha,k}(0,y,z) (\sigma\p_{\sigma})^{i}
	(\sigma\p_{\YY})^{\alpha}\p_{z}^{\beta}.
\end{gather*}
The mapping properties of $\wt{P}$ are deduced from mapping properties of the $N(\wt{P})_{y}$. In particular, to prove Proposition~\ref{thm:edgeparam} we will need Lemma~\ref{thm:normaloperatorisinvertible} below, which shows that the Fourier transform of $N_{y}(x\eth)$ is invertible on certain spaces.

Edge \textit{pseudo}dif\/ferential operators also admit normal operators. Given $A \in \pdo^{m}_{{\rm e}, {\rm bds}}(M ; \SPB)$, the restriction $N(A) :=
\rho_{\ff}^{f} K_{A} \rvert_{\ff}$ is well def\/ined, and in fact
\begin{gather*}
	N(A) \in \mathcal{A}_{a,b}I^{m}\big(\ff, \Delta_{\rm e} \rvert_{\ff} ; \beta^{*}\End(\SPB)\rvert_{\ff}\big),
\end{gather*}
meaning that $N(A) = \kappa_{1} + \kappa_{2}$ where $\kappa_{1}$ is a~distribution on $\ff$ conormal to $\Delta_{\rm e} \cap \ff$ of order $m$ and $\kappa_{2}$ is a smooth function on $\ff^{int}$ with bounds in~\eqref{eq:kernelbounds} (with the point~$p$ restricted to~$\ff$).

Using \eqref{eq:operatorwithparts} and the projective coordinates in \eqref{eq:projectivecoords}, and letting $c_{\nu}$ denote the operator induced by $c(\p_{x})$ on the bundle $\SPB_{y}$, the restriction of the spin bundle to the f\/iber over $y$, the normal operator of $x \eth$ satisf\/ies
\begin{gather}\label{eq:normaloperatorxeth}
	N(x \eth) = c_{\nu} \cdot \lp \sigma\frac{\p}{\p \sigma}	+ \frac{f}{2} + \eth^{Z}_{y'} \rp + \sigma \eth_{\YY},
\end{gather}
where $\eth_{\YY}$ can be written locally in terms of the limiting base metric $h_{y} = g_Y \rvert_{T_yY}$,, i.e., the translation invariant metric on the vector space $T_yY$ def\/ined by $g_Y$ in~\eqref{eq:incompleteedgemetric} as
\begin{gather*}
	\eth_{\YY} = \sum_{i,j = 1}^{\dim Y} c(\p_{\YY_{i}}) h_{y}^{ij}\p_{\YY_{j}}.
\end{gather*}
The operator $N(x \eth)$ acts on sections of $\SPB_{ y}$.

The remainder of this subsection consists in establishing the following theorem.
\begin{Theorem}\label{thm:edgetheorem}
 Let $0 < \delta < 1$. Under Assumption {\rm \ref{thm:assumption}}, there exist left and right parametri\-ces~$\wt{Q}_{i}$, $i = 1,2$ for $x\eth$. Precisely, there are operators $\wt{Q}_{i} \in \pdo^{-1}_{{\rm e}, {\rm bnd}}(M ; \SPB)$, satisfying
 \begin{gather}\label{eq:edgeparametrix}
\wt{Q}_{1} x\eth= {\operatorname{Id}} - \Pi_{\operatorname{ker}, \delta} \qquad \text{and}\qquad x\eth\wt{Q}_{2} = {\operatorname{Id}} - \Pi_{\operatorname{coker}, \delta},
\end{gather}
where
\begin{gather*}
	\wt{Q}_{i} \colon \ x^{\delta} H^{k}_{\rm e}(M; \SPB) \lra x^{\delta}
	H^{k + 1}_{\rm e}(M; \SPB), \qquad
	\Pi_{\operatorname{ker}, \delta}, \Pi_{\operatorname{coker}, \delta} \colon \ x^{\delta} H^{k}_{\rm e}(M; \SPB) \lra x^{\delta} H^{\infty}_{\rm e}(M; \SPB)
\end{gather*}
for any $k$. Here $\Pi_{\operatorname{ker}, \delta}$ $($resp.\ $\Pi_{\operatorname{coker}, \delta})$ is $x^{\delta}L^{2}(M; \SPB)$
orthogonal projection onto the kernel $($resp.\ cokernel$)$ of~$x \eth$. The Schwartz kernels satisfy the following bounds
\begin{gather*}
	\wt{Q}_{i} \in \pdo^{-1}_{\rm e}(M; \SPB; a, b), \qquad \Pi_{\operatorname{ker}, \delta}, \Pi_{\operatorname{coker}, \delta} \in
	\pdo^{-\infty}_{\rm e}(M ; \SPB; a, b),
\end{gather*}
where $a > \delta - f/2 - 1/2$ and $b > -\delta - f/2 - 1/2$.
Furthermore, $N(\Pi_{\operatorname{ker}, \delta}) \equiv 0 \equiv N(\Pi_{\operatorname{coker}, \delta}) $ so we have
\begin{gather}\label{eq:normaloperatorequation}
	N\big(\wt{Q_{1}}\big) N(x\eth) = N(\operatorname{Id}) = N(x\eth) N\big(\wt{Q}_{2}\big).
\end{gather}
In particular this establishes that $x\eth \colon x^{\delta} H^{k + 1}_{\rm e}(M; \SPB) \lra x^{\delta} H^{k}_{\rm e}(M; \SPB)$ is Fredholm.
\end{Theorem}

We will see that Theorem \ref{thm:edgetheorem} can be deduced from the work of Mazzeo in~\cite{Ma1991} and its modif\/ications in \cite{ALMPI, ALMPII, MV2013}. (See also~\cite{KM2014}, which is closely related to~\cite{MV2013}, though not directly used in the current work.) In order to see that the results of those papers apply, we must prove that the normal operator $N(x\eth)$ is invertible in a~suitable sense. Taking the Fourier transform of~$N_{y}(x \eth)$ in~\eqref{eq:normaloperatorxeth} in the~$\YY$ variable gives
\begin{gather}\label{eq:11}
	L(y, \eta) := \widehat{N_{y}(x \eth)}(\sigma, \eta, z) = c_{\nu} \cdot \lp \sigma
	\frac{\p}{\p \sigma} + \frac{f}{2} +	\eth^{Z}_{y} \rp + i \sigma c(\eta)
\end{gather}
and for each $y$, one considers the mapping of weighted edge Sobolev spaces def\/ined by picking a positive cutof\/f function $\phi \colon [0, \infty)_{\sigma} \lra \mathbb{R}$ that is $1$ near zero and~$0$ near inf\/inity and letting
\begin{gather}\label{eq:Fredholmnormalspaces}
	\mathcal{H}^{r, \delta, l} := \big\{ u \in
	\mathcal{D}'([0, \infty)_{\sigma} \times Z; \SPB_{y}) \colon \phi u \in
	\sigma^{\delta} H^{r},\, (1 - \phi) u \in \sigma^{-l} H^{r}\big\},
\end{gather}
where, in terms of $k_{y} = g_{N/Y}\rest{\phi^*y}$, the f\/iber metric on the f\/iber above $y \in Y$,
\begin{gather*}
	H^{r} := H^{r}\big(\sigma^{f} d\sigma d{\rm Vol}_{k_{y}}; \SPB_{y}\big),
\end{gather*}
i.e., it is the standard Sobolev space on $[0, \infty)_{\sigma} \times Z$ whose sections take values in the bundle $\SPB$ restricted to the boundary over the base point~$y$. Consider
\begin{gather*}
	L(y, \eta) \colon \ \mathcal{H}^{r, \delta, l} \lra \mathcal{H}^{r - 1, \delta, l}.
\end{gather*}

\begin{Lemma}\label{thm:normaloperatorisinvertible}
If the fiber operators $\eth^{Z}_{y}$ satisfy Assumption~{\rm \ref{thm:assumption}} for each $y$, then
\begin{gather*}
	L(y, \eta) \colon \ \mathcal{H}^{r, \delta, l} \lra \mathcal{H}^{r - 1, \delta, l}
\end{gather*}
is invertible for $0 < \delta < 1$, where $L(y, \eta)$ and $\mathcal{H}^{r, \delta, l}$ are defined in~\eqref{eq:11} and~\eqref{eq:Fredholmnormalspaces}.
\end{Lemma}

\begin{proof}[Proof of Lemma \ref{thm:normaloperatorisinvertible}] Given $y \in Y$ and $\eta \in T_{y}Y$ with $\eta \neq 0$, writing $\hat{\eta} = \eta / \absv{\eta}$, we have $(ic(\hat{\eta}))^{2} =\operatorname{id}$. Furthermore,
\begin{gather*}
	\eth^{Z}_{y} ic(\hat{\eta}) = ic(\hat{\eta}) \eth^{Z}_{y} ,
\end{gather*}
so these operators are simultaneously diagonalizable on $L^{2}(Z;\SPB_{0,y}, k_{0})$. Thus for each $y$ and $\hat{\eta}$ we have an orthonormal basis $\set{\phi_{i, \pm}}_{i =1}^{\infty}$ of $L^{2}(Z; \SPB_{0,y}, k_{y})$ satisfying
\begin{gather}\label{eq:etaspectrum}
\eth^{Z}_{y} \phi_{i, \pm} = \pm \mu_{i} \phi_{i, \pm}, \qquad ic(\hat{\eta}) \phi_{i, \pm}	= \pm \phi_{i, \pm}, \qquad c_{\nu}\phi_{i, \pm} = \pm \phi_{i, \mp}.
\end{gather}
Note that the existence of such an orthonormal basis is automatic from the existence of any simultaneous diagonalization $\wt{\phi}_{i}$. Indeed, since $c_{\nu}$ is the operator on the bundle $\SPB_{y}$ induced by~$c(\p_{x})$, we have $i c(\hat{\eta})c_{\nu}\wt{\phi}_{i} = -c_{\nu}\wt{\phi}_{i}$, so we can reindex to obtain $\phi_{i, \pm}$ satisfying the two equations on the right in~\eqref{eq:etaspectrum}. But then since $c_{\nu}\eth^{Z}_{y} = -\eth^{Z}_{y}c_{\nu}$, the f\/irst equation in~\eqref{eq:etaspectrum} follows automatically. Using the $\phi_{i, \pm}$, we def\/ine subspaces of $\mathcal{H}^{r, \delta, l}$ by
\begin{gather*}
	W_{i}^{r, \delta, l} = \spn\big\{ \lp a(\sigma) \phi_{i, +} + b(\sigma)\phi_{i, -} \rp \colon
	a, b \in \mathcal{H}^{r, \delta, l}(d\sigma) \big\},
\end{gather*}
where $\mathcal{H}^{r, \delta, l}(d\sigma)$ is def\/ined as in~\eqref{eq:Fredholmnormalspaces} in the case that $Z$ is a single point. In particular, for all~$\eta$ and~$i$,
\begin{gather*}
	W_{i}^{r, \delta, l} \subset \mathcal{H}^{r, \delta, l}.
\end{gather*}
Note that multiplication by $c_{\nu}$ def\/ines a unitary isomorphism of $ W_{i}^{r, \delta, l}$. We consider the map~$L(y, \eta)$ on each space individually. We claim that
\begin{gather}\label{eq:L0mapping}
	L(y, \eta) \colon \ W_{i}^{r, \delta, l} \lra W_{i}^{r - 1, \delta, l} \ \ \text{is invertible for $0< \delta < 1$.}
\end{gather}
 From \eqref{eq:11}, we compute
\begin{gather}\label{eq:normalactiontransform}
-c_{\nu} \cdot L(y, \eta) a(\sigma) \phi_{i, \pm}= \lp \sigma \p_{\sigma} + \frac f2 \pm \mu \rp a \phi_{i, \pm} -\sigma \absv{\eta} a \phi_{i, \mp}.
\end{gather}
Thus, writing elements in $W_{i}^{r, \delta, l}$ as vector valued functions $(a, b)^{T} = a \phi_{i, +} + b \phi_{i, -}$, we see that $L(y, \eta)$ indeed maps $W_{i}^{r, \delta, l}$ to $W_{i}^{r - 1, \delta, l}$, acting as the matrix
\begin{gather*}
	-c_{\nu} L(y, \eta) \rvert_{W_{i}} = \sigma \p_{\sigma} +f/2 + \lp
	\begin{matrix} \mu & - \sigma\absv{\eta} \\ - \sigma\absv{\eta} & - \mu \end{matrix} \rp.
\end{gather*}
From this, one checks that that the solutions to $L(y, \eta)\phi = 0$ can be written using separation of variables as superpositions of sections given, again in terms of the $\phi_{i, \pm}$ by
\begin{gather*}
	\mathcal{I}_{\mu, \eta}(\sigma) := \sigma^{-f/2 + 1/2} \lp
	\begin{matrix}
	I_{\absv{\mu + 1/2}}(\absv{\eta}\sigma) \\
	I_{\absv{\mu - 1/2}}(\absv{\eta}\sigma)
	\end{matrix}
	\rp,
\qquad
	\mathcal{K}_{\mu, \eta}(\sigma) := \sigma^{-f/2 + 1/2} \lp
	\begin{matrix}
	-K_{\absv{\mu + 1/2}}(\absv{\eta}\sigma) \\
	K_{\absv{\mu - 1/2}}(\absv{\eta}\sigma)
	\end{matrix}	\rp ,
\end{gather*}
where $I_{\nu}(z)$ and $K_{\nu}(z)$ are the modif\/ied Bessel functions~\cite{AS1964}. But neither of these lie in $W^{r,\delta, l}_i$. Indeed, by the asymptotic formulas \cite[equation~(9.7)]{AS1964}, the sections involving the $\mathcal{I}_{\mu,\eta}$ are not even tempered distributions, as the grow exponentially as $z \to \infty$. The sections involving the~$\mathcal{K}_{\mu, \eta}$ are tempered distributions, but since $K_{\nu}(z) \sim z^{-\nu}$ as $z \to 0$ for $\nu > 0$, Assump\-tion~\ref{thm:assumption} tells us that $\mu \not \in (-1/2, 1/2)$, so $\max\{|\mu - 1/2|, |\mu + 1/2|\} \ge 1$. Thus $\mathcal{K}_{\mu, \eta} \not \in \mathcal{H}^{r,\delta, l}$ for any $\delta > 0$. In other words, Assumption~\ref{thm:assumption} implies that
\begin{gather*}
	\text{\eqref{eq:L0mapping} is injective if $\delta > 0$}.
\end{gather*}
On the other hand, the ordinary dif\/ferential operator in~\eqref{eq:normalactiontransform} admits an explicit right inverse if $\delta < 1$. Specif\/ically, consider the matrix
\begin{gather}
\mathcal{M}_{\mu, \absv{\eta}}(\sigma, \wt{\sigma}) = (\sigma\wt{\sigma})^{1/2}\absv{\eta} \lp
	\begin{matrix}
	I_{\absv{\mu + 1/2}}(\absv{\eta}\sigma) & -K_{\absv{\mu + 1/2}}(\absv{\eta}\sigma) \\
	I_{\absv{\mu - 1/2}}(\absv{\eta}\sigma) & K_{\absv{\mu - 1/2}}(\absv{\eta}\sigma)
	\end{matrix}
	\rp \nonumber\\
\hphantom{\mathcal{M}_{\mu, \absv{\eta}}(\sigma, \wt{\sigma}) =}{} \times \lp
	\begin{matrix}
	- H(\wt{\sigma} - \sigma) K_{\absv{\mu - 1/2}}(\absv{\eta}\wt{\sigma}) &
	- H(\wt{\sigma} - \sigma) K_{\absv{\mu + 1/2}}(\absv{\eta}\wt{\sigma}) \\
	- H(\sigma - \wt{\sigma}) I_{\absv{\mu - 1/2}}(\absv{\eta}\wt{\sigma}) &
	H(\sigma - \wt{\sigma}) I_{\absv{\mu + 1/2}}(\absv{\eta}\wt{\sigma})
	\end{matrix}
	\rp .\label{eq:normalgreens}
\end{gather}
Then the operator $Q_{y, \mu}$ on $W^{r - 1, \delta, l}_{i}$ def\/ined by acting on elements $a(\sigma) \phi_{i, +} + b(\sigma)\phi_{i, -} $ by
\begin{gather}\label{eq:normaloperatorseparated}
	Q_{y, \mu} \left(
	\begin{matrix} a \\ b \end{matrix}
	\right) :=
	\sigma^{-f/2}\int_{0}^{\infty} \mathcal{M}_{\mu, \absv{\eta}} (\absv{\eta}, \sigma, \wt{\sigma}) \wt{\sigma}^{f/2 - 1}
	\left(
	\begin{matrix}
	-b(\wt{\sigma}) \\ a(\wt{\sigma})
	\end{matrix}
	\right) d\wt{\sigma}
\end{gather}
satisf\/ies
\begin{gather}\label{eq:10}
	L(y, \eta)Q_{y, \mu} = \operatorname{id}\qquad \text{on}\quad W^{r, \delta, l}_{i}.
\end{gather}
(One checks \eqref{eq:10} using the recurrence relations and Wronskian identity
\begin{gather}
	I_{\nu}'(z) = I_{\nu - 1}(z) - \frac{\nu}{z}I_{\nu}(z) =\frac{\nu}{z}I_{\nu}(z) + I_{\nu + 1}(z),\nonumber \\
	K_{\nu}'(z) = - K_{\nu - 1}(z) - \frac{\nu}{z}K_{\nu}(z) =\frac{\nu}{z}K_{\nu}(z) - K_{\nu + 1}(z),\nonumber \\
	1/z = I_{\nu}(z)K_{\nu + 1}(z) + I_{\nu + 1}(z)K_{\nu}(z),\label{eq:modifiedBesselrecurrencewronsk}
\end{gather}
which are equations (9.6.15) and (9.6.26) from \cite{AS1964}.) That $Q_{y, \mu} \colon W^{r, \delta, l}_{i} \lra W^{r - 1, \delta, l}_{i}$ is bounded for $\delta < 1$ can be seen using~\cite{Ma1991}, but one can also check it directly using the density of polyhomogeneous functions. Invertibility on each $W^{r, \delta, l}_{i}$ gives invertibility on $\mathcal{H}^{r, \delta, l}$. This proves Lem\-ma~\ref{thm:normaloperatorisinvertible}.
\end{proof}

Theorem \ref{thm:edgetheorem} then follows from \cite{Ma1991} as explained in \cite[Section~2]{ALMPII} using the invertibility of the normal operator from
Lemma~\ref{thm:normaloperatorisinvertible}. In the notation of those papers, one has the numbers
\begin{gather*}
	\overline{\delta} := \inf\big\{\delta\colon L(y, \eta) \colon \sigma^{\delta}
	L^{2}(d\sigma d{\rm Vol}_{z};\, \SPB_{y}) \lra L^{2}(d\sigma d{\rm Vol}_{z}; \SPB_{y}) \ \text{is injective for all $y$}\big\},\\
	\underline{\delta} := \sup\big\{\delta \colon L(y, \eta) \colon \sigma^{\delta}
	L^{2}(d\sigma d{\rm Vol}_{z}; \, \SPB_{y}) \lra L^{2}(d\sigma d{\rm Vol}_{z};\, \SPB_{y}) \ \text{is surjective for all $y$}\big\}.
\end{gather*}
By our work above, $\overline{\delta} \le 0 < 1 \le \underline{\delta}$, and thus for the map $x \eth \colon x^{\delta}H^{k}_{\rm e} \lra x^{\delta}H^{k}_{\rm e}$ with $0 < \delta < 1$, there exist $\wt{Q}_{i}$, $i = 1,2$ satisfying~\eqref{eq:edgeparametrix} for $x \eth$. In particular, by~\eqref{eq:normaloperatorequation}
\begin{gather*}
	N(x \eth)_{y}N\big(\wt{Q}_{i}\big)_{y} = N(\operatorname{Id}) = \delta_{\beta^{*}\Delta \cap \ff},
\end{gather*}
where $\beta^{*}\Delta$ is the lift of the interior of the diagonal $\Delta \subset M \times M$ to the blown up space~$M^{2}$. Thus in the coordinates~\eqref{eq:projectivecoords}, $\delta_{\beta^{*}\Delta \cap \ff} = \delta_{\sigma = 1, \YY = 0}$, so from~\eqref{eq:normalgreens} and~\eqref{eq:normaloperatorseparated}, we can write $N(\wt{Q}_{i})_{y}$ as follows. For f\/ixed $\eta$ and the basis $\phi_{i, \pm}$, $i = 1, 2, \dots$, from~\eqref{eq:etaspectrum}, let $\Pi(i, \eta)$ denote $L^{2}$ orthogonal projection onto $\phi_{i, \pm}$ and def\/ine the vectors
\begin{gather}\label{eq:thevector}
	\Pi(\eta, i) = \lp
	\begin{matrix}
	\pi(\eta, i, +) \\
	\pi(\eta, i, -)
	\end{matrix}\rp,
\end{gather}
where $\pi(\eta, i, \pm)$ is orthogonal projection in $L^{2}(Z, \SPB_{0, y}, k_{y})$ onto $\phi_{i, \pm}$. We thus have
\begin{gather}\label{eq:normalparamfinal}
\Pi(\eta, i) \widehat{N(\wt{Q}_{i})_{y}}\Pi^{*}(\eta, i)=(\wt{\sigma}/\sigma)^{f/2}\wt{\sigma}^{-1} \mathcal{M}_{\mu_{i},\absv{\eta}}(\sigma, 1) ,
\end{gather}
where $\Pi^{*}(\eta, i) \lp
\begin{matrix}
 a \\
 b
\end{matrix}
\rp = a \phi_{i, +} + b\phi_{i, -}$.

\subsection[Proof of Theorem \ref{thm:essential} and the generalized inverse of $\eth$]{Proof of Theorem \ref{thm:essential} and the generalized inverse of $\boldsymbol{\eth}$}\label{sec:proofofessential}

In this section we will prove Theorem \ref{thm:essential} and describe the properties of the integral kernel of the generalized inverse of~$\eth$. We start by recalling the statement for the convenience of the reader:

\begin{Theorem} Assume that $\eth$ is a Dirac operator on a spin incomplete edge space $(M,g)$, sa\-tis\-fying Assumption~{\rm \ref{thm:assumption}}, then the unbounded operator $\eth$ on $L^{2}(M; \SPB)$ with core domain $C^{\infty}_{c}(M; \SPB)$ is essentially self-adjoint. Moreover, letting $\mathcal{D}$ denote the domain of this self-adjoint extension, the map
\begin{gather*}
 \eth \colon \ \mathcal{D} \lra L^{2}(M; \SPB)
\end{gather*}
is Fredholm.
\end{Theorem}

\begin{proof}
The proof of Theorem \ref{thm:essential} will follow from combining various elements of \cite{ALMPI, ALMPII}. The f\/irst and main step is the construction of a left parametrix for the map $\eth \colon \mathcal{D}_{\max} \lra L^{2}(M; \SPB)$, where $\mathcal{D}_{\max}$ is the maximal domain def\/ined in \eqref{eq:maxandmin}.

Consider $\wt{Q}_{1}$ from \eqref{eq:edgeparametrix} and set $\wt{Q}_{1}x = \overline{Q}_{1}$. Then by \eqref{eq:edgeparametrix}
\begin{gather}\label{eq:incedgeparametrix1}
 \overline{Q}_{1} \eth = {\operatorname{Id}} - \Pi_{\operatorname{ker}, \delta},
\end{gather}
where both sides of this equation are thought of as maps of $x^{\delta}L^{2}_{\rm e}(M ; \SPB)$. We claim that in fact equation~\eqref{eq:incedgeparametrix1} holds not only on $x^{\delta}L^{2}_{\rm e}(M ; \SPB)$, but on the maximal domain $\mathcal{D}_{\max}$ def\/ined in~\eqref{eq:maxandmin}. This follows from \cite[Lemma~2.7]{ALMPII} as follows. In the notation of that paper, $L = \eth$ and $P = x \eth$. Taking (again, in the notation of that paper) $\mathcal{E}(L)$ to be $\mathcal{D}_{\max}$, by \cite[Lemma~2.1]{ALMPII}, $\mathcal{E}^{(\tau)}(L) = \mathcal{E}(L)$. Furthermore, $\mathcal{E}_{\tau}(L) = x^{\tau}L^{2}(M; \SPB) \cap \mathcal{D}_{\max}$. Since $\wt{Q}_{1}$ maps $x^{\delta}L^{2}(M; \SPB)$ to $x^{\delta}H^{1}_{\rm e}$, we have ${\operatorname{Id}} - \overline{Q}_{1} \eth$ is bounded on~$\mathcal{E}_{\tau}$. Futhermore, $x\eth$ maps $\mathcal{D}_{\max}$ to $xL^{2}(M; \SPB) \subset x^{\delta}L^{2}(M; \SPB)$, so $\overline{Q}_{1} \eth = \wt{Q}_{1} x\eth$ maps $\mathcal{D}_{\max}$ to $x^{\delta}L^{2}(M; \SPB)$. Thus \cite[Lemma~2.7]{ALMPII} applies and~\eqref{eq:incedgeparametrix1} holds on~$\mathcal{D}_{\max}$, as advertised.

Thus $\operatorname{Id} = \overline{Q}_{1} \eth + \Pi_{\operatorname{ker}, \delta}$ on $\mathcal{D}_{\max}$, and since the right hand side is bounded $L^{2}(M; \SPB)$ to $x^{\delta}L^{2}$, for any $\delta \in (0, 1)$, we have
\begin{gather}\label{eq:extravanishing}
	\mathcal{D}_{\max} \subset \bigcap_{\delta < 1} x^{\delta}L^{2}(M, \SPB),
\end{gather}
in particular for any $\delta < 1$, $\mathcal{D}_{\max} \subset H^{1}_{\rm loc} \cap x^{\delta}L^{2}(M, \SPB)$ which is a compact subset of $L^{2}(M; \SPB)$. It then follows from Gil--Mendoza \cite{GM2003} (see \cite[Proposition~5.11]{ALMPI}) that $\mathcal{D}_{\max} \subset \mathcal{D}_{\min}$, i.e., that~$\eth$ is essentially self-adjoint. By a standard argument, e.g.,~\cite[Lemma~4.2]{MelMicro}, the fact that $\mathcal{D}_{\max}$ includes compactly into $L^{2}(M, \SPB)$ implies that $\eth$ has f\/inite-dimensional kernel and closed range. But the self-adjointness of~$\eth$ on $\mathcal{D}$ now implies that~$\eth$ has f\/inite-dimensional cokernel, so~$\eth$ is self-adjoint and Fredholm.
\end{proof}

Thus $\eth$ admits a generalized inverse $Q \colon L^{2}(M; \SPB) \lra \mathcal{D}$ satisfying
\begin{gather*}
	\eth Q = {\operatorname{Id}} - \pi_{\operatorname{ker}} \qquad \text{and} \qquad Q = Q^{*},
\end{gather*}
where $\pi_{\operatorname{ker}} $ is $L^{2}$ orthogonal projection onto the kernel of $\eth$ in $\mathcal{D}$ with respect to the pairing induced by the Hermitian inner product on~$\SPB$. To be precise, if $\set{\phi_{i}}$, $i = 1, \dots, N$ is an orthonormal basis for the kernel of $\eth$ on $\mathcal{D}$, then $\pi_{\operatorname{ker}}$ has Schwartz kernel
\begin{gather*}
	K_{\pi_{\operatorname{ker}}}(w,w') = \sum_{i = 1}^{N} \phi_{i}(w) \otimes \bar{\phi_{i}(w')}.
\end{gather*}
From \eqref{eq:extravanishing}, we see that $\pi_{\operatorname{ker}} \in \pdo^{-\infty}_{\rm e}(M; \SPB; a, b)$. The properties of the integral kernel of~$Q$ can be deduced from those of the parametrices $\wt{Q}_{i}$ in~\eqref{eq:edgeparametrix}. Indeed, setting $\wt{Q} = Q x^{-1}$, we see that
\begin{gather*}
	\wt{Q} (x \eth) = {\operatorname{Id}} - \pi_{\operatorname{ker}} \qquad \text{and} \qquad (x\eth) \wt{Q} = {\operatorname{Id}} - x \cdot
	\pi_{\operatorname{ker}} \cdot x^{-1}.
\end{gather*}
Applying the argument from \cite[Section~4]{Ma1991}, specif\/ically equations~(4.24) and~(4.25) there, shows that $\wt{Q} \in \pdo^{-1}_{\rm e}(M; \SPB;a, b)$ for the same~$a$,~$b$ as in~\eqref{eq:edgeparametrix}, and in particular that $N(\wt{Q}) = N(\wt{Q}_{i}$). In particular, by Theorem~\ref{thm:edgetheorem}, we have the bounds
\begin{gather*}
	K_{Q}(p) = O\big(\rho_{\lf}^{a}\big) \qquad \text{as} \quad p \to \lf \qquad \text{and} \qquad K_{Q}(p) = O\big(\rho_{\rf}^{b}\big) \qquad \text{as} \quad p \to \rf,
\end{gather*}
where $a > \delta - f/2 - 1/2$ and $b > - \delta - f/2 + 1/2$, $0 < \delta < 1$, and again the bounds hold for $K_{Q}$ as a~section of $\End(\SPB)$ over $M \times M$. Finally, by self-adjointness of~$Q$, we have that
\begin{gather}\label{eq:selfadjointinverse}
	K_{Q}(w, w') = K_{Q}^{*}(w', w) \qquad \text{for all} \quad w, w' \in M^{\rm int}.
\end{gather}
By \eqref{eq:selfadjointinverse}, the bound at $\rf$, which one approaches in particular if~$w$ remains f\/ixed in the interior of~$M$ and~$w'$ goes to the boundary, gives a bound at~$\lf$. Thus we obtain the following.

\begin{Proposition}\label{thm:Qproperties}
The distributional section $K_{Q}$ of $\End(\SPB)$ over $M \times M$ with the property that $Q \phi = \int_{M} K_{Q}(w, w') \phi'(w') d{\rm Vol}_{g}(w')$ is conormal at $\Delta_{\rm e}$, and $\rho_{\ff}^{f - 1} K_{Q}$ is smoothly conormal up to $\ff$, where
\begin{gather*}
	\rho_{\ff}^{f} K_{Q} x^{-1} \rvert_{\ff} = \rho_{\ff}^{f} \wt{Q} \rvert_{\ff}
\end{gather*}
satisfies \eqref{eq:normalparamfinal}. Moreover, for coordinates $(x, y, z, x', y', z')$ on $M \times M$ as in~\eqref{eq:coordinates},
\begin{gather}\label{eq:vanishingrate}
K_{Q}(x, y, z, x', y', z') = O(x^{a}), \qquad \text{uniformly for} \quad x' \ge c > 0,
\end{gather}
where $a > - \delta - f/2 + 1/2$ for any $\delta > 0$ and $c$ is an arbitrary small positive number.
\end{Proposition}

\section{Boundary values and boundary value projectors} \label{sec:boundaryvalues}
Recall that $M_{\epsilon} = \{ x \ge \epsilon \}$ is a smooth manifold with boundary, and $M - M_{\epsilon}$ is a~tubular neighborhood of the singularity.
Consider the space of harmonic sections over $M - M_{\epsilon}$
\begin{gather*}
	\mathcal{H}_{{\rm loc}, \epsilon} = \big\{ u \in L^{2}(M - M_{\epsilon}; \SPB) \colon \eth u = 0,\, \exists \, \wt{u} \in \mathcal{D} \ \text{s.t.} \ u = \wt{u} \rvert_{M - M_{\epsilon}} \big\},
\end{gather*}
where $\mathcal{D}$ is the domain for $\eth$ from Theorem~\ref{thm:essential}; in particular, $\mathcal{D} \subset H^{1}_{\rm loc}$. By the standard restriction theorem for $H^{1}$ sections \cite[Proposition~4.5, Chapter~4]{TaylorI}, any element $u \in \mathcal{H}_{{\rm loc}, \epsilon}$ has boundary values $u \rvert_{\p M_{\epsilon}} \in H^{1/2}(\p M_{\epsilon})$. We def\/ine a~domain for $\eth$ on the cutof\/f manifold~$M_{\epsilon}$ by
\begin{gather}\label{eq:epsilondomain}
	\mathcal{D}_{\epsilon} := \big\{ u \in H^{1}(M_{\epsilon}; \SPB) \colon u
\rvert_{\p M_{\epsilon}} = v \rvert_{\p M_{\epsilon}} \ \text{for some} \ v \in \mathcal{H}_{{\rm loc}, \epsilon}\big\} \subset L^{2}(M_{\epsilon}; \SPB).
\end{gather}
Essentially, $\mathcal{D}_{\epsilon}$ consists of sections over $M_{\epsilon}$ whose boundary values correspond with the boundary values of an $L^{2}$ harmonic section over $M - M_{\epsilon}$. We also have the chirality spaces
\begin{gather*}
	\mathcal{D}_{\epsilon}^{\pm} = \mathcal{D}_{\epsilon} \cap L^{2}\big(M_{\epsilon}; \SPB^{\pm}\big),
\end{gather*}
where $\SPB^{\pm}$ are the chirality subbundles of even and odd spinors. In this section we will prove the following.
\begin{Theorem}\label{thm:naturalbvp}
For $\epsilon > 0$ sufficiently small and $\mathcal{D}_{\epsilon}$ as in~\eqref{eq:epsilondomain}, the map $\eth \colon \mathcal{D}_{\epsilon} \lra L^{2}(M_{\epsilon}; \SPB)$ is Fredholm, and
\begin{gather*}
	\ind\big(\eth \colon \mathcal{D}_{\epsilon}^{+} \lra L^{2}(M_{\epsilon}; \SPB^{-})\big) = \ind\big(\eth \colon \mathcal{D}^{+}
	\lra L^{2}(M; \SPB^{-})\big).
\end{gather*}
\end{Theorem}

In the process of proving Theorem \ref{thm:naturalbvp}, we will construct a family of boundary value projec\-tors~$\pi_{\epsilon}$ which def\/ine~$\mathcal{D}_{\epsilon}$ in the sense of Claim~\ref{thm:domainfromprojection} below, and whose microlocal structure we will use in Section~\ref{sec:equivalence} to relate the index of~$\eth$ on $M_{\epsilon}$ with domain~$\mathcal{D}_{\epsilon}$ to the index of $\eth$ on $M_{\epsilon}$ with the APS boundary condition, see Theorem~\ref{thm:APSmatchesother}.

\subsection[Boundary value projector for $\mathcal{D}_{\epsilon}$]{Boundary value projector for $\boldsymbol{\mathcal{D}_{\epsilon}}$}

As already mentioned, the main tool for proving Theorem \ref{thm:naturalbvp} and also for proving Theorem~\ref{thm:APSmatchesother} below will be to
express the boundary condition in the def\/inition of $\mathcal{D}_{\epsilon}$ in~\eqref{eq:epsilondomain} in terms of a~pseudodif\/ferential projection over~$\p M_{\epsilon}$. We discuss the construction of this projection now.

First we claim that the invertible double construction of \cite[Chapter~9]{BB1993} holds in this context in the following form: there exists an incomplete edge manifold $M'$ with spinor bundle $\SPB'$ and Dirac operator $\eth'$, together with an isomorphism
\begin{gather*}
	\Phi \colon \ (M - M_{\epsilon_{0}}, \SPB) \lra (M' - M'_{\epsilon_{0}}, \SPB')
\end{gather*}
such that, with identif\/ications induced by $\Phi$, the operators $\eth$ and $\eth'$ are equal over $M - M_{\epsilon_{0}}$ $(= \Phi^{-1}(M' - M'_{\epsilon_{0}}))$, and f\/inally \textit{such that $\eth'$ is invertible}. In particular, the inverse $Q'$ satisf\/ies
\begin{gather}\label{eq:weirdQ}
	Q' \colon \ \mathcal{D}' \lra L^{2}(M'), \qquad \eth Q' = \operatorname{id} = Q' \eth,
\end{gather}
where $\mathcal{D}'$ is the unique self-adjoint domain for $\eth'$ on $M'$ with core domain $C_{c}^{\infty}(M', \SPB')$. Moreover, $Q'$ satisf\/ies all of the properties in Proposition~\ref{thm:Qproperties}.

We describe the construction of this ``invertible double'' for the convenience of the reader, though it is essentially identical to that in \cite[Chapter~9]{BB1993}, the only dif\/ference being that we must introduce a product type boundary while they have one to begin with. Choosing any point $p \in M_{\epsilon_{0}}$, let $D_{1}$, $D_{2}$ denote open discs around $p$ with $p \in D_{1} \Subset D_{2}$ and $D_{2} \cap (M - M_{\epsilon_{0}}) = \varnothing$. We can identify the annulus $D_{2} - D_{1}$ with $[1, 2)_{s} \times \mathbb{S}^{d-1}$ by a dif\/feomorphism and the metric~$g$ is homotopic to a product metric $ds^{2} + \absv{dx}^{2}$ where $x$ is the standard coordinate on~$\mathbb{B}^{d-1}$. Furthermore, the connection can be deformed so that the induced Dirac operator $\eth'$ is of product type on the annulus (see equation~(9.4) in~\cite{BB1993}). Call the bundle over $N_{1} : = M - D_{1}$ thus obtained~$\wt{\SPB}$. Letting $N_{2} := - N_{1}$, the same incomplete edge space with the opposite orientation, let $M' = N_{1} \sqcup N_{2} /\set{s = 1}$ and consider the vector bundle $\SPB'$ over~$M'$ obtained by taking $\wt{\SPB}^{+}$ over $N_{1}$ and $\wt{\SPB}^{-}$ over~$N_{2}$ and identifying the two bundles over~$D_{2}$ using Clif\/ford multiplication by~$\p_{s}$. The resulting Dirac operator, which we still denote by~$\eth'$, is seen to be invertible on~$M'$ by the symmetry and unique continuation argument in Lemma~9.2 of \cite{BB1993}.

We will now work on a neighborhood in $M - M_{\epsilon}$ of $\p M$ (or equivalently of the singular stratum~$Y$), so we drop the distinction between~$M$ and~$M'$. Using notation as in~\eqref{eq:weirdQ}, and given $f \in C^{\infty}(\p M_{\epsilon}; \SPB)$, def\/ine the harmonic extension
\begin{gather}\label{eq:harmonicextension}
	\ext_{\epsilon} f(w) := \int_{w' \in \p M_{\epsilon}} K_{Q'}(w,w')\cl_{\nu} f(w') d{\rm Vol}_{\p M_{\epsilon}},
\end{gather}
where $K_{Q'}$ is the Schwartz kernel of $Q'$ (see~\eqref{eq:schwartzkernelacts}), and $\cl_{\nu} = \cl(\p_{x})$. Since
\begin{gather}\label{eq:kernelisharmonic}
	\eth' K_{Q'}(w, w') = 0 \qquad \text{away from} \quad w = w',
\end{gather}
$\eth' \ext_{\epsilon}f(w) = 0$ for $w \not \in \p M_{\epsilon}$. Recall Green's formula for Dirac operators; specif\/ically, for a~smoothly bounded region $\Omega$ with normal vector $\p_{\nu}$,
\begin{gather}\label{eq:byparts}
	\int_{\Omega} \lp \la \eth u, v \ra - \la u, \eth v \ra \rp d{\rm Vol}_{\Omega} =\int_{\p \Omega} \la \cl(\p_{\nu}) u , v \ra d{\rm Vol}_{\p \Omega}.
\end{gather}
Green's formula for sections $u$ satisfying $\eth u \equiv 0$ in $M - M_{\epsilon}$ gives that for $u \in \mathcal{H}_{{\rm loc}, \epsilon}$,
\begin{gather}\label{eq:greens}
	u(w) = - \int_{\p M_{\epsilon}} K_{Q'}(w,w') \cl_{\nu}u(w') d{\rm Vol}_{\p M_{\epsilon}},
\end{gather}
\emph{provided}
\begin{gather}\label{eq:annoyingthing}
	\forall \, u \in \mathcal{H}_{{\rm loc}, \epsilon}, \quad w \in M^{\rm int}, \qquad \lim_{\wt{\epsilon} \to 0} \int_{\p M_{\wt{\epsilon}}} K_{Q'}(w,w') \cl_{\nu} u(w') d{\rm Vol}_{\p M_{\wt{\epsilon}}} = 0.
\end{gather}
The identity in \eqref{eq:greens} is obtained by integrating by parts in
\begin{gather*}
\int_{M - M_{\epsilon}} \eth K_{Q'}(w,w') u(w') - K_{Q'}(w,w') \eth u(w') d{\rm Vol}_{w'}
\end{gather*}
and using \eqref{eq:kernelisharmonic}. In fact, as we will see in the proof of Claim~\ref{thm:domainfromprojection} below, \eqref{eq:annoyingthing}, and thus~\eqref{eq:greens}, hold for all $u \in \mathcal{H}_{{\rm loc}, \epsilon}$.

It follows from \eqref{eq:harmonicextension} and \eqref{eq:greens} that, for $\eth u = 0$ satisfying \eqref{eq:annoyingthing},
\begin{gather}\label{eq:layerpotential}
	u \rvert_{\p M_{\epsilon}}(w) = \mathcal{E}_{\epsilon}(u \rvert_{\p	M_{\epsilon}})(w),
\end{gather}
where
\begin{gather*}
	\mathcal{E}_{\epsilon} (f)(w) := \lim_{\substack{ \wt w \to w \\ \wt{w} \in M - M_{\epsilon} } } \ext _{\epsilon} (f)(w)
\end{gather*}
We will show that the $\mathcal{E}_{\epsilon}$ def\/ine the domains $\mathcal{D}_{\epsilon}$ as follows.

\begin{Claim}\label{thm:domainfromprojection}
The operator $\mathcal{E}_{\epsilon}$ in \eqref{eq:layerpotential} is a projection operator on $L^{2}(\p M_{\epsilon}, \SPB)$, and the do\-main~$\mathcal{D}_{\epsilon}$ in~\eqref{eq:epsilondomain} is given by
\begin{gather}\label{eq:cutoffdomain}
	\mathcal{D}_{\epsilon} = \big\{ u \in H^{1}(M_{\epsilon}; \SPB) \colon ({\operatorname{Id}} - \mathcal{E}_{\epsilon})(u \rvert_{\p M_{\epsilon}}) = 0 \big\}.
\end{gather}
Moreover, there exists $\Nop \in \pdo^{0}(\p M_{\epsilon}; \SPB)$ such that
\begin{gather}
\mathcal{E}_{\epsilon} = \frac 12 {\operatorname{Id}} + \Nop,\nonumber \\
\Nop f(w) = - \int_{\p M_{\epsilon}} K_{Q'}(w,w') \cl_{\nu} f(w') d{\rm Vol}_{\p
M_{\epsilon}} \qquad \text{for} \quad w \in \p M_{\epsilon},\label{eq:offdiagcomplete}
\end{gather}
and the principal symbol of $\mathcal{E}_{\epsilon}$ satisfies
\begin{gather}\label{eq:principalofE}
	\sigma(\mathcal{E}_{\epsilon}) = \sigma(\pi_{{\rm APS}, \epsilon}),
\end{gather}
where $\pi_{{\rm APS}, \epsilon}$ is the APS projection defined in \eqref{eq:APSprojector}.
\end{Claim}

Assuming Claim \ref{thm:domainfromprojection} for the moment, we prove Theorem~\ref{thm:naturalbvp}.

\begin{proof}[Proof of Theorem~\ref{thm:naturalbvp} assuming Claim~\ref{thm:domainfromprojection}] The main use of Claim~\ref{thm:domainfromprojection} in this context (it will be used again in Theorem \ref{thm:APSmatchesother} ) is to show that the map
\begin{gather}\label{eq:cutoffoperator}
	\eth \colon \ \mathcal{D}_{\epsilon} \lra L^{2}(M_{\epsilon}).
\end{gather}
is self-adjoint on $L^{2}(M_{\epsilon}; \SPB)$ and Fredholm. The Fredholm property follows from the principal symbol equality \eqref{eq:principalofE}, since from~\cite{BB1993} any projection in $\pdo^{0}(\p M_{\epsilon}, \SPB)$ with principal symbol equal to that of the Atiyah--Patodi--Singer boundary projection def\/ines a Fredholm problem. To see that it is self-adjoint, note that from \eqref{eq:byparts} the adjoint boundary condition is
\begin{gather*}
 \mathcal{D}_{\epsilon}^{*} = \set{ \phi \colon \la g , \cl_{\nu} \phi \rvert_{\p M_{\epsilon}} \ra_{\p M_{\epsilon}} = 0 \ \text{for all} \ g \ \text{with} \ ({\operatorname{Id}} - \mathcal{E}_{\epsilon}) g = 0 }.
\end{gather*}
Again by \eqref{eq:byparts}, for any $v \in \mathcal{H}_{{\rm loc}, \epsilon}$, $v \rvert_{\p M_{\epsilon}} \in \mathcal{D}_{\epsilon}^{*}$. Thus $\mathcal{D}_{\epsilon} \subset \mathcal{D}_{\epsilon}^{*}$, and it remains to show that $\mathcal{D}_{\epsilon}^{*} \subset \mathcal{D}_{\epsilon}$. Let $\phi \in \mathcal{D}_{\epsilon}^{*}$, and set $f := \phi \rvert_{\p M_{\epsilon}}$. We want to show that $(I - \mathcal{E}_{\epsilon}) f = 0$, or equivalently
\begin{gather}\label{eq:fdeductions}
	\la (I - \mathcal{E}_{\epsilon}) f , g \ra_{\p M_{\epsilon}} = 0 \quad \forall\, g \quad \iff \quad
\la f , (I - \mathcal{E}_{\epsilon}^{*}) g \ra_{\p M_{\epsilon}} 	= 0 \quad \forall\, g.
\end{gather}
Since $\la f , \cl_{\nu}g \ra = - \la \cl_{\nu} f, g \ra$, by~\eqref{eq:byparts} we have $\la f , \cl_{\nu} g \ra = 0$ for every $g \in \Ran \mathcal{E}$, and thus~\eqref{eq:fdeductions} will hold if $(I - \mathcal{E}_{\epsilon}^{*}) g \in \Ran \cl_{\nu} \mathcal{E}$. In fact, we claim that
\begin{gather*}
	I - \mathcal{E}_{\epsilon}^{*} = - \cl_{\nu} \mathcal{E} \cl_{\nu}.
\end{gather*}
To see that his holds, note that by Claim \ref{thm:domainfromprojection} and self-adjointness of~$Q'$, specif\/ically~\eqref{eq:selfadjointinverse}, $\Nop^{*} = \cl_\nu \Nop \cl_\nu$, so
\begin{gather*}
I - \mathcal{E}_{\epsilon}^{*} = I - \left(\frac 12 + \Nop\right)^{*} = \frac{1}{2} - \Nop^{*} = - \cl_\nu \left( \frac 12 + \Nop \right) \cl_\nu	= - \cl_\nu \mathcal{E}_{\epsilon} \cl_\nu,
\end{gather*}
which proves self-adjointness.

Now that we know that \eqref{eq:cutoffoperator} is self-adjoint, we proceed as follows. We claim that for $\epsilon > 0$ suf\/f\/iciently small, the map
\begin{align*}
\ker\big(\eth \colon \mathcal{D} \lra L^{2}(M; \SPB)\big) & \lra \ker\big(\eth \colon \mathcal{D}_{\epsilon} \lra L^{2}(M_{\epsilon}; \SPB)\big), \\
\wt{\phi} &\longmapsto \phi = \wt{\phi} \rvert_{M_{\epsilon}}
\end{align*}
is well def\/ined and an isomorphism. It is well def\/ined since by def\/inition any section $\wt{\phi} \in \ker(\eth \colon \mathcal{D} \lra L^{2}(M; \SPB))$ satisf\/ies that $\phi = \wt{\phi} \rvert_{M_{\epsilon}} \in \mathcal{D}_{\epsilon}$. It is injective by unique continuation. For surjectivity, note that for any element $\phi \in \ker(\eth \colon \mathcal{D}_{\epsilon} \lra L^{2}(M_{\epsilon}; \SPB))$, by def\/inition there is a $u \in \mathcal{H}_{{\rm loc}, \epsilon}$ such that $u \rvert_{\p M_{\epsilon}} = \phi \rvert_{\p M_{\epsilon}}$. It follows that
\begin{gather*}
	\wt{\phi}(w) :=
	\begin{cases}
	\phi(w) & \text{for} \ w \in M_{\epsilon}, \\
	u(w) & \text{for} \ w \in M - M_{\epsilon}
	\end{cases}
\end{gather*}
is in $H^{1}$ and satisf\/ies $\eth \wt{\phi} = 0$ on all of~$M$, i.e., $\wt{\phi} \in \ker(\eth \colon \mathcal{D} \lra L^{2}(M; \SPB))$. Since the full operator~$\eth$ on~$\mathcal{D}$ is self-adjoint, and since the operator in~\eqref{eq:cutoffoperator} is also, the cokernels of both maps are equal to the respective kernels. Restricting $\eth$ to a~map from sections of $\SPB^{+}$ to sections of $\SPB^{-}$ gives the theorem.

This completes the proof.\end{proof}

Thus to prove Theorem \ref{thm:naturalbvp} it remains to prove Claim~\ref{thm:domainfromprojection}.
\begin{proof}[Proof of Claim \ref{thm:domainfromprojection}]
We begin by proving \eqref{eq:offdiagcomplete}. It is a standard fact (see \cite[Section~7.11]{taylor:vol2}) that
\begin{gather*}
	\mathcal{E}_{\epsilon} \in \pdo^{0}(\p M_{\epsilon}; \SPB).
\end{gather*}
Obviously,
\begin{gather}\label{eq:offdiag}
	\mathcal{E}_{\epsilon} = A + \Nop,
\end{gather}
where
\begin{gather*}
\Nop f(w) = - \int_{\p M_{\epsilon}} K_{Q'}(w,w') \cl_{\nu} f(w') d{\rm Vol}_{\p
	M_{\epsilon}} \qquad \text{for} \quad w \in \p M_{\epsilon},\\
	\supp A \subset \diag(\p M_{\epsilon} \times \p M_{\epsilon}),
\end{gather*}
where the last containment refers to the Schwartz kernel of $A$. We claim that
\begin{gather*}
	A = \frac{1}{2} \operatorname{id} \qquad \text{and} \qquad B_{\epsilon}	\in \Psi^{0}(\p M_{\epsilon}; \SPB).
\end{gather*}
Using that $Q'$ has principal symbol $\sigma(Q') = \sigma(\eth)^{-1}$ we can write $Q'$ in local coordinates $w$ as
\begin{gather*}
	Q' =\frac{1}{(2\pi)^{n}} \int_{\mathbb{R}^{n}} e^{- (w - \wt{w}) \cdot \xi} a(w, \wt{w}, \xi) d\xi \quad \text{locally},
\end{gather*}
 where
\begin{gather*}
	a(w, \wt{w}, \xi) = \absv{\xi}_{g(w)}^{-2} i c(\xi) + O(1) \qquad \text{for} \quad \absv{\xi} \ge c > 0.
\end{gather*}
Given a bump function $\chi$ supported near $w_{0} \in \p M_{\epsilon}$, let $Q'_{\chi} := \chi Q' \chi$ and def\/ine the distributions
\begin{gather}\label{eq:localizedQ}
	K_{Q'_{\chi}} = K_{1} + K_{2},
\end{gather}
where
\begin{gather*}
K_{1} = \mathcal{F}_{\xi}^{-1}\big( \absv{\xi}_{g(w)}^{-2} i c(\xi) \big) \mathcal{F}_{\wt{x}},
\end{gather*}
where, as in \eqref{eq:schwartzkernelacts}, $K_{Q'_{\chi}}$ denotes the Schwartz kernel of $Q'_{\chi}$. The distribution $K_{2}$ is that of a~pseudodif\/ferential operator of order~$-2$ on~$M$, and it follows from the theory of homogeneous distributions (see \cite[Chapter~7]{taylor:vol2}) that the distribution~$K_{2}$ restricts to $\p M_{\epsilon}$ to be the Schwartz kernel of a~pseudodif\/ferential operator of order $-1$. The distribution $K_{1}$ is that of a~pseudodif\/ferential operator on~$M$ of order $-1$. It is smooth in $\wt{x}$ with values in homogeneous distributions in $x - \wt{x}$ of order $-n + 1$, and it follows that the restriction of the Schwartz kernel $K_{1}(w, w')$ to $\p M_{\epsilon}$ gives a~pseudodif\/ferential operator of order zero. Letting $\Nop$ in \eqref{eq:offdiag} be the operator def\/ined by the restriction of $K_{1}$ to $\p M_{\epsilon}$, we have that $\Nop$ is in~$\pdo^{0}(\p M_{\epsilon}; \SPB)$ and it remains to calculate $A$. Choosing coordinates of the form $w = (x, x')$ and $\wt{w} = (\wt{x}, \wt{x}')$ of the form in~\eqref{eq:symbolofDirac} and such that at the f\/ixed value $w_{0} = (\epsilon, x'_{0}) \in
\p M_{\epsilon}$ the metric satisf\/ies $g(x) = \operatorname{id}$, it follows (see~\cite{AHM2006}) that the Schwartz kernel of~$B$ in~\eqref{eq:localizedQ} satisf\/ies
\begin{gather*}
	B(x_{0},\wt{x}) = - \frac{1}{\omega_{n-1}} \frac{c(x_{0}) - c(\wt{x})}{\absv{x_{0} - \wt{x}}^{n}} + O\big(\absv{x_{0} - \wt{x}}^{2 - n}\big),
\end{gather*}
where $\omega_{n-1}$ is the volume of the unit sphere $\mathbb{S}^{n-1}$. If we let $\wt{B}(x', \wt{x}') = B(0, x', 0 ,\wt{x}')$, then near~$x_{0}$
\begin{gather*}
	\ext_{\epsilon} f(\delta, x') = - \frac{1}{\omega_{n-1}} \int
	\lp \frac{c((\delta, x')) - c((0, \wt{x}'))}{\absv{(\delta, x') -
	(0, \wt{x}')}^{n}} \rp \cl_{\nu} f(\wt{x}') d\wt{x}' \\
\hphantom{\ext_{\epsilon} f(\delta, x')}{} = - \frac{1}{\omega_{n-1}} \int
	\lp \frac{\delta \cl_{\nu}}{\absv{(\delta, x') -
	(0, \wt{x}')}^{n}} + \frac{c((0, x')) - c((0, \wt{x}'))}{\absv{(\delta, x') -
	(0, \wt{x}')}^{n}} \rp \cl_{\nu} f(\wt{x}')
	d\wt{x}' \\
\hphantom{\ext_{\epsilon} f(\delta, x')}{} \to \frac{1}{2} f(x') + \int \wt{B}(x', y') \cl_{\nu}f(y') dy' \qquad \text{as} \quad \delta \to 0.
\end{gather*}
This proves that $A = 1/2 $.

The principal symbol of $\mathcal{E}_{\epsilon}$ (again see \cite[Section~7.11]{taylor:vol2}) is given by the integral
\begin{gather*}
	\sigma(\mathcal{E})(x', \xi') = \frac{- 1}{2\pi}\lim_{x \to \epsilon^{-}}\int_{\mathbb{R}}
	e^{i(x - \epsilon) \xi} \frac{1}{|(\xi, \xi')|_{g}^{2}} (ic(\xi
	\p_{x}) + ic(\xi' \cdot \p_{x'})) c(\p_{x}) d\xi \\
\hphantom{\sigma(\mathcal{E})(x', \xi')}{} = \frac{-1}{2\pi}
	\lim_{x \to \epsilon^{-}} \lp \int_{\mathbb{R}}
	e^{i(x - \epsilon) \xi} \frac{\xi}{|(\xi, \xi')|_{g}^{2}}
	d\xi \rp ic(\p_{x})^{2} - \frac 12 ic(\hat{\xi}' \cdot \p_{x'})	c(\p_{x}) \\
\hphantom{\sigma(\mathcal{E})(x', \xi')}{} = \frac{-1}{2\pi}
	\lp - \frac{2\pi i}{2}\rp ic(\p_{x})^{2} - \frac 12 (- i c(\p_{x} c(\hat{\xi}' \cdot \p_{x'}))
	\\
\hphantom{\sigma(\mathcal{E})(x', \xi')}{}= \frac 12 - \frac 12 (-i c(\p_{x}) c(\hat{\xi}' \cdot \p_{x'}))).
\end{gather*}
where in the third line we used the residue theorem. Now recall that the term $-i c(\p_{x}) c(\hat{\xi}' \cdot \p_{x'}))$ is precisely the endomorphism appearing in~\eqref{eq:boundarysymbolasprojections}, so
\begin{gather*}
	\sigma(\mathcal{E})(x', \xi') = \frac 12 - \frac 12 \big( \pi_{\epsilon, + , \hat{\xi}'}(x') - \pi_{\epsilon, - , \hat{\xi}'}(x')\big) = \pi_{\epsilon, - ,\hat{\xi}'}(x'),
\end{gather*}
where the projections are those in~\eqref{eq:boundarysymbolasprojections}. Thus, Theorem~\ref{thm:APSbasics} implies the desired formula for the principal symbol of $\mathcal{E}_{\epsilon}$, \eqref{eq:principalofE}.

To f\/inish the claim, we must show the equivalence of domains in~\eqref{eq:cutoffdomain}. We f\/irst show that for any $u \in \mathcal{H}_{{\rm loc}, \epsilon}$, the formula in~\eqref{eq:greens} holds. This will show that any $f \in H^{1/2}(\p M_{\epsilon}; \SPB)$ with $f = u \rvert_{\p M_{\epsilon}}$ for some $u \in \mathcal{H}_{{\rm loc}, \epsilon}$ satisf\/ies $({\operatorname{Id}} - \mathcal{E}_{\epsilon}) f = 0$, i.e., that
\begin{gather*}
\mathcal{D}_{\epsilon} \subset \big\{ u \in H^{1}(M_{\epsilon}; \SPB) \colon ({\operatorname{Id}} -	\mathcal{E}_{\epsilon})(u \rvert_{\p M_{\epsilon}}) = 0 \big\}.
\end{gather*}
Thus we must show that \eqref{eq:annoyingthing} holds for $u \in \mathcal{H}_{{\rm loc}, \epsilon}$. For such $u$, we claim that for some $\delta > 0$, as $\epsilon \to 0$
\begin{gather}\label{eq:decayofl2bdrynorm}
	\int_{\p M_{\epsilon}} \norm{u}^{2} d{\rm Vol}_{\p M_{\epsilon}} = O\big(\epsilon^{- f - \delta}\big).
\end{gather}
To see this, note f\/irst that $u \in x^{1 - \delta}H^{1}_{\rm e}(M - M_{\epsilon}, \SPB)$ for every $\delta > 0$, which follows since $u$ has an extension to a section in $\mathcal{D}_{\max} \subset H^{1}_{\rm loc} \cap_{\delta > 0} x^{1 - \delta}L^{2}(M; \SPB)$. In particular, $x^{\delta + f/2}u \in H^{1}(M, dxdydz)$, the standard Sobolev space of order $1$ on the manifold with boundary $M$. Using the restriction theorem \cite[Proposition~4.5, Chapter~4]{TaylorI}, $x^{\delta + f/2}u = \epsilon^{\delta + f/2}u \in H^{1/2}(\p M)$ uniformly in $\epsilon$, so
\eqref{eq:decayofl2bdrynorm} holds. Thus, for f\/ixed $w \in M - M_{\epsilon}$, writing $d{\rm Vol}_{g} = x^{f} a dx dy dz$ for some $a = a(x, y, z)$ with $a(0, y, z) \neq 0$, we can use the bound for $K_{Q'}$ in~\eqref{eq:vanishingrate} with $x'$ f\/ixed and $x = \epsilon$ to conclude
\begin{gather*}
	\lp \int_{\p M_{\epsilon}} K_{Q'}(w,w') \cl_{\nu} u(w') d{\rm Vol}_{\p M_{\epsilon}} \rp^{2} = \lp \int_{\p M_{\epsilon}}
	K_{Q'}(w,w') \cl_{\nu} u(w') \epsilon^{f} a dy'dz' \rp^{2} \\
\qquad {} \le \epsilon^{2f} \left(\int_{\p M_{\epsilon}} \norm{K_{Q'}(w,w')}^{2} a dy'dz'\right)
\left(\int_{\p M_{\epsilon}}\norm{u}^{2} ady'dz' \right) \\
\qquad {} = \epsilon^{2f}o\big(\epsilon^{-2\delta - f + 1}\big)
	o\big(\epsilon^{- f - 2\delta}\big) \qquad \text{for all} \quad \delta > 0 \\
\qquad {}= o\big(\epsilon^{- 4\delta + 1}\big) \to 0 \qquad \text{as} \quad \epsilon \to 0.
\end{gather*}
To prove the other direction of containment in~\eqref{eq:cutoffdomain}, we need to know that for $f \in H^{1/2}(\p M_{\epsilon})$ satisfying $({\operatorname{Id}} - \mathcal{E}_{\epsilon})f = 0$, the section $u := \ext_{\epsilon}f \rvert_{M - M_{\epsilon}} \in \mathcal{H}_{{\rm loc}, \epsilon}$, where $\ext_{\epsilon}$ is the extension operator in~\eqref{eq:harmonicextension}. This is true since for any $H^{1/2}$ section $h$ over $\p M_{\epsilon}$, there is an $H^{1}$ extension~$v$ to the manifold~$M'$ def\/ined above, that can be taken with support away from the singular locus. If $1_{M'_{\epsilon}}$ is the indicator function of $M'_{\epsilon}$, then $\eth'(\ext_{\epsilon} f + 1_{M'_{\epsilon}} v) = \delta_{\p M_{\epsilon}}(f + h) + 1_{M'_{\epsilon}} \eth' v$. Taking $h$ to cancel $f$ gives that $\eth'(\ext_{\epsilon} f + 1_{M'_{\epsilon}} v) \in L^{2}(M'; \SPB)$. Since $1_{M'_{\epsilon}} v$ is an extendible $H^{1}$ distribution on $M_{\epsilon}'$ near $\p M_{\epsilon}$, $\ext f \rvert_{M - M_{\epsilon}}$ is an extendible $H_{\rm loc}^{1}$ distribution on $M - M_{\epsilon}$ near~$\p M_{\epsilon}$. This completes the proof of Claim~\ref{thm:domainfromprojection}.
\end{proof}

\section{Equivalence of indices}\label{sec:equivalence}
In the previous section we have shown
\begin{gather*}
\ind\big(\eth \colon \mathcal{D}_{\epsilon}^{+} \lra L^{2}(M_{\epsilon}; \SPB^{-})\big) = \ind\big(\eth \colon \mathcal{D}^{+}\lra L^{2}(M; \SPB^{-})\big),
\end{gather*}
in this section we will prove the following.
\begin{Theorem}\label{thm:APSmatchesother}
Let $\pi_{{\rm APS}, \epsilon}$ denote the APS projector from \eqref{eq:APSproblem}. Then for $\epsilon > 0$ sufficiently small,
\begin{gather*}
\ind\big(\eth \colon \mathcal{D}^{+}_{{\rm APS}, \epsilon} \lra L^{2}(M_{\epsilon}; \SPB^{-})\big) = \ind\big(\eth \colon \mathcal{D}^{+}_{\epsilon} \lra L^{2}(M_{\epsilon}; \SPB^{-})\big),
\end{gather*}
where $\mathcal{D}_{\epsilon}$ is the domain in \eqref{eq:cutoffdomain} and $\mathcal{D}_{{\rm APS}, \epsilon}$ is the domain in~\eqref{eq:APSproblem}.
\end{Theorem}

The main tool for proving Theorem \ref{thm:APSmatchesother} is the following theorem from~\cite{BB1993}. We def\/ine the `pseudodif\/ferential Grassmanians'
\begin{gather}\label{eq:grassman}
	\Gr_{{\rm APS}, \epsilon} = \big\{ \pi \in \Psi^{0}(\p M_{\epsilon}; \SPB) \colon
	\pi^{2} = \pi \ \text{and} \ \sigma(\pi) = \sigma(\pi_{{\rm APS}, \epsilon})\big\}.
\end{gather}
We endow $\Gr_{{\rm APS}, \epsilon}$ with the norm topology. If $\pi \in\Gr_{{\rm APS}, \epsilon}$, then def\/ining the domain $\mathcal{D}_{\pi, \epsilon} =\set{ u \in H^{1}(M_{\epsilon}; \SPB) \colon ({\operatorname{Id}} - \pi) (u \rvert_{\p M_{\epsilon}}) = 0}$, the map
\begin{gather*}
	\eth \colon \ \mathcal{D}_{\pi, \epsilon} \lra L^{2}(M_{\epsilon}; \SPB)
\end{gather*}
is Fredholm. The following follows from \cite[Theorem~20.8]{BB1993} and \cite[Theorem~15.12]{BB1993}
\begin{Theorem}\label{thm:samecomponent}
If $\pi_{i} \in \Gr_{{\rm APS}, \epsilon}$, $i = 1,2$ lie in the same connected component of $\Gr_{{\rm APS}, \epsilon}$ then the elliptic boundary problems
\begin{gather*}
	\eth \colon \ \mathcal{D}^{+}_{\pi_{i}, \epsilon} \lra L^{2}(M_{\epsilon}, \SPB^{-})
\end{gather*}
have equal indices.
\end{Theorem}

To apply Theorem \ref{thm:samecomponent} in our case, we will study the two families of boundary values projec\-tors~$\pi_{{\rm APS}, \epsilon}$ and $\mathcal{E}_{\epsilon}$ using the adiabatic calculus of Mazzeo and Melrose~\cite{MaMe1990}.

\subsection{Review of the adiabatic calculus}
Consider a f\/iber bundle $Z \hookrightarrow \wt{X} \xlra{\pi} Y$. The adiabatic double space $\wt{X}^{2}_{\rm ad}$ is formed by radial blow up of $\wt{X}^{2} \times [0, \epsilon_{0})_{\epsilon}$ along the f\/iber diagonal, $\diag_{\rm f\/ib}(\wt{X})$ (see~\eqref{eq:fiberdiagonal}) at $\epsilon = 0$. That is,
\begin{gather*}
	\wt{X}^{2}_{\rm ad} = \big[\wt{X}^{2} \times [0, \epsilon_{0})_{\epsilon} ;	\diag_{\rm f\/ib}\big(\wt X\big) \times \set{\epsilon = 0}\big].
\end{gather*}
Thus, $\wt{X}^{2}_{\rm ad}$ is a manifold with corners with two boundary hypersurfaces: the lift of $\set{\epsilon = 0}$, which we continue to denote by $\set{\epsilon = 0}$, and the one introduced by the blowup, which we call~$\ff$. Similar to the edge front face above, $\ff$ is a bundle over~$Y$ whose f\/ibers are isomorphic to $Z^{2} \times \mathbb{R}^{b}$ where $b = \dim Y$, and in fact this is the f\/iber product of $\pi^*T^*Y$ and~$\wt X$.
\begin{gather*}
	\text{We def\/ine $\ff_{y}$ to be the f\/iber of $\ff$ lying above $y$}.
\end{gather*}

 \begin{figure} \centering
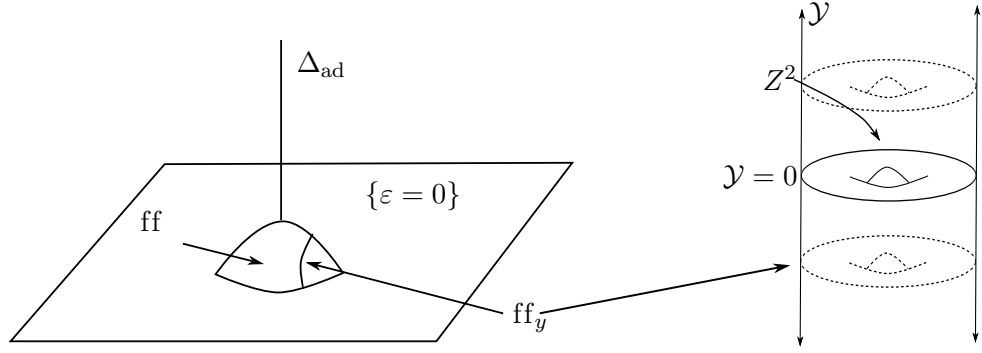
 \caption{The adiabatic double space.}\label{fig:adiabatic}
 \end{figure}

The adiabatic vector f\/ields on the f\/ibration $\wt{X}$ are families of vector f\/ields $V_{\epsilon}$ parametrized smoothly in $\epsilon \in [0, \epsilon_0)$, such that $V_{0}$ is a vertical vector f\/ield, i.e., a section of $T\wt X/Y$. Locally these are $C^{\infty}(\wt{X} \times [0, \epsilon_{0})_{\epsilon})$ linear combinations of the vector f\/ields
\begin{gather*}
	\p_{z}, \quad \epsilon \p_{y}.
\end{gather*}
Such families of vector f\/ields are in fact sections of a vector bundle
\begin{gather}\label{eq:adiabatictangentbundle}
	T_{\rm ad}(\p M) \lra \p M \times [0, \epsilon_{0})_{\epsilon}.
\end{gather}
We will now def\/ine adiabatic dif\/ferential operators on sections of $\SPB$. The space of $m^{\rm th}$ order adiabatic dif\/ferential operators $\Diff^{m}_{\rm ad}(\wt{X}; \SPB)$ is the space of dif\/ferential operators obtained by taking $C^{\infty}(\wt{X}; \End(\SPB))$ combinations of powers (up to order $m$) of adiabatic vector f\/ields. An adiabatic dif\/ferential operator $P$ admits a normal operator $N(P)$, obtained by letting~$P$ act on $\wt{X} \times \wt{X} \times [0, \epsilon_{0})_{\epsilon}$, pulling back $P$ to $\wt{X}^{2}_{\rm ad}$, and restricing it to~$\ff$. The normal operator acts tangentially along the f\/ibers of~$\ff$ over~$Y$, and $N(P)_{y}$ will denote the operator on sections of~$\SPB$ restricted to over~$\ff_{y}$. More concretely, if~$P$ is an adiabatic operator of order~$m$, then near a~point~$y_{0}$ in~$Y$, we can write
\begin{gather*}
	P = \sum_{\absv{\alpha} + \absv{\beta} \le m} a_{\alpha, \beta}(z,y,\epsilon)
	\p_{z}^{\alpha}(\epsilon \p_{y})^{\beta},
\end{gather*}
for $y$ near $y_{0}$, where $a_{\alpha, \beta}(z,y,\epsilon)$ is a smooth family of endomorphisms of $\SPB$. The normal operator is given by
\begin{gather*}
	N(P)_{y_{0}} = \sum_{\absv{\alpha} + \absv{\beta} \le m} a_{\alpha, \beta}(z,y_{0},0) \p_{z}^{\alpha}\p_{\YY}^{\beta},
\end{gather*}
where $\YY$ are coordinates on $\mathbb{R}^{b}$. Thus $N(P)_{y_{0}}$ is a dif\/ferential operator on $Z_{y_0} \times T_{y_0}Y$ that is constant coef\/f\/icient in the $TY$ direction.

Returning to the case that $\wt{X} = \p M$ for $M$ a compact manifold with boundary, we take a~collar neighborhood $\mathcal{U} \simeq \p M \times [0, \epsilon_{0})_{x}$ as in \eqref{eq:collarneighborhood}, and treating the boundary def\/ining function,~$x$, as the parameter $\epsilon$ in the previous paragraph, identify the adiabatic double space $(\p M)^{2}_{\rm ad}$ with a~blow up of $\set{x = x'} \subset \mathcal{U} \times \mathcal{U}$.

\begin{Lemma}\label{thm:adiabaticdiff}
The tangential operator $\wt{\eth}_{\epsilon}$ defined in \eqref{eq:boundaryoperator} lies in $\Diff^{m}_{\rm ad}(\p M; \SPB)$. The normal operator of $\wt{\eth}_{\epsilon}$ satisfies
\begin{gather*}
	N\big(\wt{\eth}_{\epsilon}\big)_{y} = \eth^{Z}_{y} - c_{\nu} \eth_{\YY},
\end{gather*}
where $\eth^{Z}_{y}$ is as in \eqref{eq:fiberop}, and $\eth_{\YY}$ is the standard Dirac operator on $T_{y}Y$.
\end{Lemma}
\begin{proof}
 This follows from equation \eqref{eq:normaloperatorxeth} above.
\end{proof}

The space of \textit{adiabatic pseudodifferential operators} with bounds on $\wt{X}$ of order $m$ acting on sections of $\SPB$, denoted $\Psi^{m}_{{\rm ad}, {\rm bnd}}(\wt{X}; \SPB)$, is the space of families of pseudodif\/ferential operators $\set{A_{\epsilon}}_{0 < \epsilon < \epsilon_{0}}$, where $A_{\epsilon}$ is a (standard) $\pdo$ of order $m$ for each $\epsilon$, and whose integral kernel of $A_{\epsilon}$ is conormal to the lifted diagonal $\Delta_{\rm ad} := \overline{\diag_{\wt{X}} \times (0, \epsilon_{0})_{\epsilon}}$, smoothly up to $\ff$. To be precise, the Schwartz kernel of an operator $A \in \pdo^{m}(\p M; \SPB)$ is given by a family of Schwartz kernels $K_{A_{\epsilon}} = K_{1, \epsilon} + K_{2, \epsilon}$ where $K_{1, \epsilon}$ is conormal of order~$m$ at $\Delta_{\rm ad}$ smoothly down to~$\ff$ and supported near $\Delta_{\rm ad}$, and $K_{2, \epsilon}$ is smooth on the interior and bounded at the boundary hypersurfaces.

An adiabatic pseudodif\/ferential operator $A \in \pdo^{m}(\p M; \SPB)$ with bounds comes with two crucial pieces of data: a principal symbol and a normal operator. The principal symbol $\sigma(A)(\epsilon)$ is the standard one def\/ined for a~conormal distribution, i.e., as a homogeneous section of $N^{*}(\Delta; \End(\SPB)) \otimes \Omega^{1/2}$, the conormal bundle to the lifted diagonal (with coef\/f\/icients in half-densities). In our case $N^{*}(\Delta)$ is canonically isomorphic to $T^{*}_{\rm ad}(\p M)$, the dual bundle to $T_{\rm ad}(\p M)$ def\/ined in~\eqref{eq:adiabatictangentbundle}; in particular, the symbol of $A$ is a map $\sigma(A) \colon T^{*}_{\rm ad}(\p M) \lra C^{\infty}(M; \End(\SPB))$, well def\/ined only to leading order, and smooth down to $\ff$. The normal operator is the restriction of the Schwartz kernel of $A$ to the front face
\begin{gather*}
	N(A) = K_{A}\rvert_{\ff}.
\end{gather*}
We thus have maps
\begin{gather*}
\pdo^{m-1}_{\rm ad}(\p M; \SPB) \hookrightarrow \pdo^{m}_{\rm ad}(\p M; \SPB)\xlra{\sigma} S^{m}(T^{*}_{\rm ad}(\p M; \SPB)) \otimes \Omega^{1/2},
\end{gather*}
and
\begin{gather*}
	N \colon \ \pdo^{m}_{\rm ad}(\p M; \SPB) \lra \pdo^{m}_{\ff, {\rm ad}}(\p M; \SPB).
\end{gather*}

For f\/ixed $\eta$, we use the same eigenvectors $\phi_{i, \pm}$, $\eth^{Z}_{y}$ and $ic(\hat{\eta})$ as in~\eqref{eq:etaspectrum} above, and consider the spaces
\begin{gather}\label{eq:realspacesnormalfaceadiabatic}
	\mathcal{W}_{i} = \spn\set{ \phi_{i, +} , \phi_{i, -} }.
\end{gather}
We have the following.
 \begin{Lemma}\label{thm:layerpotentialadiabatic}
 The layer potential $\mathcal{E}_{\epsilon}$ is a zero-th order adiabatic family $($with bounds$)$, i.e., $\mathcal{E}_{\epsilon} \in \Psi^{0}_{\rm ad}(\p M; \SPB)$. Using the vectors $\Pi(\eta, i)$ from \eqref{eq:thevector}, the normal symbol of~$\mathcal{E}_{\epsilon}$ satisfies
\begin{gather}\label{eq:layerpotentialnormalsymbol}
\Pi(\eta, i) \widehat{N_{y}(\mathcal{E}_{\epsilon})} (y, \eta)\Pi^{*}(\eta, i)= \mathcal{N}_{\mu_{i}, \absv{\eta}},
\end{gather}
where
\begin{gather}\label{eq:Nmatrixcalderon}
	\mathcal{N}_{ \mu, \absv{\eta}} = \absv{\eta} \lp
	\begin{matrix}
	I_{\absv{\mu + 1/2}}(\absv{\eta}) K_{\absv{\mu - 1/2}}(\absv{\eta})
	& I_{\absv{\mu + 1/2}}(\absv{\eta}) K_{\absv{\mu + 1/2}}(\absv{\eta}) \\
	I_{\absv{\mu - 1/2}}(\absv{\eta}) K_{\absv{\mu - 1/2}}(\absv{\eta})
	& I_{\absv{\mu - 1/2}}(\absv{\eta}) K_{\absv{\mu + 1/2}}(\absv{\eta})
	\end{matrix}
	\rp
\end{gather}
and $I$, $K$ 
denote modified Bessel functions.
\end{Lemma}

\begin{proof} That $\mathcal{E}_{\epsilon}$ is an adiabatic pseudodif\/ferential operator follows from \eqref{eq:offdiagcomplete}. The formula in~\eqref{eq:layerpotentialnormalsymbol},~\eqref{eq:Nmatrixcalderon} follows from the Fourier decomposition of the normal operator of the generalized inverse $Q$ in Proposition \ref{thm:Qproperties}, since by~\eqref{eq:harmonicextension} the operator $\mathcal{E}_{\epsilon}$ is obtained by taking the limit in~\eqref{eq:normalparamfinal} as $\sigma = x/x' \uparrow 1$ and checking that $\lim\limits_{\sigma \uparrow 1}\mathcal{M}_{\mu,\absv{\eta}}(\sigma, 1) = \mathcal{N}_{ \mu, \absv{\eta}}$.
 \end{proof}

\subsection{APS projections as an adiabatic family}
To study the integral kernel of the projector $\pi_{{\rm APS}, \epsilon}$ we will make use of the fact that the boundary Dirac operator~$\wt{\eth}_{\epsilon}$ from \eqref{eq:boundaryoperator} is invertible for small $\epsilon$. This is a general fact about adiabatic pseudodif\/ferential operators: invertibility at $\epsilon=0$ implies invertibility for small epsilon, or formally
\begin{Theorem}\label{thm:adiabaticinverse}
Let $A_{\epsilon} \in \pdo^{m}_{\rm ad}(M; \SPB)$ and assume that on each fiber $\ff_{y}$ the Fourier transform of the normal operator $\widehat{N(A_{\epsilon})}_{y}$ is invertible on $L^{2}(Z;\SPB_{y}, k_{y})$, with $\SPB_{y}$ the restriction of the spinor bundle to the fiber over $y$ and $k_{y} =g_{N/Y}\big|_{y}$. Then $A_{\epsilon}$ is invertible for small~$\epsilon$.
 \end{Theorem}

It is well known \cite{APSI} that for each f\/ixed $\epsilon > 0$, $\pi_{{\rm APS}, \epsilon}$ is a pseudodif\/ferential operator of order~$0$. As $\epsilon$ varies, these operators form an adiabatic family:

\begin{Lemma} The family $\pi_{{\rm APS}, \epsilon}$ lies in $\Psi^{0}_{\rm ad}(\p M; \SPB)$. Its normal symbol $N(\pi_{{\rm APS}, \epsilon})$ satisfies
\begin{gather}\label{eq:APSnormalsymbol}
N(\pi_{{\rm APS}, \epsilon}) = \frac 12 N\big(\wt{\eth}_{\epsilon}\big)^{-1} \big(N\big(\wt{\eth}_{\epsilon}\big) -
\big|N\big(\wt{\eth}_{\epsilon}\big)\big| \big).
\end{gather}
 \end{Lemma}
 \begin{proof}
By Assumption \ref{thm:assumption}, $\wt{\eth}_{\epsilon}$ is invertible for small $\epsilon$. Indeed, by~\eqref{eq:adiabaticnormalDiracmatrix}, $N (\wt{\eth}_{\epsilon})_{y}$ does not have zero as an eigenvalue. The projectors $\pi_{{\rm APS}, \epsilon}$ can be expressed in terms of functions of the tangential operators $\wt{\eth}_{\epsilon}$ \cite{APSI} via the formula
\begin{gather}\label{eq:APSasfunction}
	\pi_{{\rm APS}, \epsilon} = \frac 12 \wt{\eth}_{\epsilon}^{-1} \big(\wt{\eth}_{\epsilon} - \big|\wt{\eth}_{\epsilon}\big| \big).
\end{gather}
Following \cite{S1966}, the operator $\wt{\eth}_{\epsilon}^{-1}|\wt{\eth}_{\epsilon}|$ is in $\Psi_{\rm ad}^{1}(\p X; \SPB)$ and has the expected normal operator, namely the one obtained by applying the appropriate functions to the normal operator of $N(\wt{\eth}_{\epsilon})_{y}$ and composing them.
\end{proof}

We compute that the operator $\widehat{N(\wt{\eth}_{\epsilon})}_{y}$ acts on the spaces $\mathcal{W}_{i}$ from~\eqref{eq:realspacesnormalfaceadiabatic} by
\begin{gather*}
	\widehat{N(\wt{\eth}_{\epsilon})}_{y}(\eta) \phi_{i, \pm} = \pm \mu \phi_{i, \pm} - \absv{\eta} \phi_{i, \mp}.
\end{gather*}
That is to say, with $\Pi(\eta, i)$ as in \eqref{eq:thevector},
\begin{gather}\label{eq:adiabaticnormalDiracmatrix}
\Pi(\eta, i)\widehat{N(\wt{\eth}_{\epsilon})}_{y}(\eta)\Pi^{*}(\eta, i) =
\lp	\begin{matrix}
	\mu & - \absv{\eta} \\
	- \absv{\eta} & - \mu
	\end{matrix} \rp.
\end{gather}
Thus,
\begin{gather*}
\Pi(\eta, i)\widehat{N(\wt{\eth}_{\epsilon}^{-1}|\wt{\eth}_{\epsilon}|)}_{y}(\eta) \Pi^{*}(\eta, i)=\frac{1}{(\mu^{2} + \absv{\eta}^{2})^{1/2}}
\lp	\begin{matrix}
\mu & - \absv{\eta} \\
- \absv{\eta} & - \mu
\end{matrix}\rp	\Pi_{\mu_{i}, j}.
\end{gather*}
Using \eqref{eq:APSasfunction}, we obtain
\begin{gather*}
	\Pi(\eta, i) \widehat{N}(\pi_{{\rm APS}, \epsilon})_{y}(\eta) \Pi^{*}(\eta, i) = \mathcal{N}^{\rm APS}_{\mu, \absv{\eta}},
\end{gather*}
where
\begin{gather*}
\mathcal{N}^{\rm APS}_{\mu, \absv{\eta}} = 	\frac 12 \lp \operatorname{Id}_{2 \times 2} + \frac{1}{(\mu^{2} +
	\absv{\eta}^{2})^{1/2}}\lp
	\begin{matrix}
	- \mu & \absv{\eta} \\
	\absv{\eta} & \mu
	\end{matrix}
	\rp\rp.
\end{gather*}

\begin{Theorem}\label{thm:homotopy}
There exists a smooth family $\pi_{\epsilon,t}$, parametrized by $t \in [0, 1]$ satisfying:
\begin{enumerate}\itemsep=0pt
\item[$1)$] for fixed $t$, $\pi_{\epsilon, t} \in \Gr_{{\rm APS}, \epsilon}$, the Grassmanians defined in \eqref{eq:grassman}, and
\item[$2)$] $\pi_{\epsilon, 0} = \mathcal{E}_{\epsilon}$ and $\pi_{\epsilon, 1} = \pi_{{\rm APS}, \epsilon}$.
\end{enumerate}
\end{Theorem}

\begin{proof} The proof proceeds in two main steps. First, we construct a homotopy from the normal operators $N(\mathcal{E}_{\epsilon})$ to
$N(\pi_{{\rm APS},\epsilon})$. Then we extend this homotopy to a homotopy of the adiabatic families as claimed in the theorem.

For the homotopy of the normal operators, the main lemma will be the following
\begin{Claim}\label{thm:normalwithinhalf}
For each $y \in Y$, the normal operators $N(\mathcal{E}_{\epsilon})_{y}$ and $N(\pi_{{\rm APS},\epsilon})$, acting on $L^{2}(Z \times T_{y}Y; \SPB_{y})$, satisfy
\begin{gather*}
\norm[L^{2} \lra L^{2}]{N(\mathcal{E}_{\epsilon})_{y} -	N(\pi_{{\rm APS},\epsilon})} < 1 - \delta,
\end{gather*}
for some $\delta > 0$ independent of $y$.
\end{Claim}

Assuming the claim for the moment, the following argument from \cite[Chapter~15]{BB1993} furnishes a homotopy. In general, let~$P$ and~$Q$ be projections on a separable Hilbert space. Def\/ine $T_{t} = {\operatorname{Id}} + t(Q - P)(2P - {\operatorname{Id}})$, and note that $T_{1}P = QT_{1}$. Now assume that $T_{t}$ is invertible for all~$t$. Then the operator
\begin{gather*}
	F_{t} = T_{t}^{-1} P T_{t}
\end{gather*}
is a homotopy from $P$ to $Q$, i.e., $F_{0} = P$, $F_{1} = Q$. This holds in particular if $\norm{P - Q} < 1$, in which case $T_{t}$ is invertible by Neumann series for $t \in [0, 1]$.

To apply this in our context, we f\/irst take $P = N(\mathcal{E}_{\epsilon})_{y}$ and $Q = N(\pi_{{\rm APS}, \epsilon})_{y}$, and see that the corresponding operator $T_{t}$ is invertible by Claim~\ref{thm:normalwithinhalf}. Now taking $P =\mathcal{E}_{\epsilon}$ and $Q = \pi_{{\rm APS}, \epsilon}$ (so $P$, $Q$, and $T_{t}$ depend on $\epsilon$) by Theorem~\ref{thm:adiabaticinverse}, $T_{t}$ is invertible for small~$\epsilon$. Thus the homotopy $F_{t} = F_{t}(\epsilon) =: \pi_{\epsilon, t}$ is well def\/ined for small~$\epsilon$. In fact, $\pi_{\epsilon, t}$ is a smooth family of adiabatic pseudodif\/ferential projections with principal symbol equal to that of $\pi_{{\rm APS}, \epsilon}$ for all~$\epsilon$.

Thus it remains to prove Claim \ref{thm:normalwithinhalf}. By the formulas for the normal operators given in~\eqref{eq:layerpotentialnormalsymbol} and~\eqref{eq:APSnormalsymbol} and Plancherel, the claim will follow if we can show that for each $\mu$ with $\absv{\mu} > 1/2$, and all $\absv{\eta}$, that
\begin{gather}\label{eq:matrixbound}
	\big\|\mathcal{N}_{\mu, \absv{\eta}} - \mathcal{N}^{\rm APS}_{\mu, \absv{\eta}} \big\| < 1 - \delta,
\end{gather}
for some $\delta$ independent of $\mu \ge 1/2$ and $\absv{\eta}$. Here the norm is as a map of $\mathbb{R}^{2}$ with the standard Euclidean norm.
We prove the bound in~\eqref{eq:matrixbound} using standard bounds on modif\/ied Bessel functions in Appendix~\ref{sec:Bessel}.
\end{proof}

\section{Proof of Main Theorem: limit of the index formula}\label{sec:maintheorem}
Recall (e.g., \cite[Section~2.14]{tapsit}) that if $E\lra M$ is a real vector bundle of rank $k$, with connec\-tion~$\nabla^E$ and curvature tensor~$R^E$ then every smooth function (or formal power series)
\begin{gather*}
	P\colon \ \mathfrak{so}(k)\lra \bbC,
\end{gather*}
that is invariant under the adjoint action of $SO(k)$, determines a closed dif\/ferential form $P(R^E) \in \CI(M;\Lambda^*T^*M)$. If $\nabla^E_1$ is another connection on $E$, with curvature tensor $R^E_1$ then $P(R^E)$ and $P(R^E_1)$ dif\/fer by an exact form. Indeed, def\/ine a family of connections on~$E$ by
\begin{gather*}
	\theta = \nabla^E_1 - \nabla^E \in \CI(M;T^*M \otimes \Hom(E)), \qquad \nabla_t^E = (1-t)\nabla^E + t\nabla^E_1 = \nabla^E + t\theta,
\end{gather*}
denote the curvature of $\nabla^E_t$ by $R^E_t$, and let
\begin{gather*}
	P'(A;B) = \left.\frac{\pa}{\pa s}\right|_{s=0} P(A+sB).
\end{gather*}
The dif\/ferential form
\begin{gather*}
	TP\big(\nabla^E, \nabla^E_1\big) = \int_0^1 P'\big(R^E_t; \theta\big) \, dt
\end{gather*}
satisf\/ies
\begin{gather*}
	dTP\big(\nabla^E, \nabla^E_1\big) = P\big(R^E\big) - P\big(R^E_1\big).
\end{gather*}

Now consider for $\eps<1$ the truncated manifold $M_{\eps} = \{ x \geq \eps \}$ and the corresponding truncated collar neighborhood $\sC_{\eps} = [\eps,1] \times N$. Let $\nabla^{\pt}$ be the Levi-Civita connection of the metric
\begin{gather*}
g_{\pt} = dx^2 + \eps^2 g_Z + \phi^*g_Y.
\end{gather*}
The Atiyah--Patodi--Singer index theorem on $M_{\eps}$ has the form \cite{G1993, Grubb1992}, cf.~\cite{DG2007}
\begin{gather*}
	\int_{M_{\eps}} AS(\nabla) + \int_{\pa M_{\eps}} TAS(\nabla, \nabla^{\pt}) -\tfrac12\eta(\pa M_{\eps})
\end{gather*}
where $AS$ is a characteristic form associated to a connection $\nabla$ and $TAS(\nabla, \nabla^{\pt})$ is its transgression form with respect to the connection $\nabla^{\pt}$.

The Levi-Civita connection of $g_{\pt}$ induces a connection on $T_{\ie}M_{\eps}$, which we continue to denote $\nabla^{\pt}$. Let $\theta^{\eps} = \nabla - \nabla^{\pt}$. Since $g_{\ie}$ and $g_{\pt}$ coincide on $\{ x = \eps \}$ we have
\begin{gather*}
	g_{\pt}\big(\nabla^{\pt}_AB,C\big)\rest{x=\eps} = g_{\ie}\big(\nabla^{\ie}_AB,C\big)\rest{x=\eps} \qquad \text{if}\quad A,B,C \in \CI(\sC_{\eps};TN).
\end{gather*}
On the other hand, if $A,B,C \in \{ \pa_x, \tfrac1xV, \wt U\}$, we have
\begin{gather*}
	g_{\pt}\big(\nabla^{\pt}_{A}B,C\big) = 0 \qquad \text{if}\quad \pa_x \in \{A,B,C\}
\end{gather*}
except for
\begin{gather*}
 g_{\pt}\big(\nabla^{\pt}_{\pa_x}\tfrac1xV_1, \tfrac1x V_2\big) = -\frac{\eps^2}{x^3}g_Z(V_1, V_2).
\end{gather*}
Note that, analogously to \eqref{eq:AsympSplittingConn}, we have
\begin{gather*}
j_0^*\nabla^{\pt} = j_0^*\nabla^{v} \oplus j_0^*\nabla^h,
\end{gather*}
where, as above, $j_{\eps}\colon N \hookrightarrow \sC$ is the inclusion of $\{x=\eps\}$, and $\nabla^v = \bv \circ \nabla \circ \bv$ is the restriction of the Levi-Civita connection to $TN/Y$.

Thus
\begin{gather}\label{eq:DescTheta}
	\theta^{\eps}_A(B)\rest{x=\eps} = 0
\end{gather}
except for
\begin{gather*}
	\theta^{\eps}_{\pa_x}\big(\tfrac1x V\big)\rest{x=\eps} = \tfrac1\eps \tfrac1xV, \qquad
	\theta^{\eps}_V(\pa_x)\rest{x=\eps} = \tfrac1x V, \qquad
	\theta^{\eps}_{V_1}\big(\tfrac1x V_2\big)\rest{x=\eps} = -g_Z(V_1, V_2) \pa_x.
\end{gather*}
In particular note that $j_{\eps}^*\theta^{\eps}$ is independent of $\eps$ and is equal to
\begin{gather*}
j_\eps^*\theta^{\eps} = j_0^*\nabla^{v_+} - j_0^*\nabla^v.
\end{gather*}

Next we need to compute the restriction to $x=\eps$ of the curvature $\Omega_t$ of the connection
$(1-t)\nabla + t\nabla^{\pt} = \nabla + t\theta^{\eps}$.
Locally, with $\omega$ the local connection one-form of $\nabla$ \eqref{eq:Schemeomega}, the curvature $\Omega_t$ is given by
\begin{gather*}
	\Omega_t = d(\omega + t\theta^{\eps}) + (\omega + t\theta^{\eps})\wedge (\omega+t\theta^{\eps})
	= \Omega + t( d\theta^{\eps} + [\omega,\theta^{\eps}]_s) + t^2 \theta^{\eps}\wedge\theta^{\eps},
\end{gather*}
where $[\cdot, \cdot]_s$ denotes the supercommutator with respect to form parity, so that $[\omega,\theta^{\eps}]_s = \omega \wedge \theta^{\eps} + \theta^{\eps} \wedge \omega$. In terms of the splitting \eqref{eq:SecondSplittingIE} we have
\begin{gather*}
	\Omega\rest{x =\eps} = \begin{pmatrix} \Omega_{v_+} & \cO(\eps) \\ \cO(\eps) & \phi^*\Omega_Y \end{pmatrix}, \qquad
	\omega = \begin{pmatrix} \omega_{v_+} & \cO(x) \\ \cO(x) & \phi^*\omega_Y + \cO(x^2) \end{pmatrix}, \qquad
	j_{\eps}^*\theta^{\eps} = \begin{pmatrix} \wt\theta & 0 \\ 0 & 0 \end{pmatrix},
\end{gather*}
and hence
\begin{gather*}
	\Omega_t\rest{x=\eps} =
	\begin{pmatrix} \Omega_{v_+}
		+t(d\wt\theta + [\omega_{v_+}, \wt\theta]_s) + t^2 \wt\theta\wedge\wt\theta
		& \cO(\eps) \\ \cO(\eps) & \phi^*\Omega_Y \end{pmatrix}.
\end{gather*}
In particular, if we denote $\Omega_{v_+,t}$ the curvature of the connection $(1-t)\nabla^{v_+} + t\nabla^v$ on the bundle $\ang{\pa_x}+TN/Y$, we have
\begin{gather*}
j_{\eps}^*\Omega_t = j_0^*\Omega_t + \cO(\eps), \qquad \text{with}\quad j_0^*\Omega_t= \begin{pmatrix} \Omega_{v_+,t} & 0 \\ 0 & \phi^*\Omega_Y \end{pmatrix}.
\end{gather*}
It follows that
\begin{gather*}
	\lim_{\eps\to 0} j_{\eps}^*T\hat A\big(\nabla, \nabla^{\pt}\big)
	= \int_0^1\left.\frac{\pa}{\pa s}\right|_{s=0}j_{0}^*\hat A(\Omega_Y) \hat A\big(\Omega_{v_+,t}+s\wt \theta\big) \, dt
	= \hat A(Y) \wedge T\hat A(\nabla^{v_+},\nabla^v)
\end{gather*}
and similarly for any multiplicative characteristic class.

We can now prove the main theorem, whose statement we recall for the reader's convenience.

\begin{Theorem}Let $X$ be stratified space with a single singular stratum endowed with an incomplete edge metric $g$ and let $M$ be its resolution. If $\eth$ is a Dirac operator associated to a spin bundle $\SPB\lra M$ and $\eth$ satisfies Assumption~{\rm \ref{thm:assumption}}, then
\begin{gather*}
	\ind\big(\eth \colon \mathcal{D}^{+} \lra L^{2}(M; \SPB^{-})\big)
	= \int_{M}\! \hat{A}(M) + \int_Y \! \hat A(Y)\lrpar{{-}\frac 12\hat\eta(\eth_Z) + \int_{\pa M/Y}\! T\hat A\big(\nabla^{v_+},\nabla^{\pt}\big) },
\end{gather*}
where $\hat A$ denotes the $\hat A$-genus, $T\hat A(\nabla^{v_+}, \nabla^{\pt})$ denotes the transgression form of the $\hat A$ genus associated to the connections $\nabla^{v_+}$ and $\nabla^{\pt}$ above, and $\hat\eta$ the $\eta$-form of Bismut--Cheeger~{\rm \cite{BC1989}}.
\end{Theorem}

\begin{proof}[Proof of Main Theorem] Combining Theorems \ref{thm:naturalbvp} and \ref{thm:APSmatchesother} we know that, for $\eps$ small enough,
\begin{gather*}
	\ind\big(\eth \colon \mathcal{D}^{+} \lra L^{2}(M; \SPB^{-})\big)
	= \ind\big(\eth \colon \mathcal{D}_{\epsilon}^{+} \lra L^{2}(M_{\epsilon}; \SPB^{-})\big) \\
\hphantom{\ind\big(\eth \colon \mathcal{D}^{+} \lra L^{2}(M; \SPB^{-})\big)}{}
	= \ind\big(\eth \colon \mathcal{D}^{+}_{{\rm APS}, \epsilon} \lra L^{2}(M_{\epsilon}; \SPB^{-})\big).
\end{gather*}
Hence
\begin{gather*}
\ind\big(\eth \colon \mathcal{D}^{+} \lra L^{2}(M; \SPB^{-})\big)
= \lim_{\eps \to 0}	\ind\big(\eth \colon \mathcal{D}^{+}_{{\rm APS}, \epsilon} \lra L^{2}(M_{\epsilon}; \SPB^{-})\big) \\
\hphantom{\ind\big(\eth \colon \mathcal{D}^{+} \lra L^{2}(M; \SPB^{-})\big)}{}
	= \lim_{\eps \to 0}	\int_{M_{\eps}} \hat A(\nabla) + \int_{\pa M_{\eps}} T\hat A\big(\nabla, \nabla^{\pt}\big) -\tfrac12\eta(\pa M_{\eps}) \\
\hphantom{\ind\big(\eth \colon \mathcal{D}^{+} \lra L^{2}(M; \SPB^{-})\big)}{}
	=\int_{M}\! \hat A(M) + \!\int_{\pa M}\! \hat A(Y) \wedge T\hat A(\nabla^{v_+},\nabla^v)
	-\tfrac12\!\int_Y \!\hat A(Y) \hat\eta(\eth_Z).\!\!\!\!\!\!\!\!\tag*{\qed}
\end{gather*}
\renewcommand{\qed}{}
\end{proof}

\subsection{Four-dimensions with circle f\/ibers}

An incomplete edge space whose link is a sphere is topologically a smooth space. So let us consider a four-dimensional manifold $X$ with a~submanifold $Y$ and a Riemannian metric on $X\setminus Y$ that in a tubular neighborhood of~$Y$ takes the form
\begin{gather*}
	dx^2+x^2\beta^2 d\theta^2 + \phi^*g_Y.
\end{gather*}
Here $\beta$ is a constant and $2\pi\beta$ is the `cone angle' along the edge.

Recall that the circle has two distinct spin structures, and with the round metric the corresponding Dirac operators have spectra equal to either the even or odd integer multiples of~$\pi$. The non-trivial spin structure on the circle is the one that extends to the disk, and so any spin structure on~$X$ will induce non-trivial spin structures on its link circles. Thus, cf.~\cite[Proposition~2.1]{Ch1989}, the generalized Witt assumption will be satisf\/ied as long as~$\beta\leq 1$.

In this setting the relevant characteristic class is the f\/irst Pontryagin class: for a two-by-two anti-symmetric matrix~$A$, let
\begin{gather*}
	p_1(A) = -c_2(A) = -\frac1{8\pi^2}\Tr\big(A^2\big).
\end{gather*}
Note that $p_1'(A;B) = -\frac1{(2\pi)^2} \Tr(AB)$, and so
\begin{gather*}
	Tp_1\big(\nabla, \nabla^{\pt}\big) = -\frac1{(2\pi)^2} \int_0^1 \Tr j_0^*(\theta \wedge \Omega_t) \, dt
\end{gather*}
with $\Omega_t = \Omega + t( d\theta^{\eps} + [\omega,\theta^{\eps}]_s) + t^2 \theta^{\eps}\wedge\theta^{\eps}$. We can simplify this formula. Indeed, note that if~$\{V_i\}$ are an orthonormal frame for $TN/Y$ then
\begin{gather*}
j_{\eps}^*(\theta^{\eps}\wedge\theta^{\eps}) = \sum \Theta_{ij} V_i^{\flat} \wedge V_j^{\flat}
	\qquad \text{with} \quad \Theta_{ij}\big(\tfrac1x V_k\big) = -\delta_{kj}\tfrac1x V_i,
\end{gather*}
and so in particular $\dim Z = 1$ implies $j_{\eps}^*(\theta^{\eps}\wedge\theta^{\eps})=0$. Moreover with respect to the splitting~\eqref{eq:SplittingIE}, $\theta$ is of\/f-diagonal and $\Omega$ is on-diagonal, hence $\Tr j_0^*(\theta \wedge \Omega)=0$ and
\begin{gather*}
Tp_1(\nabla, \nabla^{\pt}) 	= -\frac 1{4\pi^2} \int_0^1 t \Tr j_0^*(\theta \wedge d\theta) \, dt = -\frac 1{8\pi^2} \Tr j_0^*(\theta \wedge d\theta).
\end{gather*}
Next let us consider $\theta$ in more detail. From \eqref{eq:DescTheta}, with respect to the splitting \eqref{eq:SplittingIE}, we have
\begin{gather*}
	j_0^*\theta = \begin{pmatrix} 0&\Id & 0 \\ -\Id & 0 &0 \\ 0 & 0 & 0 \end{pmatrix} \alpha,
\end{gather*}
where $\alpha$ is a vertical one-form of $g_{Z}$ length one. This form is closely related to the `global angular form' described in \cite[p.~70]{Bott-Tu}. Indeed, $\alpha$ restricts to each f\/iber to be $\beta d\theta$ which integrates out to $2\pi\beta$. It follows that $d\alpha = -2\pi\beta \phi^*e$, where $e \in \CI(Y;T^*Y)$ is the Euler class of $Y$ as a~submanifold of~$X$, and hence
\begin{gather*}
	j_0^*(\theta \wedge d\theta) = \begin{pmatrix} -\Id &0 & 0 \\ 0 & -\Id &0 \\ 0 & 0 & 0 \end{pmatrix} \alpha \wedge (-2\pi\beta \phi^*e).
\end{gather*}
Thus we f\/ind
\begin{gather}
	\int_{\pa M} Tp_1\big(\nabla, \nabla^{\pt}\big)
	= -\frac 1{8\pi^2} \int_{\pa M} \Tr j_0^*(\theta \wedge d\theta)\nonumber\\
\hphantom{\int_{\pa M} Tp_1\big(\nabla, \nabla^{\pt}\big)}{}
	= -\frac 1{8\pi^2} \int_{\pa M} (4\pi\beta \alpha\wedge\phi^*e)
	= -\beta^2 \int_Y e = -\beta^2 [Y]^2.\label{eq:Trans4D}
\end{gather}
This computation yields a formula for the index of the Dirac operator and, combined with results of Dai and Dai--Zhang, also a proof of the signature theorem of Atiyah--LeBrun.

\begin{Theorem}\label{thm:signature}
Let $X$ be an oriented four-dimensional manifold, $Y$ a smooth compact oriented embedded surface, and~$g$ an incomplete edge metric on $X\setminus Y$ with cone angle $2\pi \beta$ along~$Y$.
\begin{enumerate}\itemsep=0pt
\item[$1)$] If $X$ is spin and $\beta\in (0,1]$,
\begin{gather*}
	\ind\big(\eth\colon \mathcal D^+\lra L^2(M;\cS^-)\big) = -\frac1{24}\int_M p_1(M) + \frac1{24}\big(\beta^2 - 1\big)[Y]^2.
\end{gather*}

\item[$2)$] The signature of $X$ is given by $($see Atiyah--LeBrun~{\rm \cite{AL2013})}
\begin{gather*}
	\sign(X) = \frac1{12\pi^2}\int_M \big(|W_+|^2 - |W_-|^2\big) \, d\mu + \frac{1-\beta^2}3[Y]^2.
\end{gather*}
\end{enumerate}
\end{Theorem}

\begin{proof} 1) As mentioned above, the fact that the spin structure extends to all of $X$ and $\beta \in (0,1]$ implies that the generalized Witt assumption for $\eth$ is satisf\/ied. The degree four term of the $\hat A$ genus is $-p_1/24$, so applying our index formula~\eqref{eq:finalindexformula} and using the derivation of the local boundary term for $p_1$ in~\eqref{eq:Trans4D} gives
\begin{gather*}
\ind\big(\eth\colon \mathcal D^+\lra L^2(M;\cS^-)\big) = -\frac1{24}\int_M p_1(M) - \frac1{24}\beta^2[Y]^2 + \int_Y \hat A(Y)\lrpar{-\frac 12\hat\eta(\eth_Z)},
\end{gather*}
where the f\/inal term on the right is the limit $(1/2)\lim\limits_{\eps\to 0}\eta_{\eps}$ where $\epsilon_{\epsilon}$ is the eta-invariants induced on the boundary of $M_\epsilon$ as $\epsilon \to 0$. Thus we claim (and it remains to prove) that the adiabatic limit of the eta-invariant for the spin Dirac operator is
\begin{gather}\label{eq:explicit-adiabatic-eta}
	\lim_{\eps\to 0}\tfrac12\eta_{\eps} = \frac1{24}[Y]^2,
\end{gather}
i.e., the limit of the eta-invariants is the \textit{opposite} of the local boundary term when $\beta = 1$, which indeed it should be since in that case the metric is smooth across $x = 0$.

Although other derivations of the adiabatic eta invariant exist \cite{DZ1995}, we prefer to give on here which we f\/ind intuitive and which f\/its nicely with arguments above. To this end, we consider~$N$, a disc bundle over a smooth manifold~$Y$, and we assume~$N$ is spin. We will show below that~$N$ admits a positive scalar curvature metric. Thus, given a spin structure and metric, the index of~$\eth$ vanishes on~$N$. If we furthermore note that~$N$ is dif\/feomorphic to $[0, 1)_x \times X$ where~$X$ is a~circle bundle over~$Y$, and let $N^\epsilon = [0, \epsilon)_x \times X$, we may consider metrics
\begin{gather}\label{eq:disc-bundle-metric}
g = dx^2 + f^2(x) k + h,
\end{gather}
where $h$ is the pullback of a metric on $Y$, $k \in \operatorname{Sym}^{0,2}(N^\epsilon)$ $x$ and $dx$-independent and restricts to a Riemannian metric on the f\/ibers of $X$. We assume $f$ is smooth across $x = 0$ with $f(x) = x + O(x^2)$ which implies that $g$ is \textit{smooth} on~$N$. Using the computation of the connection above, with respect to the orthonormal basis, $X_{i}$, $\frac{1}{f} U$, $\p_x$ the connection one form of~$g$ is
\begin{gather} \label{eq:curvoneform}
 \omega = \lp\begin{array}{c|c|c} \wt{\omega}_{\wt{h}} -
 f^{2}\frac{1}{2}g^{\p}\mathcal{R} & - f g^{\p}\big(
 \widehat{\Rmnum{2}} + \frac{1}{2}\widehat{\mathcal{R}} \big) & 0 \\
 \hline
 f g^{\p}\big(\widehat{\Rmnum{2}} +
 \frac{1}{2}\widehat{\mathcal{R}}\big) & 0 & f' U^{\sharp} \\ \hline
 0 & - f' U^{\sharp} & 0
 \end{array}\rp,
\end{gather}
where $g^{\p} = k + h$ is the metric on the circle bundle $X$. We will take
\begin{gather*}
f(x) = f_\epsilon(x) = x \chi_\epsilon(x),
\end{gather*}
where $\chi$ is a smooth positive function that is monotone decreasing with $\chi(x)= 1 $ for $x \le 1/3$ and $\chi(x) = \beta$ for $x \ge 2/3$. Then $f = f' = f'' = O(1/\epsilon)$, and using $\Omega = d\omega + \omega \wedge \omega$, we see that
\begin{gather*}
\hat{A}_g = F d{\rm Vol}_g,
\end{gather*}
where $F$ is a function that is $O(1/\epsilon)$. Since ${\rm Vol}(N_\epsilon) = O(\epsilon^2)$,
\begin{gather*}
\int_{N_\epsilon} \hat{A}_\epsilon = - \frac{1}{24} \int_{N_\epsilon} p_1 \to 0 \qquad \text{as}\quad \epsilon \to 0.
\end{gather*}
Since the index of the Dirac operator vanishes on $N^\epsilon$, applying the APS formula gives
\begin{gather*} 
 0 = 	- \frac{1}{24} \int_{N^{\epsilon}} p_1 + \int_{\pa N^{\epsilon}} Tp_1\big(\nabla, \nabla^{\pt}\big) -\tfrac12\eta\big(\pa N^{\epsilon}\big),
\end{gather*}
where $\nabla^{\pt}$ is as in~\eqref{eq:TwoConn}, and thus the limit of the trangression forms is exactly as computed above. Thus by taking the $\epsilon \to 0$ limit we obtain~\eqref{eq:explicit-adiabatic-eta}.

To prove part 1) it remains to prove the existence of a positive scalar curvature metric on $N$. To this end we take the metric $g$ as in~\eqref{eq:disc-bundle-metric} on $N^\epsilon$ now with
\begin{gather*} 
 f(x) = f_\delta(x) = \delta \sin(x/\delta).
\end{gather*}
Note that $f = O(\epsilon)$, $f' = O(\epsilon/\delta)$. Then curvature equals
\begin{gather*}
 \Omega = d\omega + \omega \wedge \omega = \lp\begin{array}{c|c|c} \wt{\Omega}_{\wt{h}}
 & 0 & 0 \\
 \hline
 0 & 0 & f'' dx \wedge U^{\sharp} \\ \hline
 0 & - f'' dx \wedge U^{\sharp} & 0
 \end{array}\rp + O(\epsilon) + O(\epsilon/\delta).
\end{gather*}
Denoting our orthonormal basis by $e_i$, $i = 1, \dots, n$ and taking traces gives
\begin{gather*}
\operatorname{scal}_g = \delta^{ik} \delta^{jl}\Omega_{ij}(e_k, e_l) = \operatorname{scal}_h + \frac{2}{\delta^2} + O(\epsilon/\delta),
\end{gather*}
and thus taking $\epsilon/\delta = 1$ and $\delta$ small gives a~positive scalar curvature metric.

2) Since $X$ is a smooth manifold we can use Novikov additivity of the signature to decompose the signature as
\begin{gather*}
	\sign(X) = \sign(X \setminus M_{\eps}) + \sign(M_{\eps}).
\end{gather*}
Identifying $X\setminus M_{\eps}$ with a disk bundle over $Y$ we have from \cite[p.~314]{Dai:adiabatic} that
\begin{gather*}
	\sign(X \setminus M_{\eps}) = \sign \lrpar{ \int_Y e },
\end{gather*}
i.e., the signature is the sign of the self-intersection number of $Y$ in $X$. In fact this is a simple exercise using the Thom isomorphism theorem.

The Atiyah--Patodi--Singer index theorem for the signature of $M_{\eps}$ yields
\begin{gather*}
	\sign(M_{\eps}) = \frac13\int_{M_{\eps}} p_1(\nabla) + \frac13\int_{\pa M_{\eps}} Tp_1(\nabla, \nabla^{\pt}) - \eta^{\ev}_{\eps},
\end{gather*}
where $\eta^{\ev}_{\eps}$ is the eta-invariant of the boundary signature operator restricted to forms of even degree. As $\eps\to0$, the eta invariant is undergoing adiabatic degeneration and its limit is computed in \cite[Theorem~3.2]{DZ1995},
\begin{gather*}
	\lim_{\eps\to 0} = -\int_Y L(TY)\big(\coth e - e^{-1}\big) + \sign \lrpar{ B_e },
\end{gather*}
where $B_e$ is the bilinear form on $H^0(Y)$ given by $H^0(Y) \ni c,c' \mapsto cc'\ang{e,Y} \in \bbR$, i.e., it is again the sign of the self-intersection of~$Y$. (In comparing with~\cite{DZ1995} note that the orientation of $\pa M_{\eps}$ is the opposite of the orientation of the spherical normal bundle of~$Y$ in~$X$, and so $\sign(B_e) = -\sign(X\setminus M_{\eps})$.) The only term in $L(TY)(\coth e - e^{-1})$ of degree two is $\tfrac13e$, and hence
\begin{gather*}
	\sign(X) = \frac13\int_X p_1 + \frac13[Y]^2 + \frac13 \big({-}\beta^2 [Y]^2\big)
\end{gather*}
as required. (Note that we could also argue as in the Dirac case to compute the limit of the eta invariants.)
\end{proof}

\section{Positive scalar curvature metrics}\label{sec:posscal}

In this short section, we prove Theorem \ref{thm:PosScal} following~\cite{Ch1985}. We recall the statement of the theorem for the convenience of the reader:

\medskip

\noindent \textbf{Theorem.} {\it Let $(M,g)$ be a spin space with an incomplete edge metric. The
`geometric Witt assumption'~\eqref{eq:GeoWittAss} holds if either:
\begin{enumerate}\itemsep=0pt
\item[$1.$] $\dim Z \ge 2$ and the scalar curvature of $g$ is non-negative in a neighborhood of
 $\p M$.
\item[$2.$] $\dim Z = 1$, the spin structure on $M$ is the lift of a spin structure on $X$, and the cone angle satisfies $2 \pi \beta \le 2
 \pi$.
\end{enumerate}

If the geometric Witt assumption holds and in addition the scalar curvature of $g$ is non-negative on all of $M$, and positive somewhere, then $\ind(\eth)=0$.}

\begin{proof}
1) Taking traces in \eqref{eq:AsympSplittingCurv}, the scalar curvature $R_{g}$ satisf\/ies
\begin{gather*}
	R_{g} = R_{\rm cone} + \cO(1),
\end{gather*}
where $R_{\rm cone}$ is the scalar curvature of the cone with metric $dx^{2} + x^{2} g_{N/Y} \rvert_{\ang{\pa_x} \oplus \tfrac1xTN/Y}$, as in
\eqref{eq:SecondSplittingIE}. On the other hand, by \cite[Section~4]{Ch1985}, the scalar curvature of an exact cone $C(Z)$ is equal to $x^{-2}(R_{Z} - \dim(Z)(\dim(Z) - 1))$, where $R_{Z}$ is the scalar curvature of $Z$. Thus $R_{g} \ge 0$ implies that $R_{Z} \ge \dim(Z)(\dim (Z) - 1)$, which by \cite[Lemma~3.5]{Ch1985} shows that Assumption~\ref{thm:assumption} holds.

2) In the circle f\/iber case, there exist local trivializations of the boundary f\/ibration such that $\eth^Z_y = \frac{1}{\beta^2}
\eth_\theta$ where $\eth_\theta$ is the Dirac operator on $\mathbb{S}^1$ for the spin structure that bounds a disk. Since $\spec(\eth_\theta) \cap (-1/2, 1/2) = \varnothing$, the assumption $\beta \le 1$ implies that also $\spec(\eth^Z_y ) \cap (-1/2, 1/2) = \varnothing$ in this case.

Now assume that the geometric Witt condition is satisf\/ied. By Theorem \ref{thm:essential}, $\eth$ is essentially self-adjoint. That is, the graph closure of $\eth$ on~$C^{\infty}_c(M)$ is self-adjoint, with domain $\mathcal{D}$ from Theorem~\ref{thm:essential}, and furthermore by the Main Theorem its index satisf\/ies~\eqref{eq:finalindexformula}.

From the Lichnerowicz formula \cite{BGV2004},
\begin{gather*}
	\eth^{*} \eth = \nabla^{*}\nabla + R/4,
\end{gather*}
where $R$ is the scalar curvature. Thus, for every $\phi \in C^{\infty}_c(M)$, $\norm[L^{2}]{\eth \phi} = \norm[L^{2}]{\nabla \phi} + \la R \phi, \phi \ra_{L^{2}}$. We conclude that for all $\phi \in C^{\infty}_c(M)$,
\begin{gather}\label{eq:onccinfty}
	\norm[L^{2}]{\eth \phi} \ge \norm[L^{2}]{\eth \phi} - \la R \phi, \phi \ra_{L^{2}} \ge \norm[L^{2}]{\nabla \phi} \ge 0.
\end{gather}
This implies in particular that $\mathcal{D}_{\min}(\eth) = \mathcal{D} \subset \mathcal{D}_{\min}(\nabla)$, where we recall that $\mathcal{D}_{\min}(P)$ refers to the graph closure of the operator $P$ with domain~$C^{\infty}_c(M)$. We claim that the index of the operator
\begin{gather*}
	\eth \colon \ \mathcal{D}^{+} \lra L^{2}(M; \SPB^{-})
\end{gather*}
vanishes, so by formula~\eqref{eq:finalindexformula}, Theorem \ref{thm:PosScal}(b) holds. In fact, the kernel of $\eth$ on $\mathcal{D}$ consists only of the zero vector, since if $\phi \in \mathcal{D}$ has $\eth \phi = 0$, then since~\eqref{eq:onccinfty} holds on $\mathcal{D}$, $\nabla \phi = 0$ also. By the Lichnerowicz formula again, $R \phi = 0$, but since by assumption $R$ is not identically zero, $\phi$ must vanish somewhere and by
virtue of its being parallel, $\phi \equiv 0$. \end{proof}

\appendix
\section{Appendix}\label{sec:Bessel}

In this appendix we prove Claim~\ref{thm:normalwithinhalf} by using standard bounds on modif\/ied Bessel functions to prove the sup norm bound~\eqref{eq:matrixbound}: for each $\mu$ with $\absv{\mu} > 1/2$, and all $\absv{\eta}$,
\begin{gather*}
	\big\|\mathcal{N}_{\mu, \absv{\eta}} - \mathcal{N}^{\rm APS}_{\mu, \absv{\eta}}\big\|
	< 1 - \delta.
\end{gather*}
Among references for modif\/ied Bessel functions we recall \cite{A1974, Bar2010, Bar2013, PBG2007}.

To begin with, using the Wronskian equation \eqref{eq:modifiedBesselrecurrencewronsk}, note that
\begin{gather*}
	\Tr \mathcal{N}_{\mu, z} = \Tr \mathcal{N}^{\rm APS}_{\mu, z} = 1.
\end{gather*}
Thus the dif\/ference $\mathcal{N}_{\mu, z} -\mathcal{N}^{\rm APS}_{\mu, z}$ has two equal eigenvalues and hence its norm is the square root of the determinant. We now assume that $\mu \ge 1/2$, since the $\mu \le -1/2$ case is treated the same way. Using \eqref{eq:modifiedBesselrecurrencewronsk} again, we see that
\begin{gather}
\det\big(\mathcal{N}_{\mu, z} -\mathcal{N}^{\rm APS}_{\mu, z}\big)= - \frac 12 + \frac 12 \frac{z}{(\mu^{2} + z^{2})^{1/2}} \big(\mu (I_{\mu-1/2}K_{\mu + 1/2} - I_{\mu+1/2}K_{\mu-1/2} ) \nonumber \\
\hphantom{\det\big(\mathcal{N}_{\mu, z} -\mathcal{N}^{\rm APS}_{\mu, z}\big)}{} + z(I_{\mu+1/2}K_{\mu + 1/2} + I_{\mu-1/2}K_{\mu-1/2} ) \big),\label{eq:determinant}
\end{gather}
and we want to show that for some $\delta > 0$ independent of $\mu \ge 1/2$ and $z \ge 0$,
\begin{gather}\label{eq:Besseltoshow}
	-1 + \delta \le \det\big(\mathcal{N}_{\mu, z} -\mathcal{N}^{\rm APS}_{\mu, z}\big) \le 1 - \delta.
\end{gather}

To begin with, we prove that
\begin{gather}\label{eq:betweenzeroandhalf}
0 \le z I_{\nu}(z) K_{\nu}(z) \le 1/2 \qquad \text{for} \quad \nu \ge 1/2, z \ge 0.
\end{gather}
In fact, we claim that for $\nu \ge 1/2$, $zK_{\nu}(z)I_{\nu}(z) $ is monotone. To see that this holds, dif\/ferentiate
\begin{gather*}
(zK_{\nu}(z)I_{\nu}(z))' = K_{\nu}I_{\nu} + z(K'_{\nu}I_{\nu} +
K_{\nu}I'_{\nu}) = K_{\nu}I_{\nu}\left( 1 + \frac{zK'_{\nu}(z)}{K_{\nu}(z)} + \frac{zI'_{\nu}(z)}{I_{\nu}(z)}\right).
\end{gather*}
Thus we want to show that $\frac{zK'_{\nu}(z)}{K_{\nu}(z)} + \frac{zI'_{\nu}(z)}{I_{\nu}(z)}) \ge -1$. Using \cite[equation~(5.1)]{Bar2013}, for $\nu \ge 1/2$
\begin{gather*}
	\lp\frac{zK'_{\nu}(z)}{K_{\nu}(z)}\rp' + \lp \frac{zI'_{\nu}(z)}{I_{\nu}(z)}\rp' \le 0,
\end{gather*}
so the quantity $\frac{zK'_{\nu}(z)}{K_{\nu}(z)} + \frac{zI'_{\nu}(z)}{I_{\nu}(z)}$ is monotone decreasing. In fact, we claim that
\begin{gather*}
	\frac{zK'_{\nu}(z)}{K_{\nu}(z)} +
	\frac{zI'_{\nu}(z)}{I_{\nu}(z)}
	\to\begin{cases}
0 &\text{as} \ z \to 0, \\
 -1 &\text{as} \ z \to \infty.
	\end{cases}
\end{gather*}
The limit as $z \to \infty$ can be seen using the large argument asymptotic formulas from \cite[Section~9.7]{AS1964}, while the limit as $z \to 0$ follows from the recurrence relations~\eqref{eq:modifiedBesselrecurrencewronsk} and the small argument asymptotics in \cite[Section~9.6]{AS1964}.
Thus $z K_{\nu}(z) I_{\nu}(z)$ is monotone on the region under consideration. Using the asymptotic formulas again shows that
\begin{gather*}
	zK_{\nu}(z)I_{\nu}(z)
\to	\begin{cases}
	0 &\text{as} \ z \to 0, \\
	 1/2 &\text{as} \ z \to \infty,
	\end{cases}
\end{gather*}
so \eqref{eq:betweenzeroandhalf} holds.

We can now show the upper bound in \eqref{eq:Besseltoshow}. Using the Wronskian relation in \eqref{eq:modifiedBesselrecurrencewronsk},
we write
\begin{gather}
\det\big(\mathcal{N}_{\mu, z} -\mathcal{N}^{\rm APS}_{\mu, z}\big)
 = - \frac 12 + \frac 12 \frac{\mu}{(\mu^{2} + z^{2})^{1/2}} + \frac 12 \frac{z}{(\mu^{2} + z^{2})^{1/2}}
	\big({-} 2\mu I_{\mu+1/2}K_{\mu-1/2} \nonumber\\
\hphantom{\det\big(\mathcal{N}_{\mu, z} -\mathcal{N}^{\rm APS}_{\mu, z}\big)=}{}
+ z(I_{\mu+1/2}K_{\mu + 1/2} +	I_{\mu-1/2}K_{\mu-1/2} ) \big) \nonumber\\
\hphantom{\det\big(\mathcal{N}_{\mu, z} -\mathcal{N}^{\rm APS}_{\mu, z}\big)}{}
\le \frac 12 \frac{z}{(\mu^{2} + z^{2})^{1/2}} \big(z(I_{\mu+1/2}K_{\mu + 1/2} + I_{\mu-1/2}K_{\mu-1/2} ) \big).\label{eq:determinanttwo}
\end{gather}
Now, if $\mu \ge 1$, by \eqref{eq:betweenzeroandhalf}, the right hand side in the f\/inal inequality is bounded by $1/2$, establishing the upper
bound in \eqref{eq:determinant} in this case (with $\delta = 1/2$). If $\mu \in [1/2, 1]$, we use the following inequalities of Barciz \cite[equations~(2.3),~(2.4)]{Bar2013}
\begin{gather*}
	\frac{zI_{\nu}'(z)}{I_{\nu}(z)} < \sqrt{z^{2} + \nu^{2}} \qquad \text{and} \qquad \frac{zK_{\nu}'(z)}{K_{\nu}(z)} < - \sqrt{z^{2} + \nu^{2}} ,
\end{gather*}
for $\nu \ge 0$, $z \ge 0$. Using these inequalities and the recurrence relation~\eqref{eq:modifiedBesselrecurrencewronsk} gives
\begin{gather*}
\frac{I_{\mu - 1/2}}{I_{\mu + 1/2}}	< \frac{\sqrt{z^{2} + (\mu - 1/2)^{2}} + \mu - 1/2}{z}, \qquad \frac{K_{\mu - 1/2}}{K_{\mu + 1/2}} < \frac{z}{\sqrt{z^{2} + (\mu + 1/2)^{2}} + \mu + 1/2},
\end{gather*}
so continuing the inequality \eqref{eq:determinanttwo} gives
\begin{gather}
	\det\big(\mathcal{N}_{\mu, z} -\mathcal{N}^{\rm APS}_{\mu, z}\big)
\le \frac 12 \frac{z}{(\mu^{2} + z^{2})^{1/2}} \lp zI_{\mu+1/2}K_{\mu + 1/2} \rp \nonumber\\
\hphantom{\det\big(\mathcal{N}_{\mu, z} -\mathcal{N}^{\rm APS}_{\mu, z}\big)\le }{}
 \times\lp 1 + \frac{\sqrt{z^{2} + (\mu + 1/2)^{2}} + \mu + 1/2}{\sqrt{z^{2} + (\mu - 1/2)^{2}} + \mu - 1/2} \rp.\label{eq:determinantthree}
\end{gather}
One checks that for $1/2 \le \mu$, the fraction in the second line is monotone decreasing in $z$, and thus by \eqref{eq:betweenzeroandhalf}, for $z \ge 1$ the determinant is bounded by
\begin{gather}\label{eq:100}
	\frac 14 \lp 1 + \frac{\sqrt{1 + (\mu + 1/2)^{2}} + \mu + 1/2}{\sqrt{1 + (\mu - 1/2)^{2}}+ \mu - 1/2} \rp \le \frac 14 \big(1 + \big(1 + \sqrt{2}\big)\big) \le 1 - \delta,
\end{gather}
where the middle bound is obtained by checking that the fraction on the left is monotone decreasing in $\mu$ for $\mu \ge 1/2$ and equal to $1 + \sqrt{2}$ at $\mu = 1/2$. Thus, we have established the upper bound in~\eqref{eq:Besseltoshow} in the region $z \ge 1$. For $z \le
1$, rewrite the bound in~\eqref{eq:determinantthree} as
\begin{gather*}
	\frac 12 \frac{z}{(\mu^{2} + z^{2})^{1/2}} \lp
	I_{\mu+1/2}K_{\mu + 1/2} \rp z \lp 1 + \frac{\sqrt{z^{2} +
	(\mu + 1/2)^{2}} + \mu +
	1/2}{\sqrt{z^{2} + (\mu - 1/2)^{2}} + \mu - 1/2} \rp.
\end{gather*}
For $\mu \ge 1/2$, by \cite{PBG2007}, the function $I_{\mu + 1/2}(z)K_{\mu + 1/2}(z)$ is monotone decreasing, and by the asymptotic formulas it is goes to $1/2$ as $z \to 0$. Thus in $0 \le z \le 1$ the determinant is bounded about by
\begin{gather*}
\frac 14 z \lp 1 + \frac{\sqrt{z^{2} + (\mu + 1/2)^{2}} + \mu + 1/2}{\sqrt{z^{2} + (\mu - 1/2)^{2}} + \mu - 1/2} \rp.
\end{gather*}
This function is monotone increasing in $z$ for $\mu \in [1/2, 1]$, so the max is obtained at $z = 1$, i.e., it is bounded by the left hand side of~\eqref{eq:100}, in particular by $1 - \delta$ for the same $\delta$. This establishes the upper bound in~\eqref{eq:Besseltoshow}.

Finally we establish the lower bound. First, we rewrite the determinant again, this time using the Wronskian relation in the opposite direction to obtain
\begin{gather}
\det\big(\mathcal{N}_{\mu, z} -\mathcal{N}^{\rm APS}_{\mu,z}\big) = - \frac 12 - \frac 12 \frac{\mu}{(\mu^{2} +
	z^{2})^{1/2}} + \frac 12 \frac{z}{(\mu^{2} + z^{2})^{1/2}}
	\big( 2\mu I_{\mu-1/2}K_{\mu+1/2} \nonumber\\
\hphantom{\det\big(\mathcal{N}_{\mu, z} -\mathcal{N}^{\rm APS}_{\mu,z}\big) =}{} + z(I_{\mu+1/2}K_{\mu + 1/2} +
	I_{\mu-1/2}K_{\mu-1/2} ) \big).\label{eq:determinantfour}
\end{gather}
Now, recalling that $zI_{\mu+1/2}(1)K_{\mu + 1/2}(1)$ is monotone increasing, using the asymptotic formulas \cite[equations~(9.7.7), (9.7.8)]{AS1964} we see that
\begin{gather*}
	I_{\mu+1/2}(1)K_{\mu + 1/2}(1) \to \frac{1}{2(\mu + 1/2)}
\end{gather*}
as $\mu \to \infty$, we use the inequality \cite[equation~(11)]{A1974}, namely
\begin{gather*}
I_{\mu - 1/2}(z) \ge \frac{\mu - 1/2 + \big(z^{2} + (\mu + 3/2)^{2}\big)^{1/2}}{z}I_{\mu + 1/2}(z).
\end{gather*}
On the region $z \in [0, 1]$, $\mu - 1/2 + (z^{2} + (\mu + 3/2)^{2})^{1/2} \ge \delta_{0} > 0$. Dropping the terms with equal order in~\eqref{eq:determinantfour}
then gives
\begin{gather*}
\det\big(\mathcal{N}_{\mu, z} -\mathcal{N}^{\rm APS}_{\mu, z}\big)
> -1 + \frac 12 \frac{z}{(\mu^{2} + z^{2})^{1/2}}	2\mu I_{\mu-1/2}K_{\mu+1/2} \\
\hphantom{\det\big(\mathcal{N}_{\mu, z} -\mathcal{N}^{\rm APS}_{\mu, z}\big)}{}
\ge -1 + \frac 12 \frac{2\mu}{(\mu^{2} + z^{2})^{1/2}}
	I_{\mu+1/2}K_{\mu+1/2} \big(\mu - 1/2 + \big(z^{2} + (\mu + 3/2)^{2}\big)^{1/2}\big) \\
\hphantom{\det\big(\mathcal{N}_{\mu, z} -\mathcal{N}^{\rm APS}_{\mu, z}\big)}{}
	\ge -1 + \delta_{0} \frac 12 \frac{\mu - 1/2 + \big(z^{2} + (\mu +	3/2)^{2}\big)^{1/2}}{(\mu^{2} +	z^{2})^{1/2}}
	\ge -1 + \delta.
\end{gather*}
This completes the proof of \eqref{eq:Besseltoshow}.

\subsection*{Acknowledgements}
P.A.~was supported by NSF grant DMS-1104533 and Simons Foundation grant \#317883. The authors are happy to thank Rafe Mazzeo and Richard Melrose for many useful and interesting discussions. They are also grateful to the comments of the anonymous referees, particularly their suggestion of Remark~\ref{thm:edge-theorem}.

\LastPageEnding

\end{document}